\documentclass[English,11pt]{amsart}

\usepackage{amssymb}
\usepackage{amsthm, amsfonts, thmtools, thm-restate}
\usepackage{dsfont}
\usepackage{mathtools} 
\usepackage{hyperref,cleveref}

\usepackage{enumitem}

\usepackage[backend=biber, style=ieee-alphabetic, maxnames=5, maxalphanames=5, giveninits=true, sorting=nyt, isbn=false, doi=false, url=false, eprint=false]{biblatex} 

\DeclareNameAlias{author}{first-last}
\AtEveryBibitem{\clearlist{language}}
\AtEveryBibitem{\clearlist{month}}
\newtheorem{thm}{Theorem}[section]
\newtheorem{prop}[thm]{Proposition}
\newtheorem{lem}[thm]{Lemma}
\newtheorem{cor}[thm]{Corollary}
\newtheorem{definition}[thm]{Definition}

\newtheorem{rmk}[thm]{Remark}
\usepackage[all]{xy} 

\usepackage{lscape} 

\usepackage{graphicx}
\usepackage[left=3cm,right=3cm,top=4cm,bottom=3cm]{geometry} 
\usepackage{tabularx}
\usepackage{mathrsfs} 
\usepackage{frcursive}
\usepackage[usenames, dvipsnames]{xcolor}

\usepackage{caption}
\usepackage{subcaption}

\usepackage{pgf,tikz}
\usetikzlibrary{arrows}
\usepackage{xcolor} 
\usepackage{cleveref}
\usepackage{color}

\newcommand{\CC}{\mathbb{C}}
\newcommand{\RR}{\mathbb{R}}
\newcommand{\NN}{\mathbb{N}}

\renewcommand{\Im}{\mathrm{Im}}

\title{Deviations for generalized tiling billiards in cyclic polygons}
\author{Magali Jay }

\begin{document}

\maketitle

\textbf{Abstract.}
This work continues the study of tiling billiards, a class of dynamical system introduced by Davis et al. in 2018. We develop the study of generalized tiling billiards in a cyclic polygon. This work shows that the behavior of generalized tiling billiards in cyclic $N$-gons with $N\geqslant 5$ is considerably different from that of triangular and quadrilateral tiling billiards studied before. Indeed,
we exhibit an open set of generalized tiling billiard trajectories deviating sublinearly from their asymptotic direction, whereas for $N=3$ or $4$ almost every trajectory stays at a bounded distance from a line. Moreover, we establish the rate of deviations both in the generic case and in some non generic cases.

\section{Introduction}

\subsection{A quick bibliographical orverview of tiling billiards}
Davis, DiPietro, Rustad and St-Laurent \cite{DDiPRStL} defined and studied a class of dynamical systems that they called \emph{tiling billiards}. In a tiling billiard, the trajectory is refracted across the boundary of a plane tiling (the angle of incidence and the angle of refraction are opposite), see Figure~\ref{Fig:ex_polTB}. 
Nogueira had already suggested to study such billiards at the end of \cite{Nogueira}, in relation with interval exchange transformations with flips (IETFs). 
\vspace*{-0.5cm}
\begin{figure}[h!]
\centering
\begin{subfigure}[b]{0.27\textwidth}
         \centering
\definecolor{zzttqq}{rgb}{0.6,0.2,0}
\begin{tikzpicture}[line cap=round,line join=round,>=triangle 45,x=1.0cm,y=1.0cm]
\clip(4.26,-1.24) rectangle (9.85,1.49);
\draw[color=zzttqq,fill=zzttqq,fill opacity=0.1] (6,0) -- (7,0) -- (6.5,0.87) -- cycle;
\draw[color=zzttqq,fill=zzttqq,fill opacity=0.1] (7.5,0.87) -- (7,0) -- (6.5,0.87) -- cycle;
\draw[color=zzttqq,fill=zzttqq,fill opacity=0.1] (6,0) -- (7,0) -- (6.5,-0.87) -- cycle;
\draw[color=zzttqq,fill=zzttqq,fill opacity=0.1] (7.5,0.87) -- (7,0) -- (8,0) -- cycle;
\draw[color=zzttqq,fill=zzttqq,fill opacity=0.1] (7.5,-0.87) -- (7,0) -- (6.5,-0.87) -- cycle;
\draw[color=zzttqq,fill=zzttqq,fill opacity=0.1] (7.5,-0.87) -- (7,0) -- (8,0) -- cycle;
\draw (6.81,0)-- (6.74,0.45);
\draw (7.09,0.16)-- (6.74,0.45);
\draw (6.81,0)-- (6.74,-0.45);
\draw (7.09,0.16)-- (7.52,0);
\draw (7.09,-0.16)-- (6.74,-0.45);
\draw (7.09,-0.16)-- (7.52,0);
\end{tikzpicture}
\end{subfigure}
\hfill
\begin{subfigure}[b]{0.27\textwidth}
         \centering
\definecolor{zzttqq}{rgb}{0.6,0.2,0}
\begin{tikzpicture}[line cap=round,line join=round,>=triangle 45,x=1.0cm,y=1.0cm]
\fill[color=zzttqq,fill=zzttqq,fill opacity=0.1] (0,0) -- (1,0) -- (1,1) -- (0,1) -- cycle;
\fill[color=zzttqq,fill=zzttqq,fill opacity=0.1] (2,0) -- (1,0) -- (1,1) -- (2,1) -- cycle;
\fill[color=zzttqq,fill=zzttqq,fill opacity=0.1] (0,2) -- (1,2) -- (1,1) -- (0,1) -- cycle;
\fill[color=zzttqq,fill=zzttqq,fill opacity=0.1] (2,2) -- (1,2) -- (1,1) -- (2,1) -- cycle;
\fill[color=zzttqq,fill=zzttqq,fill opacity=0.1] (4,2) -- (3,2) -- (3,1) -- (4,1) -- cycle;
\fill[color=zzttqq,fill=zzttqq,fill opacity=0.1] (2,2) -- (3,2) -- (3,1) -- (2,1) -- cycle;
\fill[color=zzttqq,fill=zzttqq,fill opacity=0.1] (2,0) -- (3,0) -- (3,1) -- (2,1) -- cycle;
\fill[color=zzttqq,fill=zzttqq,fill opacity=0.1] (4,0) -- (3,0) -- (3,1) -- (4,1) -- cycle;
\draw [color=zzttqq] (0,0)-- (4,0);
\draw [color=zzttqq] (0,1)-- (4,1);
\draw [color=zzttqq] (0,2)-- (4,2);
\draw [color=zzttqq] (0,2)-- (0,0);
\draw [color=zzttqq] (1,2)-- (1,0);
\draw [color=zzttqq] (2,2)-- (2,0);
\draw [color=zzttqq] (3,2)-- (3,0);
\draw [color=zzttqq] (4,2)-- (4,0);
\draw (0.28,1)-- (1,1.28);
\draw (1.72,1)-- (1,1.28);
\draw (1.72,1)-- (1,0.72);
\draw (0.28,1)-- (1,0.72);
\draw (0,1.34)-- (1,1.73);
\draw (2,1.34)-- (1,1.73);
\draw (2,1.34)-- (3,1.73);
\draw (4,1.34)-- (3,1.73);
\end{tikzpicture}
\end{subfigure}
\hfill
\begin{subfigure}[b]{0.27\textwidth}
         \centering
\definecolor{zzttqq}{rgb}{0.6,0.2,0}
\begin{tikzpicture}[line cap=round,line join=round,>=triangle 45,x=0.5cm,y=0.5cm]
\draw[color=zzttqq,fill=zzttqq,fill opacity=0.1] (14.72,-6.58) -- (15.58,-6.08) -- (15.58,-5.08) -- (14.72,-4.58) -- (13.85,-5.08) -- (13.85,-6.08) -- cycle;
\draw[color=zzttqq,fill=zzttqq,fill opacity=0.1] (16.45,-6.58) -- (15.58,-6.08) -- (15.58,-5.08) -- (16.45,-4.58) -- (17.32,-5.08) -- (17.32,-6.08) -- cycle;
\draw[color=zzttqq,fill=zzttqq,fill opacity=0.1] (16.45,-3.58) -- (16.45,-4.58) -- (15.58,-5.08) -- (14.72,-4.58) -- (14.72,-3.58) -- (15.58,-3.08) -- cycle;
\draw[color=zzttqq,fill=zzttqq,fill opacity=0.1] (12.99,-3.58) -- (13.85,-3.08) -- (14.72,-3.58) -- (14.72,-4.58) -- (13.85,-5.08) -- (12.99,-4.58) -- cycle;
\draw[color=zzttqq,fill=zzttqq,fill opacity=0.1] (12.99,-6.58) -- (12.12,-6.08) -- (12.12,-5.08) -- (12.99,-4.58) -- (13.85,-5.08) -- (13.85,-6.08) -- cycle;
\draw[color=zzttqq,fill=zzttqq,fill opacity=0.1] (14.72,-6.58) -- (15.58,-6.08) -- (16.45,-6.58) -- (16.45,-7.58) -- (15.58,-8.08) -- (14.72,-7.58) -- cycle;
\draw[color=zzttqq,fill=zzttqq,fill opacity=0.1] (14.72,-6.58) -- (14.72,-7.58) -- (13.85,-8.08) -- (12.99,-7.58) -- (12.99,-6.58) -- (13.85,-6.08) -- cycle;
\draw[color=zzttqq,fill=zzttqq,fill opacity=0.1] (18.18,-6.58) -- (17.32,-6.08) -- (16.45,-6.58) -- (16.45,-7.58) -- (17.32,-8.08) -- (18.18,-7.58) -- cycle;
\draw[color=zzttqq,fill=zzttqq,fill opacity=0.1] (19.05,-5.08) -- (19.05,-6.08) -- (18.18,-6.58) -- (17.32,-6.08) -- (17.32,-5.08) -- (18.18,-4.58) -- cycle;
\draw[color=zzttqq,fill=zzttqq,fill opacity=0.1] (16.45,-3.58) -- (17.32,-3.08) -- (18.18,-3.58) -- (18.18,-4.58) -- (17.32,-5.08) -- (16.45,-4.58) -- cycle;
\draw (14.51,-4.69)-- (13.85,-5.32);
\draw (13.19,-4.69)-- (13.85,-5.32);
\draw (13.19,-4.69)-- (14.07,-4.95);
\draw (13.85,-5.84)-- (14.07,-4.95);
\draw (13.85,-5.84)-- (13.64,-4.95);
\draw (14.51,-4.69)-- (13.64,-4.95);
\draw (14.9,-4.68)-- (14.34,-6.36);
\draw (13.16,-7.68)-- (14.34,-6.36);
\draw (14.9,-4.68)-- (16.07,-3.36);
\draw (15.37,-4.95)-- (14.95,-6.44);
\draw (16.45,-6.82)-- (14.95,-6.44);
\draw (15.37,-4.95)-- (16.45,-3.85);
\draw (17.52,-4.95)-- (16.45,-3.85);
\draw (17.52,-4.95)-- (17.95,-6.44);
\draw (16.45,-6.82)-- (17.95,-6.44);
\end{tikzpicture}
\end{subfigure}
\caption{Trajectories of regular tiling billiards}
\label{fig:reg_TB}
\end{figure}
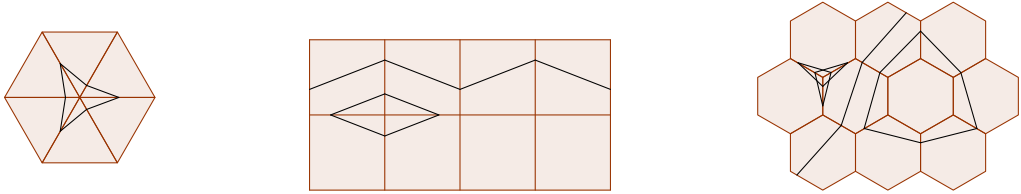

\begin{figure}[h!]
\definecolor{zzttqq}{rgb}{1,0.4,0} 
\begin{tikzpicture}[line cap=round,line join=round,>=triangle 45,x=.5cm,y=.5cm]
\fill[color=zzttqq,fill=zzttqq,fill opacity=0.05] (0.94,0.36) -- (3.44,-0.14) -- (3.64,1.36) -- (2.1,2.78) -- (0.24,1.46) -- cycle;
\draw[color=zzttqq] (0.94,0.36) -- (3.44,-0.14) -- (3.64,1.36) -- (2.1,2.78) -- (0.24,1.46) -- cycle;
\fill[color=zzttqq,fill=zzttqq,fill opacity=0.05] (3.64,1.36) -- (3.44,-0.14) -- (4.64,-1.06) -- (6.04,-0.49) -- (6.24,1.01) -- (5.04,1.94) -- cycle;
\draw[color=zzttqq] (3.64,1.36) -- (3.44,-0.14) -- (4.64,-1.06) -- (6.04,-0.49) -- (6.24,1.01) -- (5.04,1.94) -- cycle;
\fill[color=zzttqq,fill=zzttqq,fill opacity=0.05] (2.1,2.78) -- (3.64,1.36) -- (5.04,1.94) -- (5.72,3.18) -- cycle;
\draw[color=zzttqq] (2.1,2.78) -- (3.64,1.36) -- (5.04,1.94) -- (5.72,3.18) -- cycle;
\fill[color=zzttqq,fill=zzttqq,fill opacity=0.05] (5.72,3.18) -- (5.04,1.94) -- (6.24,1.01) -- (8,0.74) -- (7.94,-0.42) -- (9.7,1.04) -- (7.46,1.32) -- (7.96,2.94) -- cycle;
\draw[color=zzttqq] (5.72,3.18) -- (5.04,1.94) -- (6.24,1.01) -- (8,0.74) -- (7.94,-0.42) -- (9.7,1.04) -- (7.46,1.32) -- (7.96,2.94) -- cycle;
\fill[color=zzttqq,fill=zzttqq,fill opacity=0.05] (7.94,-0.42) -- (8,0.74) -- (6.24,1.01) -- (6.04,-0.49) -- (7.12,-1.82) -- cycle;
\draw[color=zzttqq] (7.94,-0.42) -- (8,0.74) -- (6.24,1.01) -- (6.04,-0.49) -- (7.12,-1.82) -- cycle;
\fill[color=zzttqq,fill=zzttqq,fill opacity=0.05] (7.12,-1.82) -- (6.04,-0.49) -- (4.64,-1.06) -- (3.88,-2.22) -- (6.16,-3) -- (7.06,-3.02) -- cycle;
\draw[color=zzttqq] (7.12,-1.82) -- (6.04,-0.49) -- (4.64,-1.06) -- (3.88,-2.22) -- (6.16,-3) -- (7.06,-3.02) -- cycle;
\fill[color=zzttqq,fill=zzttqq,fill opacity=0.05] (3.88,-2.22) -- (4.64,-1.06) -- (3.44,-0.14) -- (0.94,0.36) -- (-0.28,-0.28) -- (1.08,-1.86) -- cycle;
\draw[color=zzttqq] (3.88,-2.22) -- (4.64,-1.06) -- (3.44,-0.14) -- (0.94,0.36) -- (-0.28,-0.28) -- (1.08,-1.86) -- cycle;
\fill[color=zzttqq,fill=zzttqq,fill opacity=0.05] (7.06,-3.02) -- (7.12,-1.82) -- (7.94,-0.42) -- (9.7,1.04) -- (9.86,-1.54) -- cycle;
\draw[color=zzttqq] (7.06,-3.02) -- (7.12,-1.82) -- (7.94,-0.42) -- (9.7,1.04) -- (9.86,-1.54) -- cycle;
\fill[color=zzttqq,fill=zzttqq,fill opacity=0.05] (9.7,1.04) -- (7.46,1.32) -- (7.96,2.94) -- cycle;
\draw[color=zzttqq] (9.7,1.04) -- (7.46,1.32) -- (7.96,2.94) -- cycle;
\fill[color=zzttqq,fill=zzttqq,fill opacity=0.05] (0.94,0.36) -- (0.24,1.46) -- (-0.28,-0.28) -- cycle;
\draw[color=zzttqq] (0.94,0.36) -- (0.24,1.46) -- (-0.28,-0.28) -- cycle;
\fill[color=zzttqq,fill=zzttqq,fill opacity=0.05] (7.06,-3.02) -- (6.7,-4.14) -- (6.16,-3) -- cycle;
\draw[color=zzttqq] (7.06,-3.02) -- (6.7,-4.14) -- (6.16,-3) -- cycle;
\fill[color=zzttqq,fill=zzttqq,fill opacity=0.05] (6.7,-4.14) -- (6.16,-3) -- (3.88,-2.22) -- (1.08,-1.86) -- cycle;
\draw[color=zzttqq] (6.7,-4.14) -- (6.16,-3) -- (3.88,-2.22) -- (1.08,-1.86) -- cycle;
\fill[color=zzttqq,fill=zzttqq,fill opacity=0.05] (6.7,-4.14) -- (7.06,-3.02) -- (9.86,-1.54) -- cycle;
\draw[color=zzttqq] (6.7,-4.14) -- (7.06,-3.02) -- (9.86,-1.54) -- cycle;
\draw [line width=1pt] (1.89,2.63) -- (3.57,0.82) -- (5.18,1.83) -- (5.45,2.68) -- (4.69,1.79) -- (5.29,1.74) -- (5.15,2.14) -- (3.27,1.71) -- (3.01,-0.05) -- (2.5,-2.04) -- (2.5,-2.44);
\draw [line width=1pt,blue] (4.72,-3.34)-- (4.57,-2.46) -- (6.33,-0.84) -- (7.98,0.34) -- (8.4,-0.04) -- (9.44,-1.76) -- (9.48,-1.85);
\end{tikzpicture}
\caption{Example of two trajectories of a polygonal tiling billiard}\label{Fig:ex_polTB}
\small
Only a part of the tiling is drawn, hence only a part of each trajectory.
\end{figure}

For a fixed tiling of the plane, we are interested in the following questions. What sort of trajectories exist? Periodic, unbounded, dense in open sets of the plane? Are they stable under perturbations of initial conditions of the trajectory (starting point and direction) or of the form of the tiling?

There exist examples of tiling billiards with \emph{periodic}, \emph{drift-periodic} or \emph{linearly escaping} trajectories.  We call \textbf{drift-periodic} a trajectory with a translation symmetry. We say that a trajectory is \textbf{linearly escaping} if it is not drift-periodic and the main part of the asymptotic expansion of its coodinates (by arclength parametrization) is linear and the remaining term is bounded by a constant.
Another behavior is known: in the trihexagonal tiling billiard, the generic trajectories are dense in open sets of the plane \cite{DH}.
\subsubsection*{Triangle and cyclic quadrilateral tilings}
Let us consider a triangle tiling where each tile is a centrally symmetric copy of each of its neighboring tiles. 
The qualitative behavior of triangle tiling billiards has been fully classified in \cite{PR19}, building upon \cite{BSDFI18} and \cite{HPR22}. 
We say that a trajectory is \textbf{stable} is a small perturbation of it leads to another trajectory that crosses the same sides (in the same order).
\begin{thm}[\cite{BSDFI18}, \cite{PR19}, \cite{HPR22}]
Let $\delta$ be a triangle. Let $\gamma$ be a trajectory in the $\delta$-tiling billiard. Then exactly one of the following four cases holds:
\begin{enumerate}
\item 
The trajectory $\gamma$ is periodic and stable under pertubation of the tiling or of the initial conditions.
\item
The trajectory $\gamma$ is drift-periodic and the angles of $\delta$ are rationnaly dependant.
\item
The trajectory $\gamma$ is linearly escaping.
\item
The trajectory $\gamma$ is not bounded, has no asymptotic direction and passes through circumcenter of triangles. Moreover $\delta$ belongs to the Rauzy gasket, which has a zero Lebesgue measure.
\end{enumerate}  
\end{thm}

\begin{figure}[h!]
     \centering
\begin{subfigure}[b]{0.4\textwidth}
	\centering
	\includegraphics[trim=1.5cm 1cm 1cm 0.5cm,clip,scale=0.3]{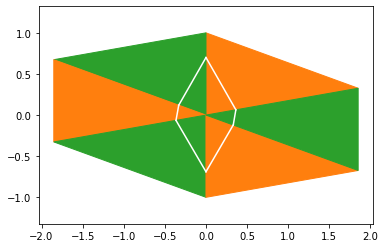}
	
	\includegraphics[trim=1.5cm 1cm 1cm 0.5cm,clip,scale=0.3]{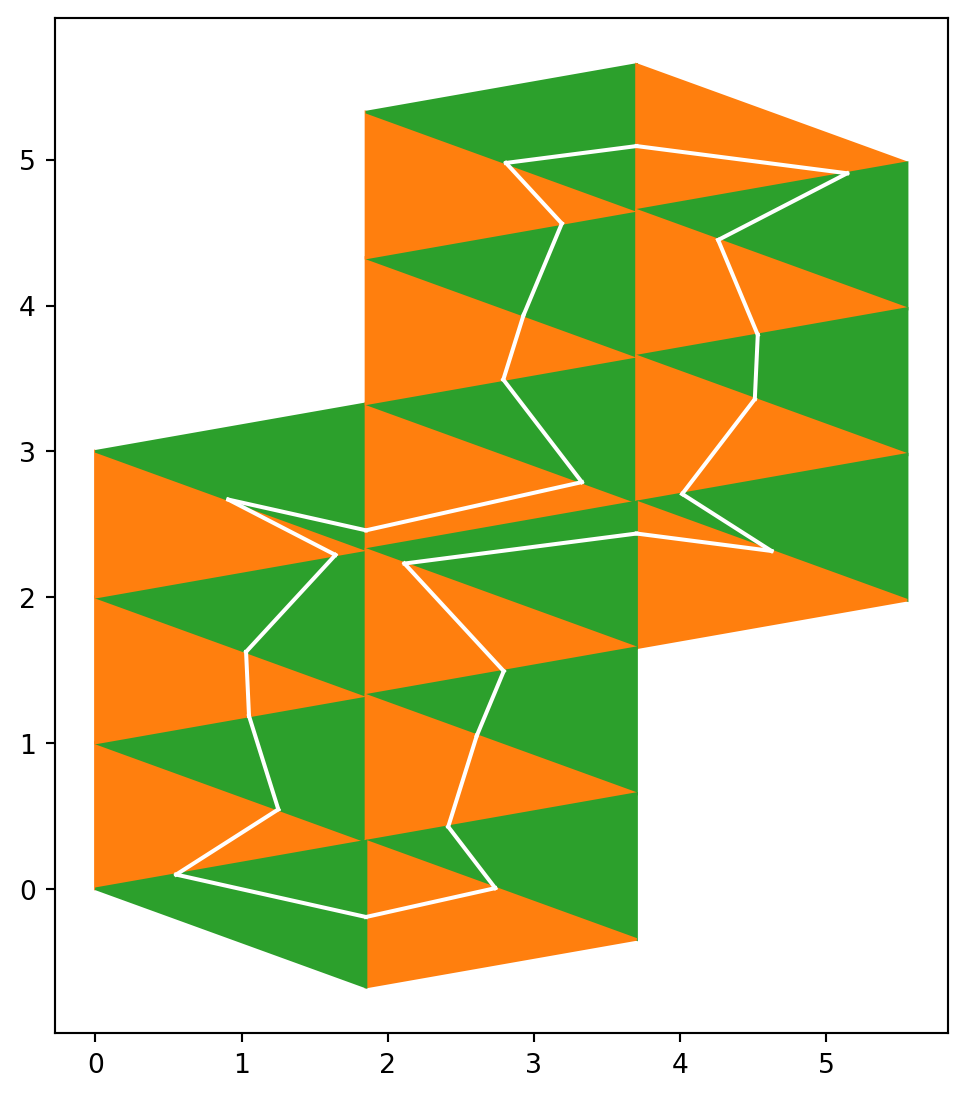}

	\includegraphics[trim=1.5cm 1cm 1cm 0.5cm,clip,scale=0.3]{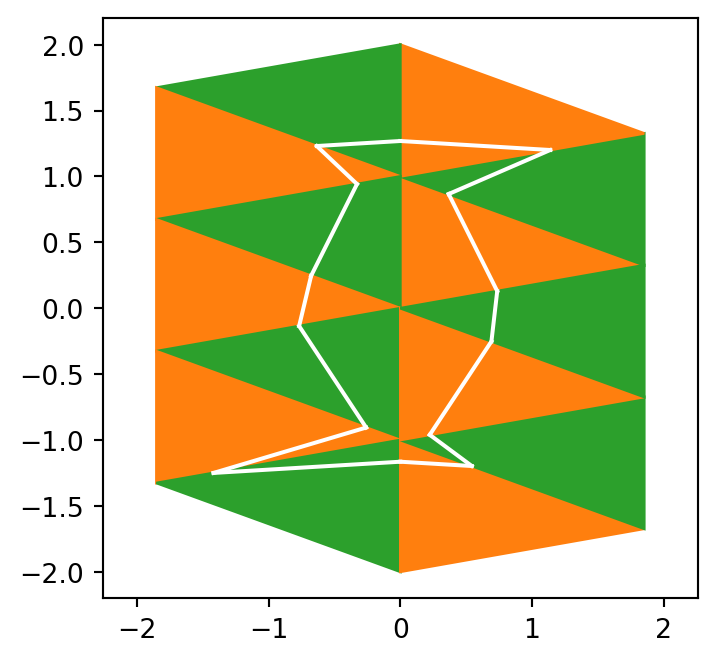}

	\includegraphics[trim=1.5cm 1cm 1cm 0.5cm,clip,scale=0.3]{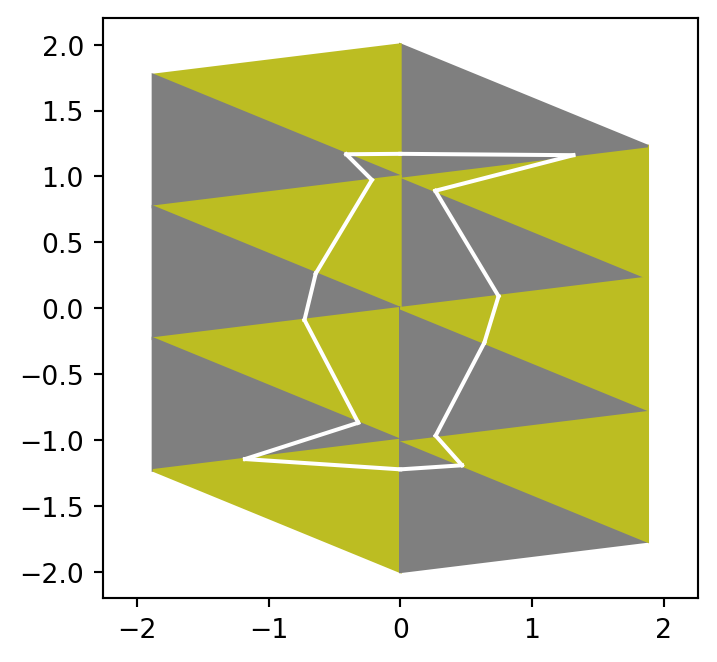}
	\includegraphics[trim=1.5cm 1cm 1cm 0.5cm,clip,scale=0.3]{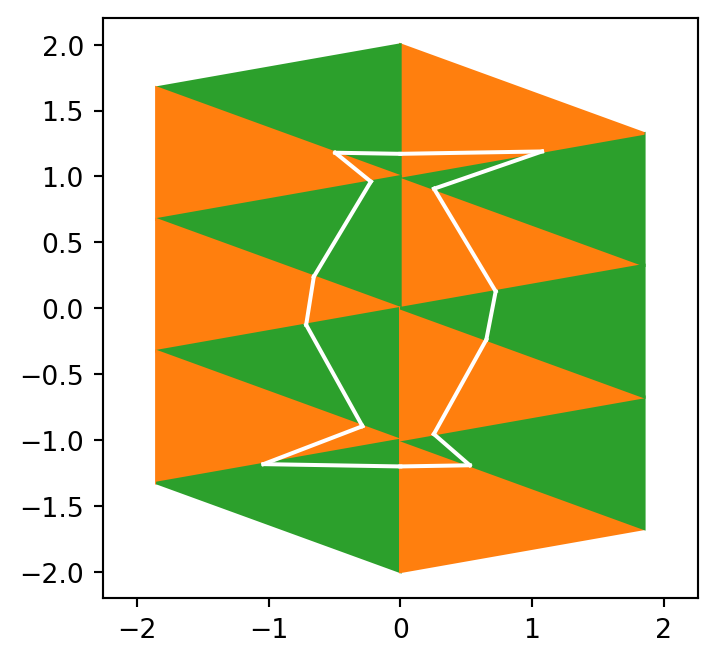}
     \caption{Some periodic trajectories and two (stable) pertubations}
\begin{footnotesize}\begin{flushleft}
The green and grey tiling is a perturbation of the green and orange one.
\end{flushleft}
 \end{footnotesize} 
\end{subfigure}
\hfill   
\begin{subfigure}[b]{0.55\textwidth}
	\centering
\vspace*{-0.2cm}
	\includegraphics[trim=1.5cm 1cm 1cm 1.5cm,clip,scale=0.3]{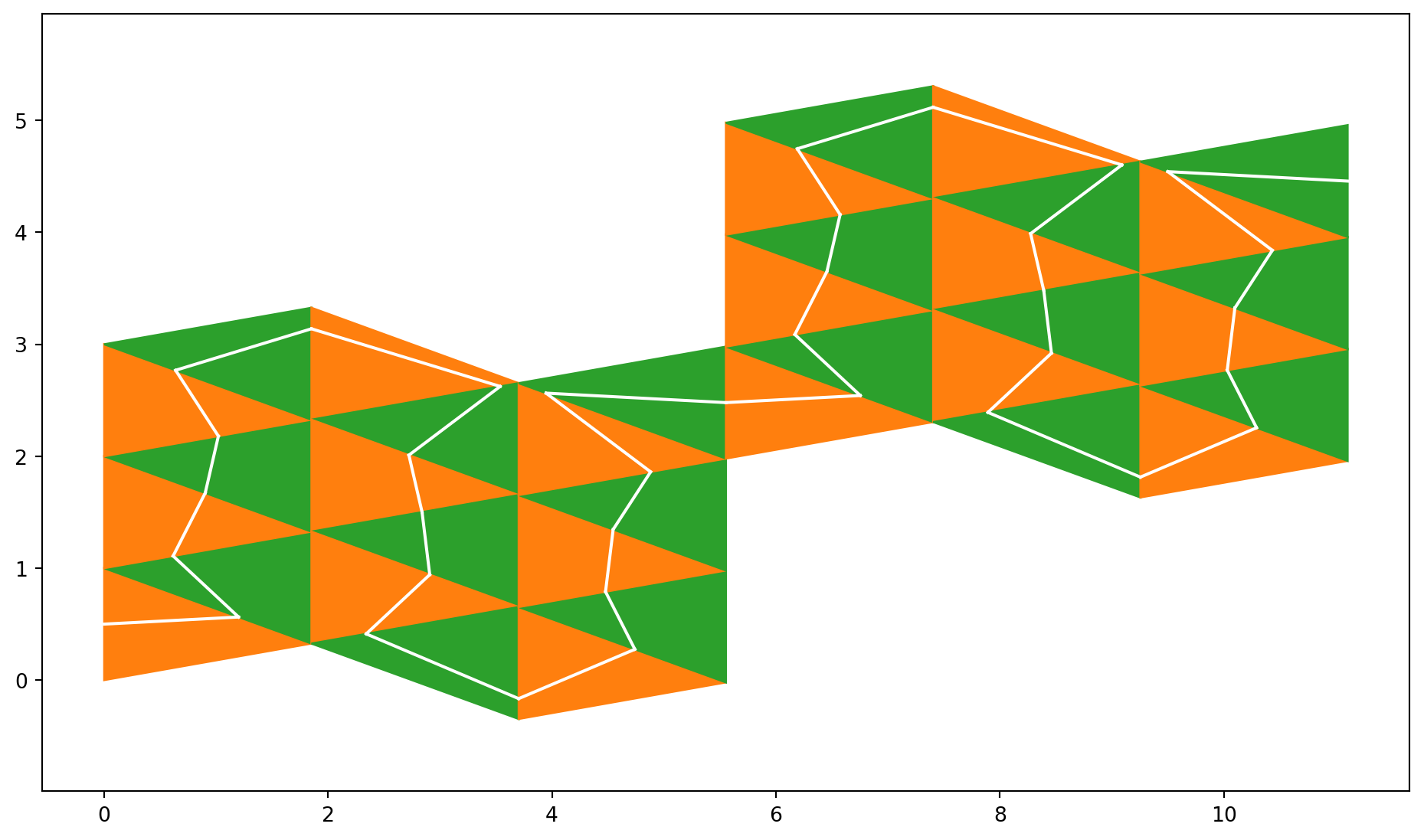}
\vspace*{-0.2cm}
	\includegraphics[trim=1.5cm 1cm 1cm 1.5cm,clip,scale=0.3]{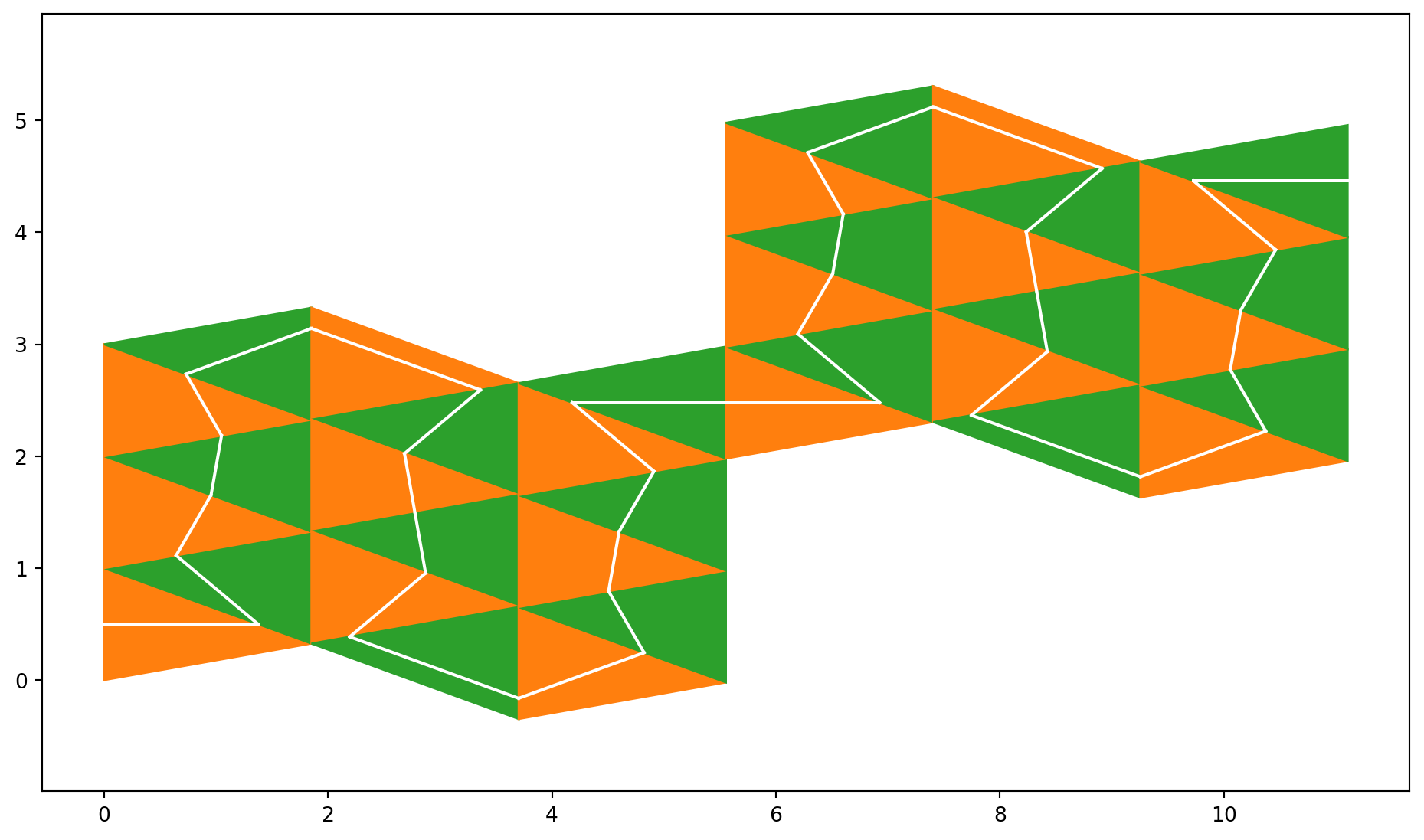}
\vspace*{-0.2cm}
	\includegraphics[trim=1.5cm 1cm 1cm 1cm,clip,scale=0.3]{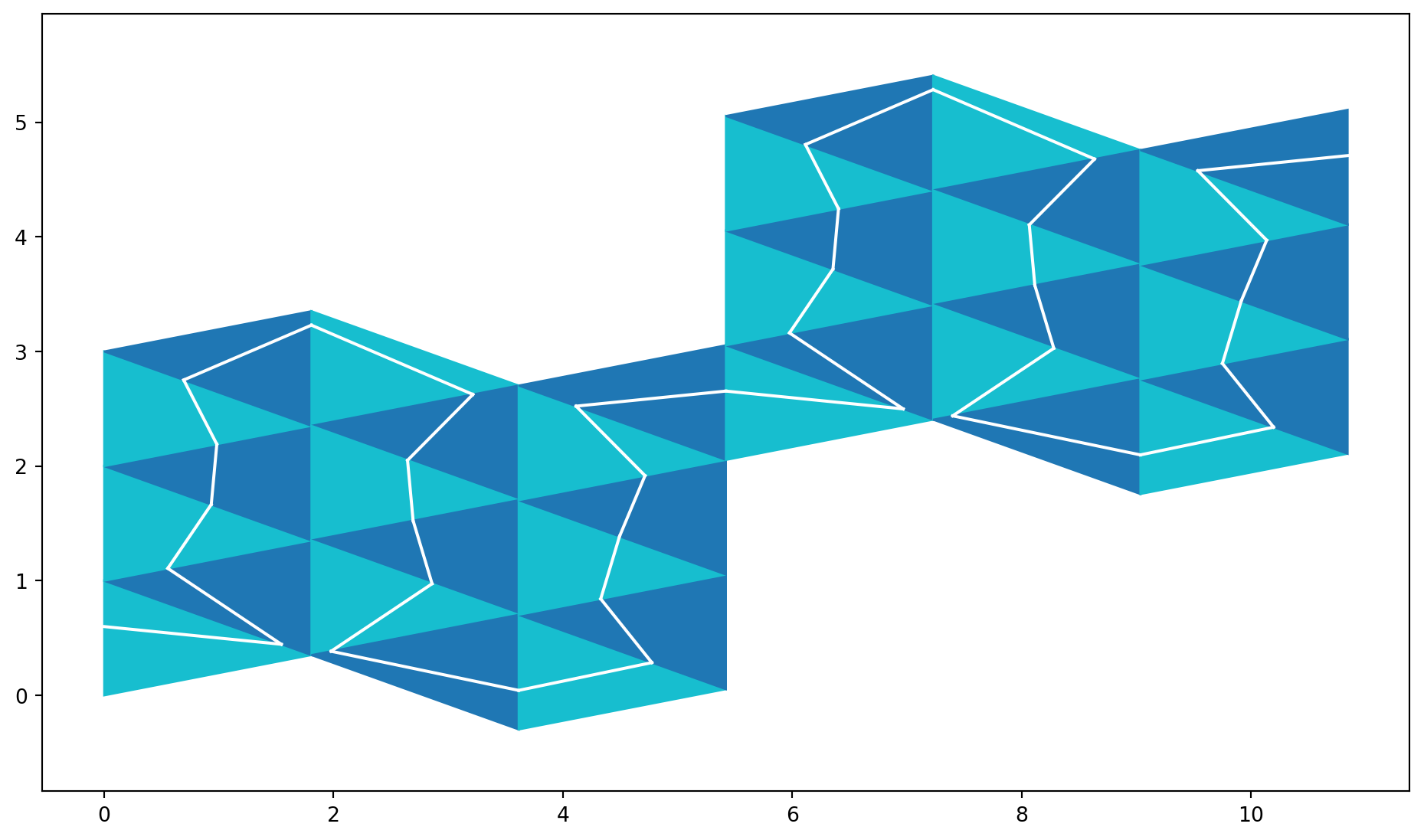}
     \caption{A part of a drift-periodic trajectory and two of its (stable) perturbations}
\begin{footnotesize}
\begin{flushleft}
 The light and deep blue tiling is a perturbation of the green and orange one.
\end{flushleft}
 \end{footnotesize}
\end{subfigure}

\begin{subfigure}[b]{0.98\textwidth}
	\centering
	\includegraphics[trim=0.5cm 1cm 0cm 1.5cm,clip,scale=0.34,angle=90]{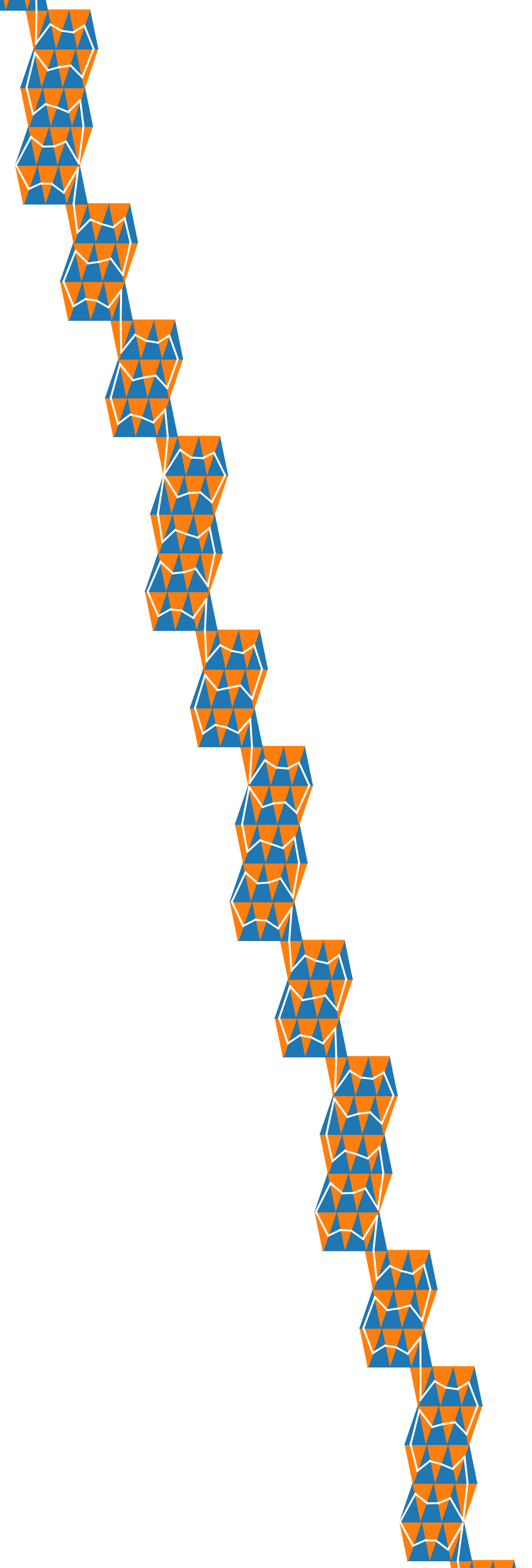}
	     \caption{A part of a linearly escaping trajectory}
\end{subfigure}
        \caption{Examples of trajectories in triangle tiling billiards}
        \label{fig:tri_tb}
\end{figure}

See Figure \ref{fig:tri_tb} for an illustration of the first three cases. 
The exceptional case of unbounded trajectories with no asymptotic direction is linked to minimal IETFs. 
Hubert and Paris-Romaskevich with Dynnikov, Mercat and Skripchenko
\cite{DHMPRS} obtained the analogous result for classication of trajectories on cyclic quadrilateral tilings, and related the dynamics of such billiards to the Novikov’s problem of the study of foliations on subsurfaces of genus 3 of the 3-torus given by 1-forms, see also \cite{PR21}.

For triangle and cyclic quadrilateral tiling, the dynamical systems have an integral of motion: the oriented distance of a trajectory to the circumcenters of the crossed polygons is preserved. 
Baird-Smith, Davis, Fromm, Iyer explain this for triangles in \cite{BSDFI18} (Lemma 3.2), with a beautiful argument of folding. It was then generalized in \cite{PR19}.
This allows us to study such tiling billiards via IETFs.
See Section \ref{subsec_system-coord} for another explanation of the link between tiling billiards and IETFs.

\vspace{0.5cm}
The behaviors of trajectories in tiling billiards tend to differ from the ones in classical inner billiards, as dynamics of IETFs differs from those of interval exchange transformations without flips (IETs).
No general study of tiling billiards exists yet.
In this article, we will handle open sets of parameters in generalized tiling billiards, made with any cyclic $N$-gons with $N\geqslant 5$, although they will not tile the plane.
The results are true for polygons with more than five sides but the reader can think of pentagons in a first reading. All polygons will be drawn as pentagons.

\subsection{The studied system}
\subsubsection*{The generalized tiling}
\begin{figure}[h]
 \centering
\begin{subfigure}[b]{0.4\textwidth}
         \centering
\definecolor{qqzzqq}{rgb}{0,0.6,0}
\definecolor{ffccqq}{rgb}{1,0.6,0}
\definecolor{ffttqq}{rgb}{1,0.3,0.4}
\definecolor{qqqqcc}{rgb}{0,0,0.8}
\definecolor{tttttt}{rgb}{0.2,0.2,0.2}
\begin{tikzpicture}[line cap=round,line join=round,>=triangle 45,x=0.35cm,y=0.35cm]
\fill[color=tttttt,fill=tttttt,fill opacity=0.1] (4.18,1.4) -- (-3.28,0.85) -- (0.14,-2.21) -- (2.47,-1.62) -- (3.59,-0.49) -- cycle;
\fill[color=qqqqcc,fill=qqqqcc,fill opacity=0.1] (-3.28,0.85) -- (4.18,1.4) -- (0.76,4.46) -- (-1.57,3.86) -- (-2.68,2.74) -- cycle;
\fill[color=ffttqq,fill=ffttqq,fill opacity=0.1] (8.22,5.01) -- (0.76,4.46) -- (4.18,1.4) -- (6.51,2) -- (7.62,3.12) -- cycle;
\fill[color=ffccqq,fill=ffccqq,fill opacity=0.1] (2.47,-1.62) -- (9.93,-1.06) -- (6.51,2) -- (4.18,1.4) -- (3.06,0.28) -- cycle;
\fill[color=qqzzqq,fill=qqzzqq,fill opacity=0.1] (4.77,3.29) -- (-2.68,2.74) -- (0.74,-0.32) -- (3.06,0.28) -- (4.18,1.4) -- cycle;

\draw[color=tttttt] (4.18,1.4) -- (-3.28,0.85) -- (0.14,-2.21) -- (2.47,-1.62) -- (3.59,-0.49) -- cycle;
\draw[color=qqqqcc] (-3.28,0.85) -- (4.18,1.4) -- (0.76,4.46) -- (-1.57,3.86) -- (-2.68,2.74) -- cycle;
\draw[color=ffttqq] (8.22,5.01) -- (0.76,4.46) -- (4.18,1.4) -- (6.51,2) -- (7.62,3.12) -- cycle;
\draw[color=ffccqq] (2.47,-1.62) -- (9.93,-1.06) -- (6.51,2) -- (4.18,1.4) -- (3.06,0.28) -- cycle;
\draw[color=qqzzqq] (4.77,3.29) -- (-2.68,2.74) -- (0.74,-0.32) -- (3.06,0.28) -- (4.18,1.4) -- cycle;

\draw [shift={(4.18,1.4)},fill=black,fill opacity=0.1] (0,0) -- (-175.78:0.3) arc (-175.78:-107.42:0.3) -- cycle;
\draw [shift={(3.59,-0.49)},color=qqzzqq,fill=qqzzqq,fill opacity=0.1] (0,0) -- (72.58:0.7) arc (72.58:225.15:0.7) -- cycle;
\draw [shift={(2.47,-1.62)},color=ffccqq,fill=ffccqq,fill opacity=0.1] (0,0) -- (45.15:0.6) arc (45.15:194.4:0.6) -- cycle;
\draw [shift={(0.14,-2.21)},color=ffttqq,fill=ffttqq,fill opacity=0.1] (0,0) -- (14.4:0.5) arc (14.4:138.18:0.5) -- cycle;
\draw [shift={(-3.28,0.85)},color=qqqqcc,fill=qqqqcc,fill opacity=0.1] (0,0) -- (-41.82:0.4) arc (-41.82:4.22:0.4) -- cycle;
\draw [shift={(4.18,1.4)},color=qqqqcc,fill=qqqqcc,fill opacity=0.1] (0,0) -- (138.18:0.4) arc (138.18:184.22:0.4) -- cycle;
\draw [shift={(4.18,1.4)},color=ffttqq,fill=ffttqq,fill opacity=0.1] (0,0) -- (14.4:0.5) arc (14.4:138.18:0.5) -- cycle;
\draw [shift={(4.18,1.4)},color=ffccqq,fill=ffccqq,fill opacity=0.1] (0,0) -- (-134.85:0.6) arc (-134.85:14.4:0.6) -- cycle;
\draw [shift={(4.18,1.4)},color=qqzzqq,fill=qqzzqq,fill opacity=0.1] (0,0) -- (72.58:0.7) arc (72.58:225.15:0.7) -- cycle;
\end{tikzpicture}
\caption{Polygons either overlap...}
\end{subfigure}
\hfill
\begin{subfigure}[b]{0.4\textwidth}
         \centering
\definecolor{ffccqq}{rgb}{1,0.6,0}
\definecolor{qqzzqq}{rgb}{0,0.6,0}
\definecolor{ffttqq}{rgb}{1,0.3,0.4}
\definecolor{qqqqcc}{rgb}{0,0,0.8}
\definecolor{tttttt}{rgb}{0.2,0.2,0.2}
\begin{tikzpicture}[line cap=round,line join=round,>=triangle 45,x=0.35cm,y=0.35cm]
\fill[color=tttttt,fill=tttttt,fill opacity=0.1] (4.18,1.4) -- (-3.28,0.85) -- (0.14,-2.21) -- (2.47,-1.62) -- (3.59,-0.49) -- cycle;
\fill[color=qqqqcc,fill=qqqqcc,fill opacity=0.1] (-3.28,0.85) -- (4.18,1.4) -- (0.76,4.46) -- (-1.57,3.86) -- (-2.68,2.74) -- cycle;
\fill[color=ffttqq,fill=ffttqq,fill opacity=0.1] (8.22,5.01) -- (0.76,4.46) -- (4.18,1.4) -- (6.51,2) -- (7.62,3.12) -- cycle;
\fill[color=ffccqq,fill=ffccqq,fill opacity=0.1] (5.91,0.1) -- (13.37,0.65) -- (9.95,3.72) -- (7.62,3.12) -- (6.51,2) -- cycle;

\draw[color=tttttt] (4.18,1.4) -- (-3.28,0.85) -- (0.14,-2.21) -- (2.47,-1.62) -- (3.59,-0.49) -- cycle;
\draw[color=qqqqcc] (-3.28,0.85) -- (4.18,1.4) -- (0.76,4.46) -- (-1.57,3.86) -- (-2.68,2.74) -- cycle;
\draw[color=ffttqq] (8.22,5.01) -- (0.76,4.46) -- (4.18,1.4) -- (6.51,2) -- (7.62,3.12) -- cycle;
\draw[color=ffccqq] (5.91,0.1) -- (13.37,0.65) -- (9.95,3.72) -- (7.62,3.12) -- (6.51,2) -- cycle;

\draw [shift={(4.18,1.4)},fill=black,fill opacity=0.1] (0,0) -- (-175.78:0.3) arc (-175.78:-107.42:0.3) -- cycle;
\draw [shift={(3.59,-0.49)},color=qqzzqq,fill=qqzzqq,fill opacity=0.1] (0,0) -- (72.58:0.7) arc (72.58:225.15:0.7) -- cycle;
\draw [shift={(2.47,-1.62)},color=ffccqq,fill=ffccqq,fill opacity=0.1] (0,0) -- (45.15:0.6) arc (45.15:194.4:0.6) -- cycle;
\draw [shift={(0.14,-2.21)},color=ffttqq,fill=ffttqq,fill opacity=0.1] (0,0) -- (14.4:0.5) arc (14.4:138.18:0.5) -- cycle;
\draw [shift={(-3.28,0.85)},color=qqqqcc,fill=qqqqcc,fill opacity=0.1] (0,0) -- (-41.82:0.4) arc (-41.82:4.22:0.4) -- cycle;
\draw [shift={(4.18,1.4)},color=qqqqcc,fill=qqqqcc,fill opacity=0.1] (0,0) -- (138.18:0.4) arc (138.18:184.22:0.4) -- cycle;
\end{tikzpicture}
\caption{... or leave empty space}
\end{subfigure}
\caption{The way we put polygons together to play tiling billiards will not tile the plane with $N$-gons, $N\geqslant 5$.}\label{Fig:penta_around_vertex}
\end{figure}
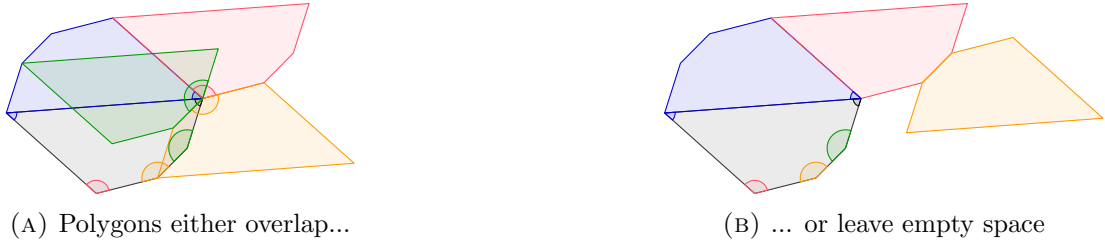

We first define how we build the generalized tiling with the $N$-gons ($N\geqslant 5$). Each one of its tiles is congruent to a fixed cyclic polygon. Moreover if $T$ and $T'$ are two neighboring tiles then $T'$ is the image of $T$ via a central symmetry with respect to the middle of their shared side. This does not define a proper tiling: polygons will either overlap each other or not fill the entire plane. Indeed when turning around a vertex, all the angles of the $N$-gon appear and sum up to $(N-2)\pi$, which is not a divisor of $2\pi$, see Figure~\ref{Fig:penta_around_vertex}.
But we can still define a billiard trajectory in the plane, generalizing tiling billiards for any polygon.

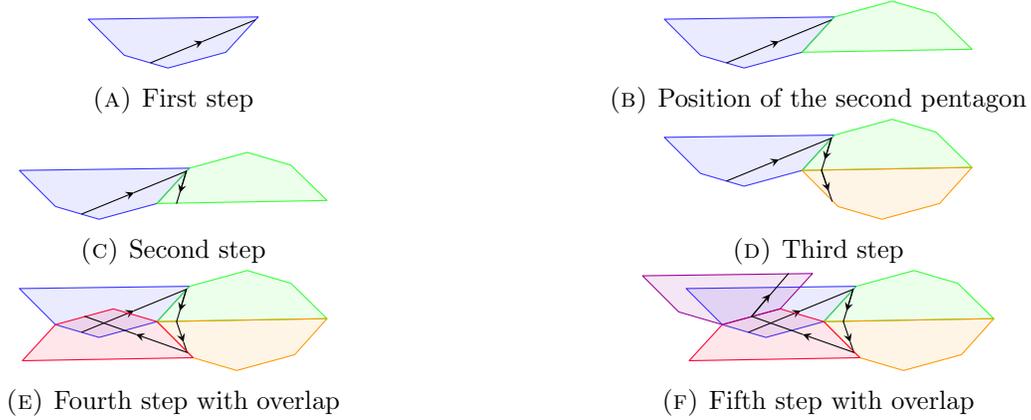
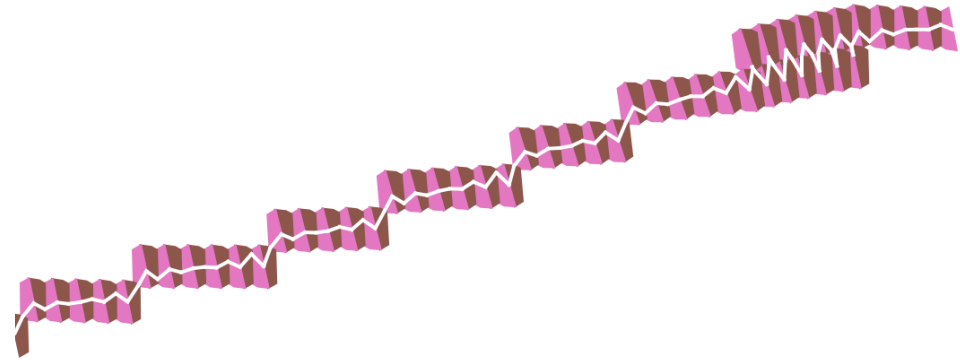
\begin{figure}[h!]
     \centering
\begin{subfigure}[b]{0.4\textwidth}
         \centering
         \definecolor{zzqqzz}{rgb}{0.6,0,0.6}
\definecolor{ffqqtt}{rgb}{1,0,0.2}
\definecolor{ffzzqq}{rgb}{1,0.6,0}
\definecolor{ttfftt}{rgb}{0.2,1,0.2}
\definecolor{ttttff}{rgb}{0.2,0.2,1}
\begin{tikzpicture}[line cap=round,line join=round,>=stealth,x=1.0cm,y=1.0cm]
\fill[color=ttttff,fill=ttttff,fill opacity=0.1] (4.18,1.39) -- (3.75,0.92) -- (2.98,0.71) -- (2.4,0.88) -- (1.92,1.36) -- cycle;
\draw [color=ttttff] (4.18,1.39)-- (3.75,0.92) -- (2.98,0.71) -- (2.4,0.88) -- (1.92,1.36) -- cycle;

\draw [->] (2.75,0.78) -- (3.445,1.065);
\draw (3.445,1.065) -- (4.14,1.35);
	\end{tikzpicture}
	\caption{First step}
    \label{fig_overlap1}
	\end{subfigure} 
\hfill
    \begin{subfigure}[b]{0.5\textwidth}
\centering
\definecolor{zzqqzz}{rgb}{0.6,0,0.6}
\definecolor{ffqqtt}{rgb}{1,0,0.2}
\definecolor{ffzzqq}{rgb}{1,0.6,0}
\definecolor{ttfftt}{rgb}{0.2,1,0.2}
\definecolor{ttttff}{rgb}{0.2,0.2,1}
	\begin{tikzpicture}[line cap=round,line join=round,>=stealth,x=1.0cm,y=1.0cm]
\draw [color=ttttff] (4.18,1.39)-- (3.75,0.92) -- (2.98,0.71) -- (2.4,0.88) -- (1.92,1.36) -- cycle;
\fill [color=ttttff,fill=ttttff,fill opacity=0.1] (4.18,1.39) -- (3.75,0.92) -- (2.98,0.71) -- (2.4,0.88) -- (1.92,1.36) -- cycle;

\draw [color=ttfftt] (3.75,0.92)-- (4.18,1.39) -- (4.95,1.6) -- (5.53,1.43) -- (6.01,0.96) -- cycle;
\fill [color=ttfftt,fill=ttfftt,fill opacity=0.1] (3.75,0.92) -- (4.18,1.39) -- (4.95,1.6) -- (5.53,1.43) -- (6.01,0.96) -- cycle;

\draw [->] (2.75,0.78) -- (3.445,1.065);
\draw (3.445,1.065) -- (4.14,1.35);
	\end{tikzpicture}
    \caption{Position of the second pentagon}
    \label{fig_overlap2penta}
	\end{subfigure} 
\hfill
    \begin{subfigure}[b]{0.4\textwidth}
\centering
\definecolor{zzqqzz}{rgb}{0.6,0,0.6}
\definecolor{ffqqtt}{rgb}{1,0,0.2}
\definecolor{ffzzqq}{rgb}{1,0.6,0}
\definecolor{ttfftt}{rgb}{0.2,1,0.2}
\definecolor{ttttff}{rgb}{0.2,0.2,1}
	\begin{tikzpicture}[line cap=round,line join=round,>=stealth,x=1.0cm,y=1.0cm]
\draw [color=ttttff] (4.18,1.39) -- (3.75,0.92) -- (2.98,0.71) -- (2.4,0.88) -- (1.92,1.36) -- cycle;
\fill [color=ttttff,fill=ttttff,fill opacity=0.1] (4.18,1.39) -- (3.75,0.92) -- (2.98,0.71) -- (2.4,0.88) -- (1.92,1.36) -- cycle;

\draw [color=ttfftt] (3.75,0.92) -- (4.18,1.39) -- (4.95,1.6) -- (5.53,1.43) -- (6.01,0.96) -- cycle;
\fill[color=ttfftt,fill=ttfftt,fill opacity=0.1] (3.75,0.92) -- (4.18,1.39) -- (4.95,1.6) -- (5.53,1.43) -- (6.01,0.96) -- cycle;

\draw [->] (2.75,0.78) -- (3.445,1.065);
\draw (3.445,1.065) -- (4.14,1.35);

\draw [->] (4.14,1.35) -- (4.053,1.063);
\draw (4.14,1.35) -- (4.01,0.92);
	\end{tikzpicture}
    \caption{Second step}
    \label{fig_overlap2line}
	\end{subfigure} 
\hfill
    \begin{subfigure}[b]{0.5\textwidth}
\centering
\definecolor{zzqqzz}{rgb}{0.6,0,0.6}
\definecolor{ffqqtt}{rgb}{1,0,0.2}
\definecolor{ffzzqq}{rgb}{1,0.6,0}
\definecolor{ttfftt}{rgb}{0.2,1,0.2}
\definecolor{ttttff}{rgb}{0.2,0.2,1}
	\begin{tikzpicture}[line cap=round,line join=round,>=stealth,x=1.0cm,y=1.0cm]
\draw [color=ttttff] (4.18,1.39) -- (3.75,0.92) -- (2.98,0.71) -- (2.4,0.88) -- (1.92,1.36) -- cycle;
\fill [color=ttttff,fill=ttttff,fill opacity=0.1] (4.18,1.39) -- (3.75,0.92) -- (2.98,0.71) -- (2.4,0.88) -- (1.92,1.36) -- cycle;

\draw [color=ttfftt] (3.75,0.92) -- (4.18,1.39) -- (4.95,1.6) -- (5.53,1.43) -- (6.01,0.96) -- cycle;
\fill[color=ttfftt,fill=ttfftt,fill opacity=0.1] (3.75,0.92) -- (4.18,1.39) -- (4.95,1.6) -- (5.53,1.43) -- (6.01,0.96) -- cycle;

\draw [color=ffzzqq] (6.01,0.96) -- (5.58,0.48) -- (4.81,0.27) -- (4.23,0.44) -- (3.75,0.92) -- cycle;
\fill[color=ffzzqq,fill=ffzzqq,fill opacity=0.1] (6.01,0.96) -- (5.58,0.48) -- (4.81,0.27) -- (4.23,0.44) -- (3.75,0.92) -- cycle;

\draw [->] (2.75,0.78) -- (3.445,1.065);
\draw (3.445,1.065) -- (4.14,1.35);

\draw [->] (4.14,1.35) -- (4.053,1.063);
\draw (4.14,1.35) -- (4.01,0.92);

\draw [->] (4.01,0.92) -- (4.103,0.647);
\draw (4.01,0.92)-- (4.15,0.51);
	\end{tikzpicture}
    \caption{Third step}
    \label{fig_overlap3}
	\end{subfigure} 
\hfill
    \begin{subfigure}[b]{0.4\textwidth}
\centering
\definecolor{zzqqzz}{rgb}{0.6,0,0.6}
\definecolor{ffqqtt}{rgb}{1,0,0.2}
\definecolor{ffzzqq}{rgb}{1,0.6,0}
\definecolor{ttfftt}{rgb}{0.2,1,0.2}
\definecolor{ttttff}{rgb}{0.2,0.2,1}
	\begin{tikzpicture}[line cap=round,line join=round,>=stealth,x=1.0cm,y=1.0cm]
\draw [color=ttttff] (4.18,1.39) -- (3.75,0.92) -- (2.98,0.71) -- (2.4,0.88) -- (1.92,1.36) -- cycle;
\fill [color=ttttff,fill=ttttff,fill opacity=0.1] (4.18,1.39) -- (3.75,0.92) -- (2.98,0.71) -- (2.4,0.88) -- (1.92,1.36) -- cycle;

\draw [color=ttfftt] (3.75,0.92) -- (4.18,1.39) -- (4.95,1.6) -- (5.53,1.43) -- (6.01,0.96) -- cycle;
\fill[color=ttfftt,fill=ttfftt,fill opacity=0.1] (3.75,0.92) -- (4.18,1.39) -- (4.95,1.6) -- (5.53,1.43) -- (6.01,0.96) -- cycle;

\draw [color=ffzzqq] (6.01,0.96) -- (5.58,0.48) -- (4.81,0.27) -- (4.23,0.44) -- (3.75,0.92) -- cycle;
\fill[color=ffzzqq,fill=ffzzqq,fill opacity=0.1] (6.01,0.96) -- (5.58,0.48) -- (4.81,0.27) -- (4.23,0.44) -- (3.75,0.92) -- cycle;

\draw [->] (2.75,0.78) -- (3.445,1.065);
\draw (3.445,1.065) -- (4.14,1.35);

\draw [->] (4.14,1.35) -- (4.053,1.063);
\draw (4.14,1.35) -- (4.01,0.92);

\draw [->] (4.01,0.92) -- (4.103,0.647);
\draw (4.103,0.647) -- (4.15,0.51);

\draw [color=ffqqtt] (1.96,0.4) -- (2.4,0.88) -- (3.17,1.09) -- (3.75,0.92) -- (4.23,0.44) -- cycle;
\fill[color=ffqqtt,fill=ffqqtt,fill opacity=0.1] (1.96,0.4) -- (2.4,0.88) -- (3.17,1.09) -- (3.75,0.92) -- (4.23,0.44) -- cycle;

\draw [->] (4.15,0.51) -- (3.47,0.75);
\draw (3.47,0.75) -- (2.79,0.99);
	\end{tikzpicture}
    \caption{Fourth step with overlap}
    \label{fig_overlap4}
	\end{subfigure} 
    \hfill
	\begin{subfigure}[b]{0.5\textwidth}
\centering
\definecolor{zzqqzz}{rgb}{0.6,0,0.6}
\definecolor{ffqqtt}{rgb}{1,0,0.2}
\definecolor{ffzzqq}{rgb}{1,0.6,0}
\definecolor{ttfftt}{rgb}{0.2,1,0.2}
\definecolor{ttttff}{rgb}{0.2,0.2,1}
\begin{tikzpicture}[line cap=round,line join=round,>=stealth,x=1.0cm,y=1.0cm]
\draw [color=ttttff] (4.18,1.39) -- (3.75,0.92) -- (2.98,0.71) -- (2.4,0.88) -- (1.92,1.36) -- cycle;
\fill [color=ttttff,fill=ttttff,fill opacity=0.1] (4.18,1.39) -- (3.75,0.92) -- (2.98,0.71) -- (2.4,0.88) -- (1.92,1.36) -- cycle;

\draw [color=ttfftt] (3.75,0.92) -- (4.18,1.39) -- (4.95,1.6) -- (5.53,1.43) -- (6.01,0.96) -- cycle;
\fill[color=ttfftt,fill=ttfftt,fill opacity=0.1] (3.75,0.92) -- (4.18,1.39) -- (4.95,1.6) -- (5.53,1.43) -- (6.01,0.96) -- cycle;

\draw [color=ffzzqq] (6.01,0.96) -- (5.58,0.48) -- (4.81,0.27) -- (4.23,0.44) -- (3.75,0.92) -- cycle;
\fill[color=ffzzqq,fill=ffzzqq,fill opacity=0.1] (6.01,0.96) -- (5.58,0.48) -- (4.81,0.27) -- (4.23,0.44) -- (3.75,0.92) -- cycle;

\draw [->] (2.75,0.78) -- (3.445,1.065);
\draw (3.445,1.065) -- (4.14,1.35);

\draw [->] (4.14,1.35) -- (4.053,1.063);
\draw (4.14,1.35) -- (4.01,0.92);

\draw [->] (4.01,0.92) -- (4.103,0.647);
\draw (4.103,0.647)-- (4.15,0.51);

\draw [color=ffqqtt] (1.96,0.4) -- (2.4,0.88) -- (3.17,1.09) -- (3.75,0.92) -- (4.23,0.44) -- cycle;
\fill[color=ffqqtt,fill=ffqqtt,fill opacity=0.1] (1.96,0.4) -- (2.4,0.88) -- (3.17,1.09) -- (3.75,0.92) -- (4.23,0.44) -- cycle;

\draw [color=zzqqzz] (3.6,1.56) -- (3.17,1.09) -- (2.4,0.88) -- (1.82,1.05) -- (1.34,1.53) -- cycle;
\fill[color=zzqqzz,fill=zzqqzz,fill opacity=0.1] (3.6,1.56) -- (3.17,1.09) -- (2.4,0.88) -- (1.82,1.05) -- (1.34,1.53) -- cycle;

\draw [->] (4.15,0.51) -- (3.47,0.75);
\draw (3.47,0.75)-- (2.79,0.99);

\draw [->] (2.79,0.99) -- (3.035,1.285);
\draw (2.79,0.99)-- (3.28,1.56);
\end{tikzpicture}
         \caption{Fifth step with overlap}
         \label{fig_overlap5}
\end{subfigure}

\begin{subfigure}[b]{0.95\textwidth}
\center
\includegraphics[trim = 1cm 1cm 1cm 1cm, clip, angle=90, scale=0.9]{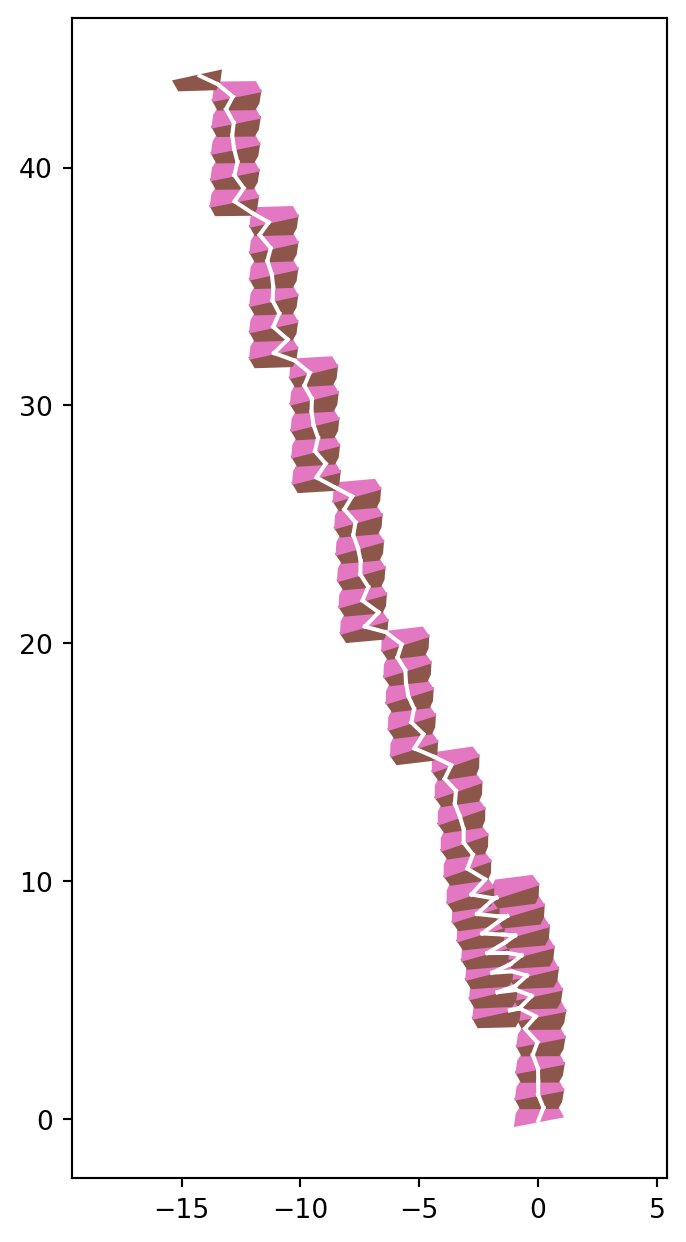}
\caption{A part of a generalized pentagon tiling billiard trajectory satisfying condition (C) }
\end{subfigure}
        \caption{Examples of trajectories in a generalized pentagonal tiling billiard}
        \label{fig:GTB}
\end{figure}

Fix a $N$-gon on the plane, a point inside it and a direction (see Figure~\ref{fig_overlap1}). Follow the direction until reaching a side. Then draw a second centrally symmetric polygon, so that the two polygons share the crossed side (see Figure~\ref{fig_overlap2penta}).
 The trajectory is refracted through the side and continues in a straight line in the second polygon (see Figure~\ref{fig_overlap2line}). Each time the trajectory reaches a side, we draw the following polygon via central symmetry. The new polygon may overlap a previous one. The trajectory takes into account only the last drawn polygon (see Figure~\ref{fig_overlap4}~and~\ref{fig_overlap5}).
A trajectory that hits a corner stops.

\subsubsection*{Assumptions}
\begin{figure}[h]
\begin{subfigure}[b]{0.4\textwidth}
\centering
\definecolor{qqzztt}{rgb}{0,0.6,0.2} 
\definecolor{ffzzqq}{rgb}{1,0.6,0} 
\definecolor{ffqqqq}{rgb}{1,0,0} 
\definecolor{zzttqq}{rgb}{0.2,0.2,0.2} 
\begin{tikzpicture}[line cap=round,line join=round,>=stealth,x=1.0cm,y=1.0cm]
\fill[color=zzttqq,fill=zzttqq,fill opacity=0.1] (0.85,0.53) -- (0.98,-0.22) -- (0.69,-0.72) -- (0.34,-0.94)-- (-0.32,-0.95) -- cycle;
\draw [color=zzttqq] (0.85,0.53) -- (0.98,-0.22) -- (0.69,-0.72) -- (0.34,-0.94)-- (-0.32,-0.95) -- cycle;

\draw [color=qqzztt]  plot[domain=-0.22:0.56,variable=\t]({1*1*cos(\t r)+0*1*sin(\t r)},{0*1*cos(\t r)+1*1*sin(\t r)});
\draw plot[domain=0.56:4.38,variable=\t]({1*1*cos(\t r)+0*1*sin(\t r)},{0*1*cos(\t r)+1*1*sin(\t r)});
\draw [color=ffqqqq]  plot[domain=4.38:5.06,variable=\t]({1*1*cos(\t r)+0*1*sin(\t r)},{0*1*cos(\t r)+1*1*sin(\t r)});
\draw [color=qqzztt]  plot[domain=5.06:5.48,variable=\t]({1*1*cos(\t r)+0*1*sin(\t r)},{0*1*cos(\t r)+1*1*sin(\t r)});
\draw [color=ffqqqq]  plot[domain=5.48:6.06,variable=\t]({1*1*cos(\t r)+0*1*sin(\t r)},{0*1*cos(\t r)+1*1*sin(\t r)});

\begin{scriptsize}
\draw[color=qqzztt] (1.2,0.26) node {$a_1$};
\draw[color=ffqqqq] (1.1,-0.5) node {$a_2$};
\draw[color=qqzztt] (0.69,-0.95) node {$a_4$};
\draw[color=ffqqqq] (0.05,-1.15) node {$...$};
\draw[color=black] (-0.9,0.74) node {$a_N$};
\end{scriptsize}
\end{tikzpicture}
\caption{Sides}
\label{fig:Notations_sides}
\end{subfigure} 
\hfill
\begin{subfigure}[b]{0.4\textwidth}
         \centering
\definecolor{ffffqq}{rgb}{1,1,0} 
\definecolor{ffzzqq}{rgb}{0.2,0.2,0.2} 
\definecolor{xdxdff}{rgb}{0.49,0.49,1}
\definecolor{ttccqq}{rgb}{0.2,0.8,0}
\definecolor{ffqqqq}{rgb}{1,0,0} 
\definecolor{uququq}{rgb}{0.25,0.25,0.25}
\begin{tikzpicture}[line cap=round,line join=round,>=stealth,x=1.0cm,y=1.0cm]

\draw [line width=1pt,color=black] (0.55,0.47) -- (-0.5,-0.19);
\draw [->] (0.55,0.47) -- (0.025,0.14);

\fill [fill=ffzzqq,fill opacity=0.1] (-0.96,-0.28) -- (-0.89,0.21) -- (-0.73,0.46) -- (0.8,0.47) -- (1,0.1) -- cycle;
\draw [ffzzqq] (-0.96,-0.28) -- (-0.89,0.21) -- (-0.73,0.46) -- (0.8,0.47) -- (1,0.1) -- cycle;

\draw [color=ffzzqq] (0.04,-0.18) circle (1cm);
\draw [dotted] (-0.5,-0.19)-- (-0.93,-0.46);
\draw [dotted] (0.55,0.47)-- (0.159+0.55,0.57);
\begin{scriptsize}
\draw [color=ffzzqq] (-0.96,-0.28)-- ++(-1.5pt,-1.5pt) -- ++(3.0pt,3.0pt) ++(-3.0pt,0) -- ++(3.0pt,-3.0pt);
\draw [color=ffzzqq] (-0.48,0.2) node {$\mathbf{P}$};
\draw [color=uququq] (-0.93,-0.46)-- ++(-1.5pt,-1.5pt) -- ++(3.0pt,3.0pt) ++(-3.0pt,0) -- ++(3.0pt,-3.0pt);
\draw [color=blue] (0.9,0.7) node {$\tau$};
\end{scriptsize}

\draw [->,>=stealth,color=blue] (-1.05,-0.53) arc (-164:50:1.145);
\end{tikzpicture}
\caption{Parameter $\tau$ of the trajectory}\label{Fig:def_tau}
\end{subfigure}
\caption{Notations for cyclic polygons}
\end{figure}
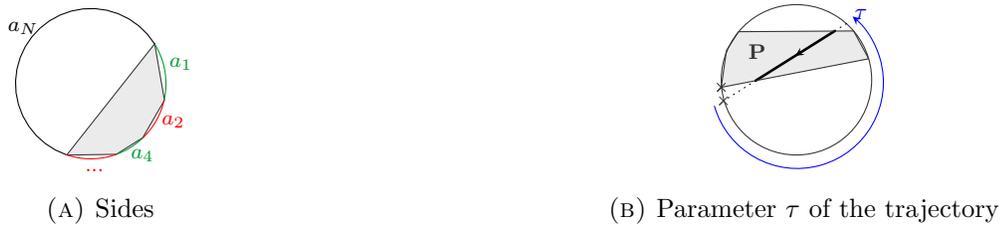

Let $P$ be a cyclic polygon.
We will denote by $\textbf{a}_\textbf{1},...,\textbf{a}_\textbf{N}$ the sides of $P$, clockwise labelled and $\textbf{a}_\textbf{N}$ being its longest side. We denote by $a_1,...,a_{N}$ the lengths of the arcs subtended by each side (see Figure~\ref{fig:Notations_sides}).
Conversely, if $a_1,\dots, a_{N}>0$, with $a_N$ being maximal, we denote by $P(a_1,\dots,a_{N})$ the cyclic polygon whose arcs subtended by its sides have lengths $a_1,\dots,a_{N}$.

We will restrict ourselves to the study of billiards satisfying the following conditions:
$$
(\text{C}) \left\{ \begin{array}{ll}
P \text{ is inscribed in a circle} & (\text{C}0) \\
a_1+..+a_{N-1}<\tau<a_{N} & (\text{C}1)
\end{array} \right.
$$

Condition (C$0$) ensures that a certain parameter $\tau$ (see Figure \ref{Fig:def_tau} and Definition~\ref{def_tau}) is preserved all along the trajectory and hence allows us to use IETFs to study the system (see Sections \ref{subsec_system-coord} and \ref{subsec_obtained-IET}), as in \cite{BSDFI18}, \cite{HPR22}, \cite{PR19}, \cite{PR21} and \cite{DHMPRS}.
This is a crucial fact for the now existing methods of study — we know nothing for tiling billiards made with non cyclic polygons.
Condition (C$1$) means both that $P$ is inscribed in a half-circle and that the trajectory is close to the circumcenter. It is crucial for the following since it allows to reduce the system to a standard IET without flips with $N$ intervals of continuity, which permits to use ergodic theory tools (see paragraphs~\ref{subsec_system-coord}~and~\ref{subsec_obtained-IET}).

\subsection{Results}
Let us identify the real plane with $\mathbb{C}$.
Let $p_n\in\mathbb{C}$ be the point of the trajectory on the $n$-th side crossing. 

\begin{restatable}{lem}{AsymptoticDir}
\label{asymptotic_dir}
For every $a_1, \dots, a_{N-1}$, almost every $a_N>\sum_{i=1}^{N-1} a_i$, and every $\tau\in \left(\sum_{i=1}^{N-1} a_i , a_N \right)$, the cyclic $N$-gon $P(a_1,\dots,a_N)$ defines a generalized tiling billiard whose trajectories having parameter $\tau$ all admit a mean displacement $m = \underset{n\rightarrow \infty}{\lim}\frac{p_n}{n}$ that does not vanish.


\end{restatable}

\begin{restatable}{theorem}{ThDev}
\label{th_dev}
Let $\mathcal{L} $ 
be a line in the asymptotic direction (almost always defined) of the trajectory.
For almost every tiling billiard satisfying the conditions (C),
the trajectory does not stay at a bounded distance from $\mathcal{L}$.
\end{restatable}

We quantify the order of magnitude of deviations from the asymptotic direction with the following two theorems. They are linked to the behavior of the Teichmüller flow on the moduli space of translation surfaces.

\begin{restatable}{theorem}{ThGeneric}
\label{th_generic}
Let $N\in \mathbb{N}, N\geqslant 5$. There exist $\theta_1 > \theta_2 > \theta_3 > 0 $ such that
for almost every 
$N$-gon $\mathbf{P}$ 
and every $\tau$,
both satisfying condition $(C)$, 
for any initial direction having parameter $\tau$ in a $\mathbf{P}$-tiling billiard,
the sequence of points $(p_n)_{n\in\NN}$ satisfies:    
$$\underset{n\rightarrow \infty}{\limsup} \frac{\log |p_{n} - nm|}{\log n} = \left\{ \begin{array}{l l}
\frac{\theta_2}{\theta_1} > 0 & \text{if } N \text{ is odd} \\
\frac{\theta_2}{\theta_1} \text{ or } \frac{\theta_3}{\theta_1} > 0  &  \text{if } N \text{ is even} \\
\end{array} \right., $$
where $ m =\underset{n\rightarrow\infty}{\lim} \frac{p_n}{n}$.
\end{restatable}

\begin{rmk}
The numbers $\theta_1$, $\theta_2$ and $\theta_3$ depend only on the number of sides of the polygon! The numbers $1+\frac{\theta_i}{\theta_1}$ ($1\leqslant i\leqslant 3$) are the first three Lyapounov exponents of Teichmüller flow in the hyperelliptic component of the stratum $\mathcal{H}(N-3)$ if $N$ is odd, and $\mathcal{H}(\frac{N}{2}-2,\frac{N}{2}-2)$ if $N$ is even.

When $N=5$, a result of Eskin, Kontsevich and Zorich \cite{EKZ} states that  $\frac{\theta_2}{\theta_1}=\frac{1}{3}$. Forni \cite{Forni} proves that $\theta_1, \theta_2$ and $\theta_3$ are distinct.
\end{rmk}

We prove these results via IETs and their generic properties. The \emph{Rauzy induction} $\mathfrak{R}$ (see paragraph \ref{subsec_Rauzy-cocycle-Rokhlin}) is a powerful tool to study orbits of an IET.
If $T$ is an IET on $I$, $\mathfrak{R}(T)$ is the first return map of $T$ on a well chosen subinterval $J\subset I$, which gives once again an IET with the same number of intervals of continuity. To track the orbit of $T$ through Rauzy induction, we need to know both the sequence of combinatorial data of the induced IETs $(\mathfrak{R}^n(T))_{n\in\NN}$ and the lengths of their intervals of continuity. This information can be encoded in an infinite path in a finite graph, called the Rauzy diagram.
The lengths of intervals of continuity of $\mathfrak{R}^n(T)$ linearly depend on those of $T$. A corresponding matrix cocycle (see paragraph \ref{subsec_Rauzy-cocycle-Rokhlin}) is important for the study of our dynamical system.
The Rauzy induction for IETs corresponds to the Teichmüller flow for their suspensions (that are translation surfaces). See \cite{Viana}, \cite{MMY},
\cite{SurveyYoccoz}.

When the path in the Rauzy diagram corresponding to some IET is a loop, which is the same as to say that the orbit of the IET through the Rauzy induction is periodic, we say that this IET is \textbf{self-similar}. The matrix associated to this loop allows us to study the symbolic coding of a point by the IET $T$. When the loop is complex enough, the matrix is primitive. 
This case is obviously not generic but important. Indeed it corresponds to the case when the vertical flow on a suspension of the IET is a pseudo-Anosov map.
In this case, one can give a more precise estimate of the deviations.

\begin{figure}[h!]
\includegraphics[scale=0.5]{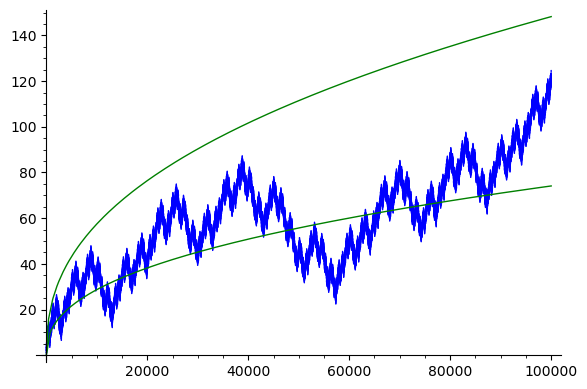}
\caption{The modulus of the deviations, depending on the number of steps of the trajectory, has a sublinear growth}
\begin{footnotesize}\begin{flushleft}
Deviations, depicted in blue, are always smaller than a (strictly) sublinear function, depicted in green, with a growth rate $\rho$, and infinitely often greater than another function, also in green, with the same growth rate $\rho$.
\end{flushleft}
\end{footnotesize}
\label{fig_graphe_dev}
\end{figure}

\begin{restatable}{theorem}{ThSelfsim}
\label{th_selfsim}
Let $N\in \mathbb{N}, N\geqslant 5$ and $a=(a_1,a_2,...,a_{N-1})$ be the (normalized) Perron eigenvector (with eigenvalue $\lambda_1$) of a primitive matrix $M\in\mathcal{M}_{N-1}(\mathbb{R})$ that corresponds to a loop in the Rauzy diagram. 
Assume that the second eigenvalue $\lambda_2$ of $M$ is such that $|\lambda_2|>1$ and that no other eigenvalue has the same modulus as $\lambda_2$.
For almost every $a_{N}>\sum_{i=1}^{N-1}a_i$,
%
%
there exist $C_1$ and $C_2$ and an infinite set $\mathcal{S} \subset \mathbb{N}$ such that
for every $\tau \in (\sum_{i=1}^{N-1}a_i, a_N)$ 
and any trajectory having parameter $\tau$ in a $P(a_1,...,a_N)$-tiling billiard,
the sequence of points $(p_n)_{n\in\NN}$ satisfies: 
$$\left\{ \begin{array}{cc}
     \forall n \in \mathbb{N} & |p_n - nm| < C_1 n^{\frac{\log\lambda_2}{\log\lambda_1}}\\
     \forall n \in \mathcal{S} &|p_n - nm| > C_2 n^{\frac{\log\lambda_2}{\log\lambda_1}}
\end{array}
\right.
.
$$
\end{restatable}
Figure \ref{fig_graphe_dev} illustrates this phenomenon of deviation and the growth rate.
This theorem quantifies the second term of the asymptotic expansion of some Birkhoff's sums for a fixed point. Another related problem is to study the uniform speed of convergence of the same Birkhoff's sums, i.e. to study the \emph{discrepancy}. This was done by Adamczewski in \cite{Adamczewski}.


To summarize, this article shows that for almost every generalized tiling billiard satisfying condition (C):
\begin{enumerate}
    \item like for triangle or cyclic quadrilateral tiling billiards,
the trajectory has an asymptotic direction (Lemma~\ref{asymptotic_dir}),
    \item unlike in triangles and cyclic quadrilaterals case, it deviates sublinearly from it (Theorems~\ref{th_dev},~\ref{th_generic}~and~\ref{th_selfsim}).
\end{enumerate}



\subsection{Outline of the proof}
The proofs of our theorems are organised as follows:
\begin{enumerate}
    \item Define the IET that corresponds to our generalized tiling billiard.
    \item Define a piecewise constant function $f$ (that depends on $P$) and approximates the trajectory. We show that $f$ has almost always a non zero mean, which means that the trajectory of the studied tiling billiard has an asymptotic direction. 
    \item Study the function $h = \Im(\frac{1}{m}f - 1) $ linked to the deviations from the asymptotic direction. We show that for almost every $a_N$, $h$ is not contained in a codimension 3 subspace (resp. 2 if $N$ is odd).
    \item Zorich's results on deviations for IETs give the rate given in  Theorem~\ref{th_generic} 
\end{enumerate}
IETs are generically uniquely ergodic \cite{Masur} \cite{Veech}. Hence we know that the orbit of a point equidistributes in the different intervals of continuity. Zorich establishes \cite{howdo} how much an orbit can deviate from this equidistribution. We can think of such results as central limit theorems. To prove this, Zorich studies Teichmüller flow of suspension surfaces (that are translation surfaces) of IETs and how it behaves with respect to Oseledets' flag decomposition.

The main contribution of this article is to prove that one can apply Zorich's results to the case of tiling billiards we study. 
This is not a straightforward verification. We have to show where some specific vector lies in the flag given by the Oseledets' theorem (see Section \ref{subsec_Oseledets}). The difficult point is that we know the expression of the vector but, in general, not of the subspaces of the flag. We can think of them as random subspaces. We can still conclude for almost every parameters because the vector has one more degree of freedom than the flag. Moreover it depends analytically on this parameter independent of the flag. This corresponds to step (3) in the following decomposition of the proof:

The most technical parts of the proofs of Theorems~\ref{th_generic} are done first for pentagons and then generalized.
We prove Theorem~\ref{th_selfsim} by direct computation.

\subsection{Outline of the paper}
We explain step (1) in Section~\ref{section_studied_system}. Then we recall briefly some background in Section~\ref{background}. We introduce functions $f$ and $g$ in Section~\ref{section_approx_traj}. Section~\ref{section_G} studies the function $g$ and show the properties that allow us to apply Zorich's results.  
Section ~\ref{section_selfsim} handles with the direct computations of deviations in case of a self-similar IET, using techniques of symbolic coding, Rauzy induction and Rokhlin towers, in order to prove Theorem~\ref{th_selfsim}.

\vspace{0.5cm}
\textbf{Acknowledgements}

I thank my PhD advisors Pascal Hubert and Olga Paris-Romaskevich for having introduced this subject to me. 

\section{The studied system}\label{section_studied_system} 
\subsection{Definition of a generalized tiling billiard}
\begin{figure}
\definecolor{ffffqq}{rgb}{1,1,0} 
\definecolor{ffzzqq}{rgb}{1,0.6,0} 
\definecolor{xdxdff}{rgb}{0.49,0.49,1}
\definecolor{ttccqq}{rgb}{0.2,0.8,0}
\definecolor{ffqqqq}{rgb}{1,0,0} 
\definecolor{uququq}{rgb}{0.25,0.25,0.25}
\begin{tikzpicture}[line cap=round,line join=round,>=stealth,x=2.0cm,y=2.0cm]

\fill[fill=ffzzqq,fill opacity=0.1] (-0.96,-0.28) -- (-0.89,0.21) -- (-0.73,0.46) -- (0.8,0.47) -- (1,0.1) -- cycle;
\draw [color=ffzzqq] (-0.96,-0.28)-- (-0.89,0.21);
\draw [color=ffzzqq] (-0.89,0.21)-- (-0.73,0.46);
\draw [color=ffzzqq] (-0.73,0.46)-- (0.8,0.47);
\draw [color=ffzzqq] (0.8,0.47)-- (1,0.1);
\draw [color=ffzzqq] (1,0.1)-- (-0.96,-0.28);

\begin{scriptsize}
\draw [color=ffzzqq] (-0.96,-0.28)-- ++(-1.5pt,-1.5pt) -- ++(3.0pt,3.0pt) ++(-3.0pt,0) -- ++(3.0pt,-3.0pt);
\draw [color=ffzzqq] (-1.08,-0.25) node {$v_0$};
\draw [color=ffzzqq] (-0.48,0.2) node {$P_0$};
\end{scriptsize}

\fill[color=ffqqqq,fill=ffqqqq,fill opacity=0.1] (1,0.1) -- (0.92,-0.39) -- (0.77,-0.64) -- (-0.76,-0.65) -- (-0.96,-0.28) -- cycle;
\draw [color=ffqqqq] (1,0.1)-- (0.92,-0.39);
\draw [color=ffqqqq] (0.92,-0.39)-- (0.77,-0.64);
\draw [color=ffqqqq] (0.77,-0.64)-- (-0.76,-0.65);
\draw [color=ffqqqq] (-0.76,-0.65)-- (-0.96,-0.28);
\draw [color=ffqqqq] (-0.96,-0.28)-- (1,0.1);

\fill[color=ttccqq,fill=ttccqq,fill opacity=0.1] (0.69,-1.13) -- (0.77,-0.64) -- (0.92,-0.39) -- (2.45,-0.38) -- (2.65,-0.75) -- cycle;
\draw [color=ttccqq] (0.69,-1.13)-- (0.77,-0.64);
\draw [color=ttccqq] (0.77,-0.64)-- (0.92,-0.39);
\draw [color=ttccqq] (0.92,-0.39)-- (2.45,-0.38);
\draw [color=ttccqq] (2.45,-0.38)-- (2.65,-0.75);
\draw [color=ttccqq] (2.65,-0.75)-- (0.69,-1.13);
\begin{scriptsize}
\draw [color=ttccqq] (0.69,-1.13)-- ++(-1.5pt,-1.5pt) -- ++(3.0pt,3.0pt) ++(-3.0pt,0) -- ++(3.0pt,-3.0pt);
\draw[color=ttccqq] (0.58,-1.11) node {$v_2$};
\draw[color=ttccqq] (1.67,-0.61) node {$P_2$};
\end{scriptsize}

\fill[color=blue,fill=blue,fill opacity=0.1] (2.65,-0.75) -- (2.58,-1.23) -- (2.42,-1.49) -- (0.89,-1.49) -- (0.69,-1.13) -- cycle;
\draw [color=blue] (2.65,-0.75)-- (2.58,-1.23);
\draw [color=blue] (2.58,-1.23)-- (2.42,-1.49);
\draw [color=blue] (2.42,-1.49)-- (0.89,-1.49);
\draw [color=blue] (0.89,-1.49)-- (0.69,-1.13);
\draw [color=blue] (0.69,-1.13)-- (2.65,-0.75);

\begin{scriptsize}
\draw [color=blue] (2.65,-0.75)-- ++(-1.5pt,-1.5pt) -- ++(3.0pt,3.0pt) ++(-3.0pt,0) -- ++(3.0pt,-3.0pt);
\draw[color=blue] (2.74,-0.66) node {$v_3$};
\draw[color=blue] (1.94,-1.32) node {$P_3$};
\end{scriptsize}

\draw [line width=1pt,color=black,->] (0.55,0.47) -- (0.025,0.14);
\draw [line width=1pt,color=black] (0.025,0.14) -- (-0.5,-0.19);
\draw [line width=1pt,color=black,->] (-0.5,-0.19) -- (0.195,-0.315);
\draw [line width=1pt,color=black] (0.195,-0.315) -- (0.89,-0.44);
\draw [line width=1pt,color=black,->] (0.89,-0.44) -- (1.105,-0.72);
\draw [line width=1pt,color=black] (1.105,-0.72) -- (1.32,-1);
\draw [line width=1pt,color=black,->] (1.32,-1) -- (1.255,-1.245);
\draw [line width=1pt,color=black] (1.255,-1.245) -- (1.19,-1.49);

\begin{scriptsize}
\draw [color=ffqqqq] (1,0.1)-- ++(-1.5pt,-1.5pt) -- ++(3.0pt,3.0pt) ++(-3.0pt,0) -- ++(3.0pt,-3.0pt);
\draw[color=ffqqqq] (1.09,0.19) node {$v_1$};
\draw[color=ffqqqq] (-0.09,-0.5) node {$P_1$};

\draw [color=black] (0.5,0.25) node {$\mathbf{T}$};
\end{scriptsize}
\end{tikzpicture}
\caption{Notations for generalized tiling billiards}\label{fig:Notations_GTB}
\end{figure}
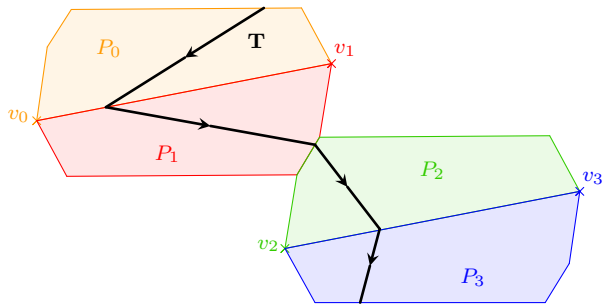

\begin{definition}[Generalized tiling billiard]
Let $P$ be a polygon and $v$ a vertex of $P$. We call a 
generalized tiling billiard trajectory the triplet $(\mathbf{T}, (P_i)_{i\in\NN}, (v_i)_{i\in\NN})$ where $\mathbf{T}$ is a piecewise linear curve in the real plane such that for every integer $i$:
\begin{itemize}
    \item $(P_i,v_i)$ is an isometric copy of $(P,v)$ in the plane,
    \item $\mathbf{T}$ crosses $P_i$ in a segment,
    \item $P_i$ and $P_{i+1}$ share a side, crossed by $\mathbf{T}$,
    \item $P_{i+1}$ is the image of $P_i$ through a central symmetry,
    \item the common side of $P_i$ and $P_{i+1}$ is the bissector of the angle made by $\mathbf{T}\cap P_i$ and $\mathbf{T}\cap P_{i+1}$. 
    \end{itemize} 
\end{definition}
This defines a $P$-generalized tiling billiard as a dynamical system.
See Figure~\ref{fig:Notations_GTB} for an illustration.

\begin{rmk}
When we think of a generalized tiling billiard trajectory, we mostly think of $\mathbf{T}$ although defining it as a triplet with crossed polygons and vertices $v_i$ is much more convenient for the study of the trajectory. The definition requires no property of $P$ although we will restrict in the following to the case where $P$ is cyclic. We denote by $\mathcal{C}_i$ the circumcircle of $P_i$. 
\end{rmk}

\subsection{A system of coordinates}\label{subsec_system-coord}
Here we describe a convenient frame of reference to study a tiling billiard trajectory. It is the one of an ant that would follow this trajectory. At each step, the ant can see the polygon it is crossing and its circumcircle. 
We need the following data to determine the position of $P_i$ with respect to the ant (in other words, the position of the ant in $P$): 
\begin{enumerate}
\item
On which chord of the circumcircle $\mathcal{C}_i$ of $P_i$ the ant is travelling,
\item
Where $v_i$ is on $\mathcal{C}_i$ (this determines how $P_i$ is placed in $\mathcal{C}_i$).
\end{enumerate}

Lemma 3.2 in \cite{BSDFI18} applies also for generalized tiling billiards made with any cyclic polygon and shows that the chord of $\mathcal{C}_i$ on which the ant is travelling is preserved (See Section \ref{subsec_obtained-IET}).
The data (1) is hence a parameter of the trajectory.
%
For example if the ant is travelling on a diameter of $\mathcal{C}_i$, then it will also travel on a diameter of any $\mathcal{C}_{j}$. 
Therefore we will need only the position of $v_i$ on $\mathcal{C}_i$ to determine the position of $P$ with respect to the ant.

\begin{definition}\label{def_tau}
We say that the trajectory has parameter $\tau$ 
if the arc between the two endpoints of the chord supporting a step of the trajectory $\mathbf{T}$ has length $\tau$ (see Figure~\ref{Fig:def_tau}).
\end{definition}

\begin{figure}[h!]
     \centering
\begin{subfigure}[b]{0.8\textwidth}
         \centering
\definecolor{ffffqq}{rgb}{1,1,0} 
\definecolor{ffzzqq}{rgb}{1,0.6,0} 
\definecolor{xdxdff}{rgb}{0.49,0.49,1}
\definecolor{ttccqq}{rgb}{0.2,0.8,0}
\definecolor{ffqqqq}{rgb}{1,0,0} 
\definecolor{uququq}{rgb}{0.25,0.25,0.25}
\begin{tikzpicture}[line cap=round,line join=round,>=stealth,x=2.0cm,y=2.0cm]

\draw [->,line width=1pt] (0.55,0.47) -- (0.025,0.14);
\draw [line width=1pt] (0.025,0.14) -- (-0.5,-0.19);

\draw [line width=1pt,color=black] (-0.5,-0.19) -- (0.89,-0.44) -- (1.32,-1) -- (1.19,-1.49);

\fill[fill=ffzzqq,fill opacity=0.1] (-0.96,-0.28) -- (-0.89,0.21) -- (-0.73,0.46) -- (0.8,0.47) -- (1,0.1) -- cycle;
\draw [color=ffzzqq] (-0.96,-0.28) -- (-0.89,0.21) -- (-0.73,0.46) -- (0.8,0.47) -- (1,0.1) -- cycle;

\draw [color=ffzzqq] (0.04,-0.18) circle (2cm);
\draw [color=ffzzqq, dotted] (-0.5,-0.19)-- (-0.93,-0.46);
\begin{scriptsize}
\draw [color=ffzzqq] (-0.96,-0.28)-- ++(-1.5pt,-1.5pt) -- ++(3.0pt,3.0pt) ++(-3.0pt,0) -- ++(3.0pt,-3.0pt);
\draw [color=ffzzqq] (-1.08,-0.25) node {$v_0$};
\draw [color=ffzzqq] (-0.48,0.2) node {$P_0$};
\draw [color=ffzzqq] (-0.01,0.7) node {$C_0$};
\draw [color=ffzzqq] (-0.93,-0.46)-- ++(-1.5pt,-1.5pt) -- ++(3.0pt,3.0pt) ++(-3.0pt,0) -- ++(3.0pt,-3.0pt); 
\draw [color=ffzzqq] (-1,-0.55) node {$o_0$}; 
\end{scriptsize}

\fill [color=ffqqqq,fill=ffqqqq,fill opacity=0.1] (1,0.1) -- (0.92,-0.39) -- (0.77,-0.64) -- (-0.76,-0.65) -- (-0.96,-0.28) -- cycle;
\draw [color=ffqqqq] (1,0.1)-- (0.92,-0.39) -- (0.77,-0.64) -- (-0.76,-0.65) -- (-0.96,-0.28) -- cycle;

\draw [color=ffqqqq] (0,0) circle (2cm);
\draw [dotted,color=ffqqqq] (0.89,-0.44)-- (0.9,-0.44);

\fill [color=ttccqq,fill=ttccqq,fill opacity=0.1] (0.69,-1.13) -- (0.77,-0.64) -- (0.92,-0.39) -- (2.45,-0.38) -- (2.65,-0.75) -- cycle;
\draw [color=ttccqq] (0.69,-1.13)-- (0.77,-0.64) -- (0.92,-0.39) -- (2.45,-0.38) -- (2.65,-0.75) -- cycle;

\draw [color=ttccqq] (1.69,-1.03) circle (2cm);
\draw [dotted,color=ttccqq] (1.32,-1)-- (2.05,-1.96);
\begin{scriptsize}
\draw [color=ttccqq] (0.69,-1.13)-- ++(-1.5pt,-1.5pt) -- ++(3.0pt,3.0pt) ++(-3.0pt,0) -- ++(3.0pt,-3.0pt);
\draw [color=ttccqq] (0.58,-1.11) node {$v_2$};
\draw [color=ttccqq] (1.67,-0.61) node {$P_2$};
\draw [color=ttccqq] (1.68,-0.12) node {$C_2$};
\draw [color=ttccqq] (2.05,-1.96)-- ++(-1.5pt,-1.5pt) -- ++(3.0pt,3.0pt) ++(-3.0pt,0) -- ++(3.0pt,-3.0pt);
\draw [color=ttccqq] (2.17,-2) node {$o_2$};
\end{scriptsize}

\fill [color=blue,fill=blue,fill opacity=0.1] (2.65,-0.75) -- (2.58,-1.23) -- (2.42,-1.49) -- (0.89,-1.49) -- (0.69,-1.13) -- cycle;
\draw [color=blue] (2.65,-0.75) -- (2.58,-1.23) -- (2.42,-1.49) -- (0.89,-1.49) -- (0.69,-1.13) -- cycle;

\draw [color=blue] (1.65,-0.85) circle (2cm);
\draw [dotted,color=blue] (1.19,-1.49)-- (1.13,-1.7);
\begin{scriptsize}
\draw [color=blue] (2.65,-0.75)-- ++(-1.5pt,-1.5pt) -- ++(3.0pt,3.0pt) ++(-3.0pt,0) -- ++(3.0pt,-3.0pt);
\draw [color=blue] (2.74,-0.66) node {$v_3$};
\draw [color=blue] (1.94,-1.32) node {$P_3$};
\draw [color=blue] (1.63,0.25) node {$C_3$};
\draw [color=blue] (1.13,-1.7)-- ++(-1.5pt,-1.5pt) -- ++(3.0pt,3.0pt) ++(-3.0pt,0) -- ++(3.0pt,-3.0pt);
\draw [color=blue] (1.25,-1.65) node {$o_3$};
\end{scriptsize}

\begin{scriptsize}
\draw [color=ffqqqq] (-0.9,0.69) node {$C_1$};
\draw [color=ffqqqq] (1,0.1)-- ++(-1.5pt,-1.5pt) -- ++(3.0pt,3.0pt) ++(-3.0pt,0) -- ++(3.0pt,-3.0pt);
\draw [color=ffqqqq] (1.09,0.19) node {$v_1$};
\draw [color=ffqqqq] (-0.09,-0.5) node {$P_1$};
\draw [color=ffqqqq] (0.9,-0.44)-- ++(-1.5pt,-1.5pt) -- ++(3.0pt,3.0pt) ++(-3.0pt,0) -- ++(3.0pt,-3.0pt);
\draw [color=ffqqqq] (0.8,-0.5) node {$o_1$};

\draw [color=black] (0.5,0.25) node {$\mathbf{T}$};
\end{scriptsize}



\end{tikzpicture}
\caption{Notations}      
\end{subfigure}

\begin{subfigure}[b]{0.4\textwidth}
         \centering
\definecolor{ffffqq}{rgb}{1,1,0} 
\definecolor{ffzzqq}{rgb}{1,0.6,0} 
\definecolor{xdxdff}{rgb}{0.49,0.49,1}
\definecolor{ttccqq}{rgb}{0.2,0.8,0}
\definecolor{ffqqqq}{rgb}{1,0,0} 
\definecolor{uququq}{rgb}{0.25,0.25,0.25}
\begin{tikzpicture}[rotate = 164, line cap=round,line join=round,>=stealth,x=2.0cm,y=2.0cm]

\draw [color=ffzzqq] (0.04,-0.18) circle (2cm);
\draw (0.04,-0.18) node { };
\draw [dotted] (-0.5,-0.19)-- (-0.93,-0.46);

\draw [->,line width=1pt] (0.55,0.47) -- (0.025,0.14);
\draw [line width=1pt] (0.025,0.14) -- (-0.5,-0.19);

\fill [fill=ffzzqq,fill opacity=0.1] (-0.96,-0.28) -- (-0.89,0.21) -- (-0.73,0.46) -- (0.8,0.47) -- (1,0.1) -- cycle;
\draw [color=ffzzqq] (-0.96,-0.28)-- (-0.89,0.21) -- (-0.73,0.46) -- (0.8,0.47) -- (1,0.1) -- cycle;

\begin{scriptsize}
\draw [color=ffzzqq] (-0.96,-0.28)-- ++(-1.5pt,-1.5pt) -- ++(3.0pt,3.0pt) ++(-3.0pt,0) -- ++(3.0pt,-3.0pt);
\draw [color=ffzzqq] (-0.825,-0.15) node {$v_0$};
\draw [color=ffzzqq] (-0.48,0.2) node {$P_0$};
\draw [color=ffzzqq] (-0.01,0.7) node {$C_0$};
\draw [color=ffzzqq] (-0.93,-0.46)-- ++(-1.5pt,-1.5pt) -- ++(3.0pt,3.0pt) ++(-3.0pt,0) -- ++(3.0pt,-3.0pt);
\draw [color=ffzzqq] (-0.8,-0.5) node {$o_0$};
\draw  (-1.2,-0.15) node {$x_0$}; 
\end{scriptsize}

\draw [ dashed,->] (-1.02,-0.5) arc (-164:184:1.1); 
\end{tikzpicture}
\caption{First crossed polygon}
\end{subfigure}
\hfill
\begin{subfigure}[b]{0.4\textwidth}
         \centering
\definecolor{ffffqq}{rgb}{1,1,0} 
\definecolor{xdxdff}{rgb}{0.49,0.49,1}
\definecolor{ttccqq}{rgb}{0.2,0.8,0}
\definecolor{ffqqqq}{rgb}{1,0,0} 
\definecolor{uququq}{rgb}{0.25,0.25,0.25}
\begin{tikzpicture}[rotate=24, line cap=round,line join=round,>=stealth,x=2.0cm,y=2.0cm]

\draw [line width=1pt,->] (-0.5,-0.19) -- (0.195,-0.315);
\draw [line width=1pt] (0.195,-0.315) -- (0.89,-0.44);

\fill[color=ffqqqq,fill=ffqqqq,fill opacity=0.1] (1,0.1) -- (0.92,-0.39) -- (0.77,-0.64) -- (-0.76,-0.65) -- (-0.96,-0.28) -- cycle;
\draw [color=ffqqqq] (1,0.1)-- (0.92,-0.39) -- (0.77,-0.64) -- (-0.76,-0.65) -- (-0.96,-0.28) -- cycle;

\draw [color=ffqqqq] (0,0) circle (2cm); 
\draw [dotted,color=ffqqqq] (0.89,-0.44)-- (0.9,-0.44);
\begin{scriptsize}
\draw [color=ffqqqq] (-0.9,0.69) node {$C_1$};
\draw [color=ffqqqq] (1,0.1)-- ++(-1.5pt,-1.5pt) -- ++(3.0pt,3.0pt) ++(-3.0pt,0) -- ++(3.0pt,-3.0pt);
\draw [color=ffqqqq] (1.09,0.19) node {$v_1$};
\draw [color=ffqqqq] (-0.09,-0.5) node {$P_1$};
\draw [color=ffqqqq] (0.9,-0.44)-- ++(-1.5pt,-1.5pt) -- ++(3.0pt,3.0pt) ++(-3.0pt,0) -- ++(3.0pt,-3.0pt);
\draw [color=ffqqqq] (1.05,-0.55) node {$o_1$};
\draw (1.25,0) node {$x_1$};
\end{scriptsize}

\draw [dashed,->] (1,-0.45) arc (-24:7:1.1);
\end{tikzpicture}
\caption{Second crossed polygon}
\end{subfigure}

\caption{Moving frame of reference}
\label{fig:Notations_frame}
\end{figure}
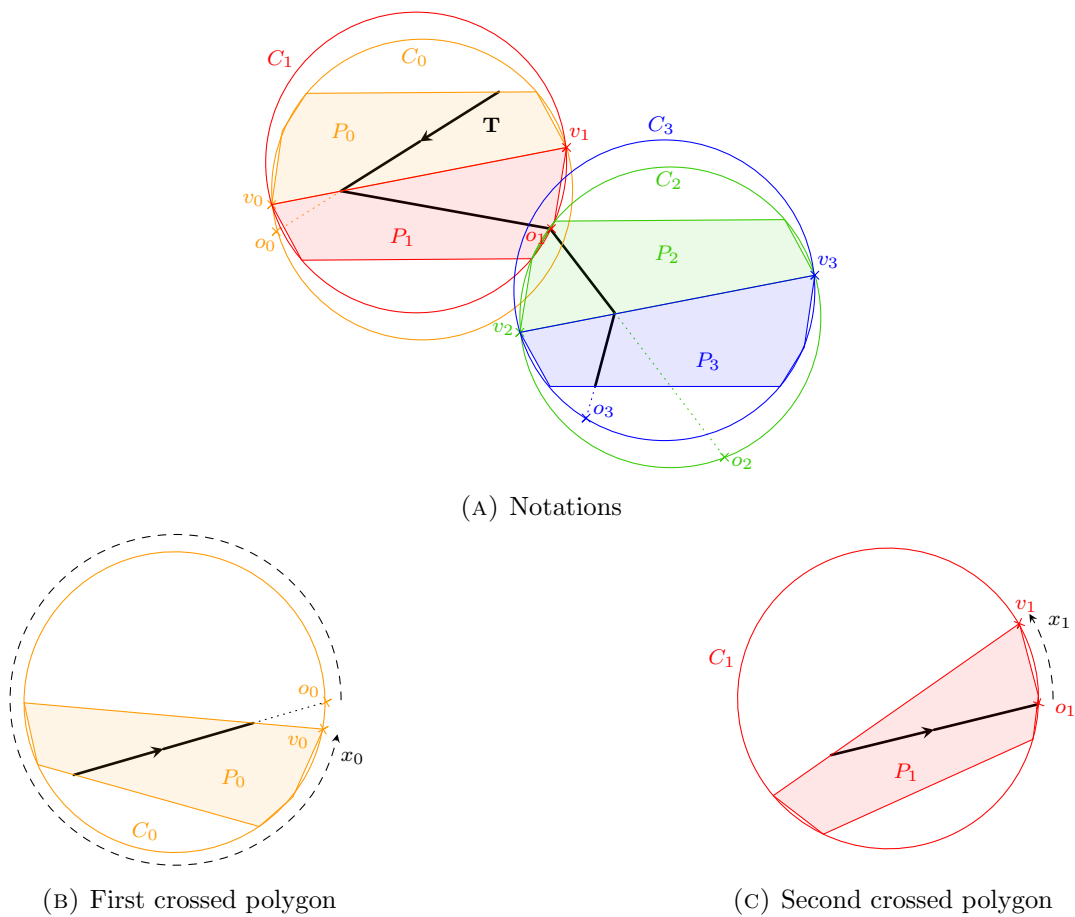

Let $(\mathbf{T}, (P_i)_{i\in\NN}, (v_i)_{i\in\NN})$ be a tiling billiard. 
We normalize the situation to $a_1~+~\dots~+~a_{N-1}~=~1$. 
We identify the circumcircle $\mathcal{C}_i$ of $P_i$ with the interval $[0,\sum_{k=1}^{N} a_{k})=[0,1+a_{N})$ by identifying 0 to the point $o_i\in\mathcal{C}_i$ where the line containing $T_i = \mathbf{T}\cap P_i$ would meet $\mathcal{C}_i$ if going on in straight line after the side crossing.
Then we locate the vertex $v_i$ of $P_i$ 
in the frame of $\mathcal{C}_i$. Let $x_i\in[0,1+a_{N})$ be the point with which $v_i$ is identified.
See Figure~\ref{fig:Notations_frame} for an illustration of notations.

\subsection{The obtained Interval Exchange Transformation} \label{subsec_obtained-IET}  

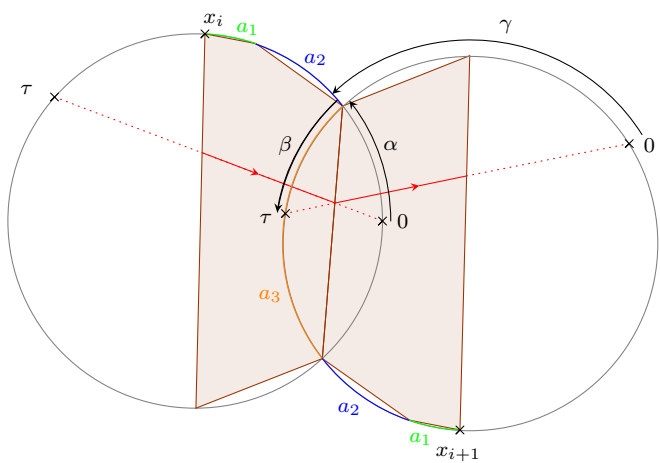
\begin{figure}
\definecolor{ffqqqq}{rgb}{1,0,0}
\definecolor{ffzzqq}{rgb}{1,0.6,0}
\definecolor{zzttqq}{rgb}{0.6,0.2,0}
\definecolor{xdxdff}{rgb}{0.49,0.49,1}
\begin{tikzpicture}[line cap=round,line join=round,>=stealth,x=1.0cm,y=1.0cm]
\fill[color=zzttqq,fill=zzttqq,fill opacity=0.1] (1.69,-1.83) -- (1.96,1.53) -- (0.81,2.36) -- (0.13,2.49) -- (0.01,-2.49) -- cycle;
\fill[color=zzttqq,fill=zzttqq,fill opacity=0.1] (1.96,1.53) -- (1.69,-1.83) -- (2.85,-2.65) -- (3.52,-2.78) -- (3.65,2.2) -- cycle;
\draw [color=gray] (0,0) circle (2.49cm);
\draw [color=zzttqq] (1.69,-1.83)-- (1.96,1.53);
\draw [color=zzttqq] (1.96,1.53)-- (0.81,2.36);
\draw [color=zzttqq] (0.81,2.36)-- (0.13,2.49);
\draw [color=zzttqq] (0.13,2.49)-- (0.01,-2.49);
\draw [color=zzttqq] (0.01,-2.49)-- (1.69,-1.83);
\draw [dotted,color=ffqqqq] (-1.87,1.65)-- (2.49,0);
\draw [->,color=ffqqqq] (0.09,0.91)-- (0.85,0.62);
\draw [-,color=ffqqqq] (0.85,0.62)-- (1.86,0.24);
\draw [color=zzttqq] (1.96,1.53)-- (1.69,-1.83);
\draw [color=zzttqq] (1.69,-1.83)-- (2.85,-2.65);
\draw [color=zzttqq] (2.85,-2.65)-- (3.52,-2.78);
\draw [color=zzttqq] (3.52,-2.78)-- (3.65,2.2);
\draw [color=zzttqq] (3.65,2.2)-- (1.96,1.53);
\draw [color=gray] (3.66,-0.3) circle (2.49cm);
\draw [dotted,color=ffqqqq] (5.77,1.03)-- (3.61,0.6);
\draw [dotted,color=ffqqqq] (1.86,0.24)-- (1.2,0.1);
\draw [->,color=ffqqqq] (1.86,0.24)-- (2.98,0.47);
\draw [color=ffqqqq] (2.98,0.47)-- (3.61,0.6);
\draw [->] (2.6,0) arc (0:37.98:2.6);
\draw [shift={(3.66,-0.3)}] plot[domain=2.32:2.98,variable=\t]({1*2.6*cos(\t r)+0*2.6*sin(\t r)},{0*2.6*cos(\t r)+1*2.6*sin(\t r)});
\draw [->] (5.94,1.15) arc (32.25:132.78:2.7);
\draw [->] (1.89,1.61) arc (132.78:170.76:2.6);
\draw [color=blue] (1.96,1.53) arc (37.98:71.02:2.49);
\draw [color=green] (0.81,2.36) arc (71.02:86.93:2.49);

\draw [color=orange] (1.96,1.53) arc (132.78:217.98:2.49);
\draw [color=blue] (1.69,-1.83) arc (217.98:251.02:2.49);
\draw [color=green] (2.85,-2.65) arc (251.02:266.17:2.49);
\begin{scriptsize}
\draw [color=black] (0.13,2.49)-- ++(-1.5pt,-1.5pt) -- ++(3.0pt,3.0pt) ++(-3.0pt,0) -- ++(3.0pt,-3.0pt);
\draw (0.25,2.68) node {$x_i$};
\draw [color=black] (-1.87,1.65)-- ++(-1.5pt,-1.5pt) -- ++(3.0pt,3.0pt) ++(-3.0pt,0) -- ++(3.0pt,-3.0pt);
\draw[color=black] (-2.24,1.75) node {$\tau$};
\draw [color=black] (3.52,-2.78)-- ++(-1.5pt,-1.5pt) -- ++(3.0pt,3.0pt) ++(-3.0pt,0) -- ++(3.0pt,-3.0pt);
\draw (3.5,-3.1) node {$x_{i+1}$}; 
\draw [color=black] (5.77,1.03)-- ++(-1.5pt,-1.5pt) -- ++(3.0pt,3.0pt) ++(-3.0pt,0) -- ++(3.0pt,-3.0pt);
\draw[color=black] (6.04,1.01) node {$0$};
\draw [color=black] (2.49,0)-- ++(-1.5pt,-1.5pt) -- ++(3.0pt,3.0pt) ++(-3.0pt,0) -- ++(3.0pt,-3.0pt);
\draw[color=black] (2.76,0) node {$0$};
\draw [color=black] (1.2,0.1)-- ++(-1.5pt,-1.5pt) -- ++(3.0pt,3.0pt) ++(-3.0pt,0) -- ++(3.0pt,-3.0pt);
\draw[color=black] (0.94,0.03) node {$\tau$};
\draw[color=black] (2.6,1) node {$\alpha$}; 
\draw[color=black] (1.2,1) node {$\beta$}; 
\draw[color=black] (4.13,2.59) node {$\gamma$};
\draw[color=blue] (1.6,2.15) node {$a_2$};
\draw[color=green] (0.7,2.55) node {$a_1$};
\draw[color=orange] (1,-1) node {$a_3$};
\draw[color=blue] (2.05,-2.5) node {$a_2$};
\draw[color=green] (3,-2.9) node {$a_1$};
\end{scriptsize}
\end{tikzpicture}
\caption{Coordinates $x_i$ and $x_{i+1}$ are related via a linear equation}\label{fig:computation_IET}
\end{figure}

We can prove the following by direct computation (see Figure \ref{fig:computation_IET}) or with an argument of folding, which was the proof of \cite{BSDFI18} for triangles. It also applies for generalized tiling billiards made with any cyclic polygon, see Lemma 1 of \cite{PR19}.

\begin{thm}[Theorem 3.3 of \cite{BSDFI18}]\label{Th_IETF}
There is an IETF $\Phi$ such that $x_{i+1}=\Phi(x_i)$  for every $i$ (see Figure~\ref{fig:Phi}).
\end{thm}

If we add the hypothesis (C1), 
the IETF $\Phi$ of Theorem~\ref{Th_IETF} has permutation 
$$
\pi(\Phi)= \left(
\begin{array}{cccccc}
\overline{I_1}&\overline{I_2}&...&\overline{I_{N-1}}&\overline{I_{N}^+}&\overline{I_{N}^-}\\
\overline{I_{N}^+}&\overline{I_1}&\overline{I_2}&...&\overline{I_{N-1}}&\overline{I_{N}^-}
\end{array}
\right)
$$ 
where we overlined flipped intervals (here, the IET flipps every interval). The intervals of continuity of the IET have lengths
$$\lambda(\Phi)= \left(
\begin{array}{cccccc}
a_1&a_2&...&a_{N-1}&a_{N} - \tau&\tau
\end{array}
\right).
$$ 

\begin{rmk}
Here we consider $\Phi$ as an IET on an interval, it has $N+1$ intervals of continuity. If we define it on a circle it has only $N$ intervals of continuity.

The choice of $v$ ensures that the trajectory crosses the side $\textbf{a}_k$ if and only if $x\in I_k$.
\end{rmk}

The condition (C$1$) implies that the trajectory crosses the longest side every other crossing.
It is hence enough to study the squared map $\Phi^2$, which is an IET without flip, see Figure \ref{fig:Phi-T}.
We consider $T=\Phi^2_{\mid I}$ the square of $\Phi$, restricted on $I = I_1 \cup ... \cup I_{N-1} $. Recall that we have normalized the situation so that $a_1+...+a_{N-1}=1$ so $T$ is defined on $(0,1)$.
It is an IET (without flip) with permutation 
$$
\pi(T) = \left(
\begin{array}{cccccc}
I_1&I_2&...&I_{N-1}\\
I_{N-1}&...&I_2&I_1
\end{array}
\right)
$$
and lengths 
$$\lambda(T) = \left(
\begin{array}{cccccc}
a_1&a_2&...&a_{N-1}
\end{array}
\right).
$$

\begin{figure}[h!]
\begin{subfigure}[b]{0.9\textwidth}
         \centering
\definecolor{ffffqq}{rgb}{1,1,0} 
\definecolor{xdxdff}{rgb}{0.49,0.49,1}
\definecolor{ttccqq}{rgb}{0.2,0.8,0}
\definecolor{ffqqqq}{rgb}{1,0,0} 
\definecolor{uququq}{rgb}{0.25,0.25,0.25}
\def \h{0.5}
\def \a{2}
\def \b{1.3}
\def \c{0.3}
\def \d{0.6}
\def \e{5.8}
\def \t{4.6}
\begin{tikzpicture}[line cap=round,line join=round,>=stealth,x=1.0cm,y=1.0cm]
\draw (0,0) -- (10,0) ;
\draw (0,-1) -- (10,-1) ;
\draw (0,0) -- (0,\h) -- (\a,0) -- (\a,\h) -- (\a+\b,0) -- (\a+\b, \h) -- (\a+\b+\c, 0) -- (\a+\b+\c,\h) -- (\a+\b+\c+\d, 0) -- (\a+\b+\c+\d,\h)-- (\a+\b+\c+\d+\e,0);

\fill[color=red,fill=red,fill opacity=0.52] (0,0) -- (0,\h) -- (\a,0)  -- cycle;
\fill[color=yellow,fill=yellow,fill opacity=0.52] (\a,0) -- (\a,\h) -- (\a+\b,0) -- cycle;
\fill[color=blue,fill=blue,fill opacity=0.52] (\a+\b,0) -- (\a+\b, \h) -- (\a+\b+\c, 0) -- cycle;
\fill[color=purple,fill=purple,fill opacity=0.52] (\a+\b+\c, 0) -- (\a+\b+\c,\h) -- (\a+\b+\c+\d, 0) -- cycle;
\fill[color=gray,fill=gray,fill opacity=0.52] (\a+\b+\c+\d, 0) -- (\a+\b+\c+\d,\h) -- (\t,0.93*\h) -- (\t,0) -- cycle;
\fill[color=green,fill=green,fill opacity=0.52] (\t,0.93*\h) -- (\t,0) --  (\t+\a+\b+\c+\d, 0) -- (\t+\a+\b+\c+\d, 0.2*\h) -- cycle;
\fill[color=gray,fill=gray,fill opacity=0.52] (\t+\a+\b+\c+\d, 0) -- (\t+\a+\b+\c+\d, 0.2*\h) -- (10,0) -- cycle;

\draw (0,\h*0.1-1) -- (\t,\h-1) -- (\t,-1) -- (\t+\a,\h-1) -- (\t+\a,-1) -- (\t+\a+\b,\h-1) -- (\t+\a+\b,-1) -- (\t+\a+\b+\c, \h-1) -- (\t+\a+\b+\c, -1) -- (\t+\a+\b+\c+\d,\h-1) -- (\t+\a+\b+\c+\d, -1) -- (10,\h*0.2-1) ;

\fill[color=red,fill=red,fill opacity=0.52] (\t,-1) -- (\t+\a,\h-1) -- (\t+\a,-1) -- cycle;
\fill[color=yellow,fill=yellow,fill opacity=0.52] (\t+\a,-1) -- (\t+\a+\b,\h-1) -- (\t+\a+\b,-1) -- cycle;
\fill[color=blue,fill=blue,fill opacity=0.52] (\t+\a+\b,-1) -- (\t+\a+\b+\c, \h-1) -- (\t+\a+\b+\c, -1) -- cycle;
\fill[color=purple,fill=purple,fill opacity=0.52] (\t+\a+\b+\c, -1) -- (\t+\a+\b+\c+\d,\h-1) -- (\t+\a+\b+\c+\d, -1) -- cycle;
\fill[color=gray,fill=gray,fill opacity=0.52] (\t+\a+\b+\c+\d, -1) -- (10,\h*0.2-1) -- (\a+\b+\c+\d+\e,-1) -- cycle;
\fill[color=green,fill=green,fill opacity=0.52] (0, -1) -- (0,\h*0.1-1) -- (\a+\b+\c+\d, -1+0.93*\h) -- (\a+\b+\c+\d, -1) -- cycle;
\fill[color=gray,fill=gray,fill opacity=0.52] (\a+\b+\c+\d, -1) -- (\a+\b+\c+\d, -1+0.93*\h) -- (\t,\h-1) -- (\t,-1) -- cycle;

\draw (1.25,0)-- ++(-1.5pt,-1.5pt) -- ++(3.0pt,3.0pt) ++(-3.0pt,0) -- ++(3.0pt,-3.0pt);
\draw (\t+\a-1.25,0)-- ++(-1.5pt,-1.5pt) -- ++(3.0pt,3.0pt) ++(-3.0pt,0) -- ++(3.0pt,-3.0pt);
\draw (\b+\c+\d+1.25,0)-- ++(-1.5pt,-1.5pt) -- ++(3.0pt,3.0pt) ++(-3.0pt,0) -- ++(3.0pt,-3.0pt);

\draw [->] (0,-1.3)--(\t,-1.3);

\node (a) at (0,0) {};
\node (b) at (0,-1) {};
\draw[->] (a) to[out=180,in=180] (b) ;
\draw (-0.7,-0.5) node {$\Phi$};

\begin{scriptsize}
\draw (1.25,-0.25) node {$x_0$};
\draw (\t+\a-1.25,-0.25) node {$\Phi(x_0)$};
\draw (\b+\c+\d+1.25,-0.25) node {$\Phi^2(x_0)$};
\draw (\t+0.2,-1.3) node {$\tau$};
\end{scriptsize}
\end{tikzpicture}
\caption{The IETF $\Phi$}
\label{fig:Phi}
\end{subfigure}

\vspace*{0.5cm}

\begin{subfigure}[b]{0.9\textwidth}
         \centering
\definecolor{ffffqq}{rgb}{1,1,0} 
\definecolor{xdxdff}{rgb}{0.49,0.49,1}
\definecolor{ttccqq}{rgb}{0.2,0.8,0}
\definecolor{ffqqqq}{rgb}{1,0,0} 
\definecolor{uququq}{rgb}{0.25,0.25,0.25}
\def \h{0.5}
\def \a{2}
\def \b{1.3}
\def \c{0.3}
\def \d{0.6}
\def \e{5.8}
\def \t{4.6}
\begin{tikzpicture}[line cap=round,line join=round,>=stealth,x=1.0cm,y=1.0cm]

\fill[color=red,fill=red,fill opacity=0.52] (0,0) -- (0,\h) -- (\a,0)  -- cycle;
\fill[color=yellow,fill=yellow,fill opacity=0.52] (\a,0) -- (\a,\h) -- (\a+\b,0) -- cycle;
\fill[color=blue,fill=blue,fill opacity=0.52] (\a+\b,0) -- (\a+\b, \h) -- (\a+\b+\c, 0) -- cycle;
\fill[color=purple,fill=purple,fill opacity=0.52] (\a+\b+\c, 0) -- (\a+\b+\c,\h) -- (\a+\b+\c+\d, 0) -- cycle;
\fill[color=gray,fill=gray,fill opacity=0.52] (\a+\b+\c+\d, 0) -- (\a+\b+\c+\d,\h) -- (\t,0.93*\h) -- (\t,0) -- cycle;
\fill[color=green,fill=green,fill opacity=0.2] (\t,0) -- (\t,0.93*\h) --  (\t+\d, 0.83*\h) -- (\t+\d, 0) -- cycle;
\fill[color=green,fill=green,fill opacity=0.4] (\t+\d,0) -- (\t+\d,0.83*\h) --  (\t+\d+\c, 0.78*\h) -- (\t+\d+\c, 0) -- cycle;
\fill[color=green,fill=green,fill opacity=0.6] (\t+\d+\c,0) -- (\t+\d+\c,0.78*\h) --  (\t+\d+\c+\b, 0.55*\h) -- (\t+\d+\c+\b, 0) -- cycle;
\fill[color=green,fill=green,fill opacity=0.9] (\t+\d+\c+\b,0) -- (\t+\d+\c+\b,0.55*\h) --  (\t+\d+\c+\b+\a, 0.21*\h) -- (\t+\d+\c+\b+\a, 0) -- cycle;
\fill[color=gray,fill=gray,fill opacity=0.52] (\t+\a+\b+\c+\d, 0) -- (\t+\a+\b+\c+\d, 0.21*\h) -- (10,0) -- cycle;

\draw (0,0) -- (0,\h) -- (\a,0) -- (\a,\h) -- (\a+\b,0) -- (\a+\b, \h) -- (\a+\b+\c, 0) -- (\a+\b+\c,\h) -- (\a+\b+\c+\d, 0) -- (\a+\b+\c+\d,\h)-- (10,0);

\fill[color=purple,fill=purple,fill opacity=0.52] (0,-1) -- (0,\h-1) -- (\d,-1)  -- cycle;
\fill[color=blue,fill=blue,fill opacity=0.52] (\d,-1) -- (\d,\h-1) -- (\d+\c,-1) -- cycle;
\fill[color=yellow,fill=yellow,fill opacity=0.52] (\d+\c,-1) -- (\d+\c, \h-1) -- (\d+\b+\c, -1) -- cycle;
\fill[color=red,fill=red,fill opacity=0.52] (\d+\b+\c, -1) -- (\d+\b+\c,\h-1) -- (\d+\b+\c+\a, -1) -- cycle;

\fill[color=gray,fill=gray,fill opacity=0.52] (\a+\b+\c+\d, -1) -- (\a+\b+\c+\d,\h-1) -- (\t,0.93*\h-1) -- (\t,-1) -- cycle;
\fill[color=green,fill=green,fill opacity=0.2] (\t+\a+\b+\c,-1) -- (\t+\a+\b+\c,0.93*\h-1) --  (\t+\a+\b+\c+\d, 0.83*\h-1) -- (\t+\a+\b+\c+\d, -1) -- cycle;
\fill[color=green,fill=green,fill opacity=0.4] (\t+\a+\b,-1) -- (\t+\a+\b,0.83*\h-1) --  (\t+\a+\b+\c, 0.78*\h-1) -- (\t+\a+\b+\c, -1) -- cycle;
\fill[color=green,fill=green,fill opacity=0.6] (\t+\a,-1) -- (\t+\a,0.78*\h-1) --  (\t+\a+\b, 0.55*\h-1) -- (\t+\a+\b, -1) -- cycle;
\fill[color=green,fill=green,fill opacity=0.9] (\t,-1) -- (\t,0.55*\h-1) --  (\t+\a, 0.21*\h-1) -- (\t+\a, -1) -- cycle;
\fill[color=gray,fill=gray,fill opacity=0.52] (\t+\a+\b+\c+\d, -1) -- (\t+\a+\b+\c+\d, 0.21*\h-1) -- (10,-1) -- cycle;

\draw (0,-1) -- (0,\h-1) -- (\d,-1) -- (\d,\h-1) -- (\d+\c,-1) -- (\d+\c, \h-1) -- (\d+\b+\c, -1) -- (\d+\b+\c,\h-1) -- (\a+\b+\c+\d, -1) -- (\a+\b+\c+\d,\h-1)--(\t,0.93*\h-1)-- (\t,-1) --(\t,0.55*\h-1)-- (\t+\a,0.21*\h-1)--(\t+\a,-1)-- (\t+\a,0.78*\h-1)--(\t+\a+\b,0.55*\h-1)--(\t+\a+\b,-1)--(\t+\a+\b,0.83*\h-1)--(\t+\a+\b+\c,0.78*\h-1)--(\t+\a+\b+\c,-1)--(\t+\a+\b+\c,0.93*\h-1)--(\t+\a+\b+\c+\d,0.83*\h-1) --(\t+\a+\b+\c+\d,-1)--(\t+\a+\b+\c+\d,0.83*\h-1) --(\t+\a+\b+\c+\d,-1) --(\t+\a+\b+\c+\d,0.21*\h-1) -- (\a+\b+\c+\d+\e,-1) ;

\draw (0,0) -- (10,0) ;
\draw (0,-1) -- (10,-1) ;

\draw (1.25,0)-- ++(-1.5pt,-1.5pt) -- ++(3.0pt,3.0pt) ++(-3.0pt,0) -- ++(3.0pt,-3.0pt);
\draw (\b+\c+\d+1.25,0)-- ++(-1.5pt,-1.5pt) -- ++(3.0pt,3.0pt) ++(-3.0pt,0) -- ++(3.0pt,-3.0pt);

\draw [->] (0,-1.3)--(\t,-1.3);

\node (a) at (0,0) {};
\node (b) at (0,-1) {};
\draw[->] (a) to[out=180,in=180] (b) ;
\draw (-0.7,-0.5) node {$\Phi^2$};

\begin{scriptsize}
\draw (1.25,-0.25) node {$x_0$};
\draw (\b+\c+\d+1.25,-0.25) node {$\Phi^2(x_0)$};
\draw (\t+0.2,-1.3) node {$\tau$};
\end{scriptsize}
\end{tikzpicture}
\caption{The IET $\Phi^2$ without flip}
\label{fig:T}
\end{subfigure}

\medskip
\small
Orbits of points inside grey parts do not correspond to a trajectory of a tiling billiard.
\caption{Example of the IETs $\Phi$ and $\Phi^2$, corresponding to a tiling billiard.}
\label{fig:Phi-T}
\end{figure}

%

We give here a straightforward way to compute the coordinate $x_{i+1}$ depending on $x_i$, and on the parameters $a_1,\dots, a_N$ and $\tau$. 
Because of the refraction law, the lengths $\alpha$ and $\beta$ on the Figure  \ref{fig:computation_IET} are equal. Since $\alpha = x_i - \sum_{k=1}^{j-1} a_k$, one gets that $\gamma = \tau - \beta = \tau - x_i + \sum_{k=1}^{j} a_k$. We conclude the expression of $x_{i+1}$ depending on $x_i$, $a_1,\dots, a_N$ and $\tau$:
\begin{align*}
x_{i+1} & = \gamma + \sum_{k=1}^j a_k \\
		& = \tau - (x_i -a_j)
\end{align*}
which corresponds to the expression of a fully flipped IET.

%
%

\section{Background}\label{background}



\subsection{Rauzy induction, Rokhlin towers, and the associated cocycle}\label{subsec_Rauzy-cocycle-Rokhlin}

The Rauzy induction $\mathfrak{R}(T)$ of an IET $T$ defined on $I$ is the first return map of $T$ to a well chosen subinterval $J\subset I$. It can be convenient to consider the renormalized Rauzy induction $\mathcal{R}$ that consists to remormalise $J$ to the size of $I$. For more details, see \cite{Viana}. 
Zorich \cite{Zorich_finite_measure} showed that one can speed up the Rauzy induction so that the induction is ergodic with respect to an absolutely continuous invariant probability measure $\mu_Z$.
The Rauzy induction acts also on the underlying permutation of the IET.
We call the \textbf{Rauzy class of the permutation $\sigma$} the permutations in the orbit of $\sigma$ by the Rauzy induction.

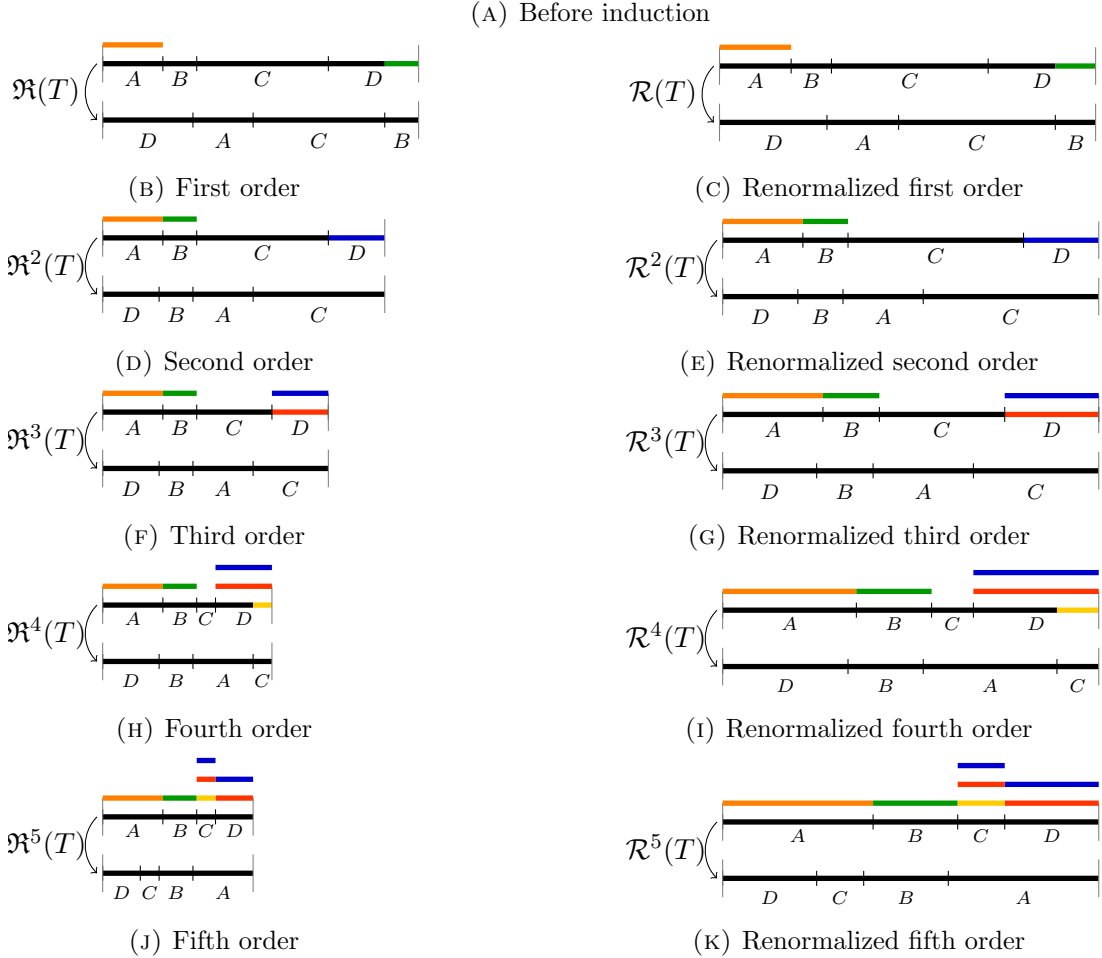
\begin{figure}[h]
\center
\begin{subfigure}[b]{0.8\textwidth}
         \centering
\definecolor{cola}{rgb}{1,0.5,0}
\begin{tikzpicture}[line join=round,x=0.5cm,y=0.5cm]
\def \a{1.6}
\def \b{0.9}
\def \c{3.5}
\def \d{4}
\def \h{1.5} 
\def \hh{2} 
\def \r{1.25} 
\def \t{0.5} 
\clip (-2,-\h-1) rectangle (10.5,0.5);
\draw[color=black, line width=2pt] (0,0) -- (\d+\c+\b,0);
\draw[color=cola, line width=2pt] (\d+\c+\b,0) -- (10,0); 
\foreach \x in {0,\a,\a+\b,\a+\b+\c}
\draw[shift={(\x,0)},color=black] (0,-0.2)--(0,0.2);
\draw[gray] (0,-0.5) -- (0,0.5);
\draw[gray] (10,-0.5) -- (10,0.5);
\draw[color=black, line width=2pt] (0,-\h) -- (10,-\h);
\foreach \x in {0,\d,\d+\c,\d+\c+\b}
\draw[shift={(\x,-\h)},color=black] (0,-0.2)--(0,0.2);
\draw[gray] (0,-\h-0.5) -- (0,-\h+0.5);
\draw[gray] (10,-\h-0.5) -- (10,-\h+0.5);

\begin{scriptsize}
\draw (\a/2,0) node (A) [below]{$A$};
\draw (\a+\b/2,0) node (B) [below]{$B$};
\draw (\a+\b+\c/2,0) node (C) [below]{$C$};
\draw (\a+\b+\c+\d/2,0) node (D) [below]{$D$};

\draw (\d+\c+\b/2,-\h) node (A') [below=2]{$B$};
\draw (\d+\c/2,-\h) node (B') [below=2]{$C$};
\draw (\d+\c+\b+\a/2,-\h) node (C') [below=2]{$A$};
\draw (\d/2,-\h) node (D') [below=2]{$D$};
\end{scriptsize}
\draw [->] (-0.15,0) to [out = -135,in =135] (-0.15,-\h);
\draw (-1.5,-0.5*\h) node {$T$}; 
\end{tikzpicture}

\medskip
\small
The colored subinterval will be cut in the following step of Rauzy induction (and stacked over the interval to form a Rokhling tower).
\caption{Before induction}
\end{subfigure}
\hfill

\begin{subfigure}[b]{0.35\textwidth}
         \centering
\definecolor{cola}{rgb}{1,0.5,0}
\definecolor{colb}{rgb}{0,0.6,0}
\definecolor{ffccqq}{rgb}{1,0.8,0}
\definecolor{ffttqq}{rgb}{1,0.2,0}
\definecolor{qqqqcc}{rgb}{0,0,0.8}
\definecolor{tttttt}{rgb}{0.2,0.2,0.2}
\begin{tikzpicture}[line join=round,x=0.5cm,y=0.5cm]
\def \a{1.6}
\def \b{0.9}
\def \c{3.5}
\def \d{4}
\def \h{1.5} 
\def \hh{2} 
\def \r{1.25} 
\def \t{0.5} 
\clip (-2.5,-\h-1) rectangle (10.5,0.8);
\draw[color=black, line width=2pt] (0,0) -- (\c+\d,0);
\draw[color=colb, line width=2pt] (\c+\d,0)--(\b+\c+\d,0);
\foreach \x in {0,\a,\a+\b,\a+\b+\c}
\draw[shift={(\x,0)},color=black] (0,-0.2)--(0,0.2);
\draw[gray] (0,-0.5) -- (0,0.5);
\draw[gray] (\b+\c+\d,-0.5) -- (\b+\c+\d,0.5);
\draw[color=black, line width=2pt] (0,-\h) -- (\b+\c+\d,-\h);
\foreach \x in {0,\d-\a,\d,\d+\c}
\draw[shift={(\x,-\h)},color=black] (0,-0.2)--(0,0.2);
\draw[gray] (0,-\h-0.5) -- (0,-\h+0.5);
\draw[gray] (\b+\c+\d,-\h-0.5) -- (\b+\c+\d,-\h+0.5);

\draw[color=cola, line width=2pt] (0,\t) -- (\a,\t); 

\begin{scriptsize}
\draw (\a/2,0) node (A) [below]{$A$};
\draw (\a+\b/2,0) node (B) [below]{$B$};
\draw (\a+\b+\c/2,0) node (C) [below]{$C$};
\draw (\a+\b+\c+\d/2-\a/2,0) node (D) [below]{$D$};

\draw (\d+\c+\b/2,-\h) node (A') [below=2]{$B$};
\draw (\d+\c/2,-\h) node (B') [below=2]{$C$};
\draw (\d-\a/2,-\h) node (C') [below=2]{$A$};
\draw (\d/2-\a/2,-\h) node (D') [below=2]{$D$};
\end{scriptsize}
\draw [->] (-0.15,0) to [out = -135,in =135] (-0.15,-\h);
\draw (-1.5,-0.5*\h) node {$\mathfrak{R}(T)$}; 
\end{tikzpicture}
\caption{First order}
\end{subfigure}
\hfill
\begin{subfigure}[b]{0.55\textwidth}
         \centering
\definecolor{cola}{rgb}{1,0.5,0}
\definecolor{colb}{rgb}{0,0.6,0}
\definecolor{ffccqq}{rgb}{1,0.8,0}
\definecolor{ffttqq}{rgb}{1,0.2,0}
\definecolor{qqqqcc}{rgb}{0,0,0.8}
\definecolor{tttttt}{rgb}{0.2,0.2,0.2}
\begin{tikzpicture}[line join=round,x=0.595238095cm,y=0.5cm]
\def \a{1.6}
\def \b{0.9}
\def \c{3.5}
\def \d{4}
\def \h{1.5} 
\def \hh{2} 
\def \r{1.19047619} 
\def \t{0.5} 
\draw[color=black, line width=2pt] (0,0) -- (\c+\d,0);
\draw[color=colb, line width=2pt] (\c+\d,0)--(\b+\c+\d,0);
\foreach \x in {0,\a,\a+\b,\a+\b+\c}
\draw[shift={(\x,0)},color=black] (0,-0.2)--(0,0.2);
\draw[gray] (0,-0.5) -- (0,0.5);
\draw[gray] (\b+\c+\d,-0.5) -- (\b+\c+\d,0.5);
\draw[color=black, line width=2pt] (0,-\h) -- (\b+\c+\d,-\h);
\foreach \x in {0,\d-\a,\d,\d+\c}
\draw[shift={(\x,-\h)},color=black] (0,-0.2)--(0,0.2);
\draw[gray] (0,-\h-0.5) -- (0,-\h+0.5);
\draw[gray] (\b+\c+\d,-\h-0.5) -- (\b+\c+\d,-\h+0.5);

\draw[color=cola, line width=2pt] (0,\t) -- (\a,\t); 

\begin{scriptsize}
\draw (\a/2,0) node (A) [below]{$A$};
\draw (\a+\b/2,0) node (B) [below]{$B$};
\draw (\a+\b+\c/2,0) node (C) [below]{$C$};
\draw (\a+\b+\c+\d/2-\a/2,0) node (D) [below]{$D$};

\draw (\d+\c+\b/2,-\h) node (A') [below=2]{$B$};
\draw (\d+\c/2,-\h) node (B') [below=2]{$C$};
\draw (\d-\a/2,-\h) node (C') [below=2]{$A$};
\draw (\d/2-\a/2,-\h) node (D') [below=2]{$D$};
\end{scriptsize}
\draw [->] (-0.15/\r,0) to [out = -135,in =135] (-0.15/\r,-\h);
\draw (-1.5/\r,-0.5*\h) node {$\mathcal{R}(T)$}; 
\end{tikzpicture}
\caption{Renormalized first order}
\end{subfigure}
\begin{subfigure}[b]{0.35\textwidth}
         \centering
\definecolor{cola}{rgb}{1,0.5,0}
\definecolor{colb}{rgb}{0,0.6,0}
\definecolor{ffccqq}{rgb}{1,0.8,0}
\definecolor{ffttqq}{rgb}{1,0.2,0}
\definecolor{colc}{rgb}{0,0,0.8}
\definecolor{tttttt}{rgb}{0.2,0.2,0.2}
\begin{tikzpicture}[line join=round,x=0.5cm,y=0.5cm]
\def \a{1.6}
\def \b{0.9}
\def \c{3.5}
\def \d{4}
\def \h{1.5} 
\def \hh{2} 
\def \r{1.25} 
\def \t{0.5} 
\clip (-2.5,-\h-1) rectangle (10.5,0.8);
\draw[color=black, line width=2pt] (0,0) -- (\a+\b+\c,0);
\draw[color=colc, line width=2pt] (\a+\b+\c,0)--(\c+\d,0);
\foreach \x in {0,\a,\a+\b,\a+\b+\c}
\draw[shift={(\x,0)},color=black] (0,-0.2)--(0,0.2);
\draw[gray] (0,-0.5) -- (0,0.5);
\draw[gray] (\c+\d,-0.5) -- (\c+\d,0.5);
\draw[color=black, line width=2pt] (0,-\h) -- (\c+\d,-\h);
\foreach \x in {0,\d-\a-\b,\d-\a,\d}
\draw[shift={(\x,-\h)},color=black] (0,-0.2)--(0,0.2);
\draw[gray] (0,-\h-0.5) -- (0,-\h+0.5);
\draw[gray] (\c+\d,-\h-0.5) -- (\c+\d,-\h+0.5);

\draw[color=cola, line width=2pt] (0,\t) -- (\a,\t); 
\draw[color=colb, line width=2pt] (\a,\t) -- (\a+\b,\t); 

\begin{scriptsize}
\draw (\a/2,0) node (A) [below]{$A$};
\draw (\a+\b/2,0) node (B) [below]{$B$};
\draw (\a+\b+\c/2,0) node (C) [below]{$C$};
\draw (\a+\b+\c+\d/2-\a/2-\b/2,0) node (D) [below]{$D$};

\draw (\d-\a-\b/2,-\h) node (A') [below=2]{$B$};
\draw (\d+\c/2,-\h) node (B') [below=2]{$C$};
\draw (\d-\a/2,-\h) node (C') [below=2]{$A$};
\draw (\d/2-\a/2-\b/2,-\h) node (D') [below=2]{$D$};
\end{scriptsize}
\draw [->] (-0.15,0) to [out = -135,in =135] (-0.15,-\h);
\draw (-1.5,-0.5*\h) node {$\mathfrak{R}^2(T)$}; 
\end{tikzpicture}
\caption{Second order}
\end{subfigure}
\hfill
\begin{subfigure}[b]{0.55\textwidth}
         \centering
\definecolor{cola}{rgb}{1,0.5,0}
\definecolor{colb}{rgb}{0,0.6,0}
\definecolor{ffccqq}{rgb}{1,0.8,0}
\definecolor{ffttqq}{rgb}{1,0.2,0}
\definecolor{colc}{rgb}{0,0,0.8}
\definecolor{tttttt}{rgb}{0.2,0.2,0.2}
\begin{tikzpicture}[line join=round,x=0.666666667cm,y=0.5cm]
\def \a{1.6}
\def \b{0.9}
\def \c{3.5}
\def \d{4}
\def \h{1.5} 
\def \hh{2} 
\def \r{1.333333333} 
\def \t{0.5} 
\draw[color=black, line width=2pt] (0,0) -- (\a+\b+\c,0);
\draw[color=colc, line width=2pt] (\a+\b+\c,0)--(\c+\d,0);
\foreach \x in {0,\a,\a+\b,\a+\b+\c}
\draw[shift={(\x,0)},color=black] (0,-0.2)--(0,0.2);
\draw[gray] (0,-0.5) -- (0,0.5);
\draw[gray] (\c+\d,-0.5) -- (\c+\d,0.5);
\draw[color=black, line width=2pt] (0,-\h) -- (\c+\d,-\h);
\foreach \x in {0,\d-\a-\b,\d-\a,\d}
\draw[shift={(\x,-\h)},color=black] (0,-0.2)--(0,0.2);
\draw[gray] (0,-\h-0.5) -- (0,-\h+0.5);
\draw[gray] (\c+\d,-\h-0.5) -- (\c+\d,-\h+0.5);

\draw[color=cola, line width=2pt] (0,\t) -- (\a,\t); 
\draw[color=colb, line width=2pt] (\a,\t) -- (\a+\b,\t); 

\begin{scriptsize}
\draw (\a/2,0) node (A) [below]{$A$};
\draw (\a+\b/2,0) node (B) [below]{$B$};
\draw (\a+\b+\c/2,0) node (C) [below]{$C$};
\draw (\a+\b+\c+\d/2-\a/2-\b/2,0) node (D) [below]{$D$};

\draw (\d-\a-\b/2,-\h) node (A') [below=2]{$B$};
\draw (\d+\c/2,-\h) node (B') [below=2]{$C$};
\draw (\d-\a/2,-\h) node (C') [below=2]{$A$};
\draw (\d/2-\a/2-\b/2,-\h) node (D') [below=2]{$D$};
\end{scriptsize}
\draw [->] (-0.15/\r,0) to [out = -135,in =135] (-0.15/\r,-\h);
\draw (-1.5/\r,-0.5*\h) node {$\mathcal{R}^2(T)$}; 
\end{tikzpicture}
\caption{Renormalized second order}
\end{subfigure}

\begin{subfigure}[b]{0.35\textwidth}
         \centering
\definecolor{cola}{rgb}{1,0.5,0}
\definecolor{colb}{rgb}{0,0.6,0}
\definecolor{ffccqq}{rgb}{1,0.8,0}
\definecolor{cold}{rgb}{1,0.2,0}
\definecolor{colc}{rgb}{0,0,0.8}
\definecolor{tttttt}{rgb}{0.2,0.2,0.2}
\begin{tikzpicture}[line join=round,x=0.5cm,y=0.5cm]
\def \a{1.6}
\def \b{0.9}
\def \c{3.5}
\def \d{4}
\def \h{1.5} 
\def \hh{2} 
\def \r{1.25} 
\def \t{0.5} 
\clip (-2.5,-\h-1) rectangle (10.5,0.8);
\draw[color=black, line width=2pt] (0,0) -- (2*\a+2*\b+\c-\d,0);
\draw[color=cold, line width=2pt] (2*\a+2*\b+\c-\d,0) -- (\a+\b+\c,0);
\foreach \x in {0,\a,\a+\b,2*\a+2*\b+\c-\d,\a+\b+\c}
\draw[shift={(\x,0)},color=black] (0,-0.2)--(0,0.2);
\draw[gray] (0,-0.5) -- (0,0.5);
\draw[gray] (\a+\b+\c,-0.5) -- (\a+\b+\c,0.5);
\draw[color=black, line width=2pt] (0,-\h) -- (\a+\b+\c,-\h);
\foreach \x in {0,\d-\a-\b,\d-\a,\d}
\draw[shift={(\x,-\h)},color=black] (0,-0.2)--(0,0.2);
\draw[gray] (0,-\h-0.5) -- (0,-\h+0.5);
\draw[gray] (\a+\b+\c,-\h-0.5) -- (\a+\b+\c,-\h+0.5);


\draw[color=cola, line width=2pt] (0,\t) -- (\a,\t); 
\draw[color=colb, line width=2pt] (\a,\t) -- (\a+\b,\t); 
\draw[color=colc, line width=2pt] (2*\a+2*\b+\c-\d,\t) -- (\a+\b+\c,\t);

\begin{scriptsize}
\draw (\a/2,0) node (A) [below]{$A$};
\draw (\a+\b/2,0) node (B) [below]{$B$};
\draw (\a+\b+\c/2-\d/2+\a/2+\b/2,0) node (C) [below]{$C$};
\draw (2*\a+2*\b+\c-\d+\d/2-\a/2-\b/2,0) node (D) [below]{$D$};


\draw (\d-\a-\b/2,-\h) node (A') [below=2]{$B$};
\draw (\d+\c/2-\d/2+\a/2+\b/2,-\h) node (B') [below=2]{$C$};
\draw (\d-\a/2,-\h) node (C') [below=2]{$A$};
\draw (\d/2-\a/2-\b/2,-\h) node (D') [below=2]{$D$};
\end{scriptsize}
\draw [->] (-0.15,0) to [out = -135,in =135] (-0.15,-\h);
\draw (-1.5,-0.5*\h) node {$\mathfrak{R}^3(T)$}; 
\end{tikzpicture}
\caption{Third order}
\end{subfigure}
\hfill
\begin{subfigure}[b]{0.55\textwidth}
         \centering
\definecolor{cola}{rgb}{1,0.5,0}
\definecolor{colb}{rgb}{0,0.6,0}
\definecolor{ffccqq}{rgb}{1,0.8,0}
\definecolor{cold}{rgb}{1,0.2,0}
\definecolor{colc}{rgb}{0,0,0.8}
\definecolor{tttttt}{rgb}{0.2,0.2,0.2}
\begin{tikzpicture}[line join=round,x=0.833333333cm,y=0.5cm]
\def \a{1.6}
\def \b{0.9}
\def \c{3.5}
\def \d{4}
\def \h{1.5} 
\def \hh{2} 
\def \r{1.66666667} 
\def \t{0.5} 
\draw[color=black, line width=2pt] (0,0) -- (2*\a+2*\b+\c-\d,0);
\draw[color=cold, line width=2pt] (2*\a+2*\b+\c-\d,0) -- (\a+\b+\c,0);
\foreach \x in {0,\a,\a+\b,2*\a+2*\b+\c-\d,\a+\b+\c}
\draw[shift={(\x,0)},color=black] (0,-0.2)--(0,0.2);
\draw[gray] (0,-0.5) -- (0,0.5);
\draw[gray] (\a+\b+\c,-0.5) -- (\a+\b+\c,0.5);
\draw[color=black, line width=2pt] (0,-\h) -- (\a+\b+\c,-\h);
\foreach \x in {0,\d-\a-\b,\d-\a,\d}
\draw[shift={(\x,-\h)},color=black] (0,-0.2)--(0,0.2);
\draw[gray] (0,-\h-0.5) -- (0,-\h+0.5);
\draw[gray] (\a+\b+\c,-\h-0.5) -- (\a+\b+\c,-\h+0.5);


\draw[color=cola, line width=2pt] (0,\t) -- (\a,\t); 
\draw[color=colb, line width=2pt] (\a,\t) -- (\a+\b,\t); 
\draw[color=colc, line width=2pt] (2*\a+2*\b+\c-\d,\t) -- (\a+\b+\c,\t);

\begin{scriptsize}
\draw (\a/2,0) node (A) [below]{$A$};
\draw (\a+\b/2,0) node (B) [below]{$B$};
\draw (\a+\b+\c/2-\d/2+\a/2+\b/2,0) node (C) [below]{$C$};
\draw (2*\a+2*\b+\c-\d+\d/2-\a/2-\b/2,0) node (D) [below]{$D$};


\draw (\d-\a-\b/2,-\h) node (A') [below=2]{$B$};
\draw (\d+\c/2-\d/2+\a/2+\b/2,-\h) node (B') [below=2]{$C$};
\draw (\d-\a/2,-\h) node (C') [below=2]{$A$};
\draw (\d/2-\a/2-\b/2,-\h) node (D') [below=2]{$D$};
\end{scriptsize}
\draw [->] (-0.15/\r,0) to [out = -135,in =135] (-0.15/\r,-\h);
\draw (-1.5/\r,-0.5*\h) node {$\mathcal{R}^3(T)$}; 
\end{tikzpicture}
\caption{Renormalized third order}
\end{subfigure}

\begin{subfigure}[b]{0.35\textwidth}
         \centering
\definecolor{cola}{rgb}{1,0.5,0}
\definecolor{colb}{rgb}{0,0.6,0}
\definecolor{cole}{rgb}{1,0.8,0}
\definecolor{cold}{rgb}{1,0.2,0}
\definecolor{colc}{rgb}{0,0,0.8}
\definecolor{tttttt}{rgb}{0.2,0.2,0.2}
\begin{tikzpicture}[line join=round,x=0.5cm,y=0.5cm]
\def \a{1.6}
\def \b{0.9}
\def \c{3.5}
\def \d{4}
\def \h{1.5} 
\def \hh{2} 
\def \r{1.25} 
\def \t{0.5} 
\clip (-2.5,-\h-1) rectangle (10.5,2*\t+0.3);
\draw[color=black, line width=2pt] (0,0) -- (\d,0);
\draw[color=cole, line width=2pt] (\d,0) -- (2*\a+2*\b+\c-\d,0);
\foreach \x in {0,\a,\a+\b,3*\a+3*\b+\c-2*\d}
\draw[shift={(\x,0)},color=black] (0,-0.2)--(0,0.2);
\draw[gray] (0,-0.5) -- (0,0.5);
\draw[gray] (2*\a+2*\b+\c-\d,-0.5) -- (2*\a+2*\b+\c-\d,0.5);
\draw[color=black, line width=2pt] (0,-\h) -- (2*\a+2*\b+\c-\d,-\h);
\foreach \x in {0,\d-\a-\b,\d-\a,\d}
\draw[shift={(\x,-\h)},color=black] (0,-0.2)--(0,0.2);
\draw[gray] (0,-\h-0.5) -- (0,-\h+0.5);
\draw[gray] (2*\a+2*\b+\c-\d,-\h-0.5) -- (2*\a+2*\b+\c-\d,-\h+0.5);

\draw[color=cola, line width=2pt] (0,\t) -- (\a,\t); 
\draw[color=colb, line width=2pt] (\a,\t) -- (\a+\b,\t); 
\draw[color=cold, line width=2pt] (3*\a+3*\b+\c-2*\d,\t) -- (2*\a+2*\b+\c-\d,\t);
\draw[color=colc, line width=2pt] (3*\a+3*\b+\c-2*\d,2*\t) -- (2*\a+2*\b+\c-\d,2*\t);

\begin{tiny}
\draw (\a/2,0) node (A) [below]{$A$};
\draw (\a+\b/2,0) node (B) [below]{$B$};
\draw (\a+\b+\c/2-\d+\a+\b,0) node (C) [below]{$C$};
\draw (3*\a+3*\b+\c-2*\d+\d/2-\a/2-\b/2,0) node (D) [below]{$D$};

\draw (\d-\a-\b/2,-\h) node (A') [below=2]{$B$};
\draw (\a+\b+\c/2,-\h) node (B') [below=2]{$C$};
\draw (\d-\a/2,-\h) node (C') [below=2]{$A$};
\draw (\d/2-\a/2-\b/2,-\h) node (D') [below=2]{$D$};
\end{tiny}
\draw [->] (-0.15,0) to [out = -135,in =135] (-0.15,-\h);
\draw (-1.5,-0.5*\h) node {$\mathfrak{R}^4(T)$}; 
\end{tikzpicture}
\caption{Fourth order}
\end{subfigure}
\hfill
\begin{subfigure}[b]{0.55\textwidth}
         \centering
\definecolor{cola}{rgb}{1,0.5,0}
\definecolor{colb}{rgb}{0,0.6,0}
\definecolor{cole}{rgb}{1,0.8,0}
\definecolor{cold}{rgb}{1,0.2,0}
\definecolor{colc}{rgb}{0,0,0.8}
\definecolor{tttttt}{rgb}{0.2,0.2,0.2}
\begin{tikzpicture}[line join=round,x=1.1111111cm,y=0.5cm]
\def \a{1.6}
\def \b{0.9}
\def \c{3.5}
\def \d{4}
\def \h{1.5} 
\def \hh{2} 
\def \r{2.2222222} 
\def \t{0.5} 
\draw[color=black, line width=2pt] (0,0) -- (\d,0);
\draw[color=cole, line width=2pt] (\d,0) -- (2*\a+2*\b+\c-\d,0);
\foreach \x in {0,\a,\a+\b,3*\a+3*\b+\c-2*\d}
\draw[shift={(\x,0)},color=black] (0,-0.2)--(0,0.2);
\draw[gray] (0,-0.5) -- (0,0.5);
\draw[gray] (2*\a+2*\b+\c-\d,-0.5) -- (2*\a+2*\b+\c-\d,0.5);
\draw[color=black, line width=2pt] (0,-\h) -- (2*\a+2*\b+\c-\d,-\h);
\foreach \x in {0,\d-\a-\b,\d-\a,\d}
\draw[shift={(\x,-\h)},color=black] (0,-0.2)--(0,0.2);
\draw[gray] (0,-\h-0.5) -- (0,-\h+0.5);
\draw[gray] (2*\a+2*\b+\c-\d,-\h-0.5) -- (2*\a+2*\b+\c-\d,-\h+0.5);

\draw[color=cola, line width=2pt] (0,\t) -- (\a,\t); 
\draw[color=colb, line width=2pt] (\a,\t) -- (\a+\b,\t); 
\draw[color=cold, line width=2pt] (3*\a+3*\b+\c-2*\d,\t) -- (2*\a+2*\b+\c-\d,\t);
\draw[color=colc, line width=2pt] (3*\a+3*\b+\c-2*\d,2*\t) -- (2*\a+2*\b+\c-\d,2*\t);

\begin{tiny}
\draw (\a/2,0) node (A) [below]{$A$};
\draw (\a+\b/2,0) node (B) [below]{$B$};
\draw (\a+\b+\c/2-\d+\a+\b,0) node (C) [below]{$C$};
\draw (3*\a+3*\b+\c-2*\d+\d/2-\a/2-\b/2,0) node (D) [below]{$D$};

\draw (\d-\a-\b/2,-\h) node (A') [below=2]{$B$};
\draw (\a+\b+\c/2,-\h) node (B') [below=2]{$C$};
\draw (\d-\a/2,-\h) node (C') [below=2]{$A$};
\draw (\d/2-\a/2-\b/2,-\h) node (D') [below=2]{$D$};
\end{tiny}
\draw [->] (-0.15/\r,0) to [out = -135,in =135] (-0.15/\r,-\h);
\draw (-1.5/\r,-0.5*\h) node {$\mathcal{R}^4(T)$}; 
\end{tikzpicture}
\caption{Renormalized fourth order}
\end{subfigure}

\begin{subfigure}[b]{0.35\textwidth}
         \centering
\definecolor{cola}{rgb}{1,0.5,0}
\definecolor{colb}{rgb}{0,0.6,0}
\definecolor{cole}{rgb}{1,0.8,0}
\definecolor{cold}{rgb}{1,0.2,0}
\definecolor{colc}{rgb}{0,0,0.8}
\definecolor{tttttt}{rgb}{0.2,0.2,0.2}
\begin{tikzpicture}[line join=round,x=0.5cm,y=0.5cm]
\def \a{1.6}
\def \b{0.9}
\def \c{3.5}
\def \d{4}
\def \h{1.5} 
\def \hh{2} 
\def \r{1.25} 
\def \t{0.5} 
\clip (-2.5,-\h-1) rectangle (10.5,3*\t+0.3);
\draw[color=black, line width=2pt] (0,0) -- (\d,0);
\foreach \x in {0,\a,\a+\b,3*\a+3*\b+\c-2*\d}
\draw[shift={(\x,0)},color=black] (0,-0.2)--(0,0.2);
\draw[gray] (0,-0.5) -- (0,0.5);
\draw[gray] (\d,-0.5) -- (\d,0.5);
\draw[color=black, line width=2pt] (0,-\h) -- (\d,-\h);
\foreach \x in {0,3*\d-3*\a-3*\b-\c,\d-\a-\b,\d-\a}
\draw[shift={(\x,-\h)},color=black] (0,-0.2)--(0,0.2);
\draw[gray] (0,-\h-0.5) -- (0,-\h+0.5);
\draw[gray] (\d,-\h-0.5) -- (\d,-\h+0.5);

\draw[color=cola, line width=2pt] (0,\t) -- (\a,\t); 
\draw[color=colb, line width=2pt] (\a,\t) -- (\a+\b,\t); 
\draw[color=cold, line width=2pt] (3*\a+3*\b+\c-2*\d,\t) -- (\d,\t);
\draw[color=colc, line width=2pt] (3*\a+3*\b+\c-2*\d,2*\t) -- (\d,2*\t);
\draw[color=cole, line width=2pt] (\a+\b,\t) -- (3*\a+3*\b+\c-2*\d,\t);
\draw[color=cold, line width=2pt] (\a+\b,2*\t) -- (3*\a+3*\b+\c-2*\d,2*\t);
\draw[color=colc, line width=2pt] (\a+\b,3*\t) -- (3*\a+3*\b+\c-2*\d,3*\t);
\begin{tiny}
\draw (\a/2,0) node (A) [below]{$A$};
\draw (\a+\b/2,0) node (B) [below]{$B$};
\draw (\a+\b+\c/2-\d+\a+\b,0) node (C) [below]{$C$};
\draw (3*\a+3*\b+\c-2*\d+3*\d/2-3*\a/2-3*\b/2-\c/2,0) node (D) [below]{$D$};

\draw (\d-\a-\b/2,-\h) node (A') [below=2]{$B$};
\draw (3*\d-3*\a-3*\b-\c+ \a+\b-\d+\c/2,-\h) node (B') [below=2]{$C$};
\draw (\d-\a/2,-\h) node (C') [below=2]{$A$};
\draw (3*\d/2-3*\a/2-3*\b/2-\c/2,-\h) node (D') [below=2]{$D$};
\end{tiny}
\draw [->] (-0.15,0) to [out = -135,in =135] (-0.15,-\h);
\draw (-1.5,-0.5*\h) node {$\mathfrak{R}^5(T)$}; 
\end{tikzpicture}
\caption{Fifth order}
\end{subfigure}
\hfill
\begin{subfigure}[b]{0.55\textwidth}
         \centering
\definecolor{cola}{rgb}{1,0.5,0}
\definecolor{colb}{rgb}{0,0.6,0}
\definecolor{cole}{rgb}{1,0.8,0}
\definecolor{cold}{rgb}{1,0.2,0}
\definecolor{colc}{rgb}{0,0,0.8}
\definecolor{tttttt}{rgb}{0.2,0.2,0.2}
\begin{tikzpicture}[line join=round,x=1.25cm,y=0.5cm]
\def \a{1.6}
\def \b{0.9}
\def \c{3.5}
\def \d{4}
\def \h{1.5} 
\def \hh{2} 
\def \r{2.5} 
\def \t{0.5} 
\draw[color=black, line width=2pt] (0,0) -- (\d,0);
\foreach \x in {0,\a,\a+\b,3*\a+3*\b+\c-2*\d}
\draw[shift={(\x,0)},color=black] (0,-0.2)--(0,0.2);
\draw[gray] (0,-0.5) -- (0,0.5);
\draw[gray] (\d,-0.5) -- (\d,0.5);
\draw[color=black, line width=2pt] (0,-\h) -- (\d,-\h);
\foreach \x in {0,3*\d-3*\a-3*\b-\c,\d-\a-\b,\d-\a}
\draw[shift={(\x,-\h)},color=black] (0,-0.2)--(0,0.2);
\draw[gray] (0,-\h-0.5) -- (0,-\h+0.5);
\draw[gray] (\d,-\h-0.5) -- (\d,-\h+0.5);

\draw[color=cola, line width=2pt] (0,\t) -- (\a,\t); 
\draw[color=colb, line width=2pt] (\a,\t) -- (\a+\b,\t); 
\draw[color=cold, line width=2pt] (3*\a+3*\b+\c-2*\d,\t) -- (\d,\t);
\draw[color=colc, line width=2pt] (3*\a+3*\b+\c-2*\d,2*\t) -- (\d,2*\t);
\draw[color=cole, line width=2pt] (\a+\b,\t) -- (3*\a+3*\b+\c-2*\d,\t);
\draw[color=cold, line width=2pt] (\a+\b,2*\t) -- (3*\a+3*\b+\c-2*\d,2*\t);
\draw[color=colc, line width=2pt] (\a+\b,3*\t) -- (3*\a+3*\b+\c-2*\d,3*\t);
\begin{tiny}
\draw (\a/2,0) node (A) [below]{$A$};
\draw (\a+\b/2,0) node (B) [below]{$B$};
\draw (\a+\b+\c/2-\d+\a+\b,0) node (C) [below]{$C$};
\draw (3*\a+3*\b+\c-2*\d+3*\d/2-3*\a/2-3*\b/2-\c/2,0) node (D) [below]{$D$};

\draw (\d-\a-\b/2,-\h) node (A') [below=2]{$B$};
\draw (3*\d-3*\a-3*\b-\c+ \a+\b-\d+\c/2,-\h) node (B') [below=2]{$C$};
\draw (\d-\a/2,-\h) node (C') [below=2]{$A$};
\draw (3*\d/2-3*\a/2-3*\b/2-\c/2,-\h) node (D') [below=2]{$D$};
\end{tiny}
\draw [->] (-0.15/\r,0) to [out = -135,in =135] (-0.15/\r,-\h);
\draw (-1.5/\r,-0.5*\h) node {$\mathcal{R}^5(T)$}; 
\end{tikzpicture}
\caption{Renormalized fifth order}
\end{subfigure}

\caption{Rokhlin towers of orders 1 to 5 for the Rauzy induction of an IET}\label{fig:Rokhlin}

\end{figure}

When we cut off $I'=I\setminus J$ from $I$, we can think of it being stacked over its preimage by $T$. The choice of $J$ ensures that the subinterval on which $I'$ is being stacked corresponds exactly to an interval of continuity of $\mathfrak{R}(T)$.
We can do that at each iteration of Rauzy induction. We obtain a stack of subintervals of $I$, that we call a \textbf{Rokhlin tower}.
Let $J^{(k)}$ be the subinterval on which $\mathfrak{R}^k(T)$ is defined and $x\in J^{(k)}$. The subintervals stacked over $x$ allow to keep track of the orbit of $x$ under $T$. Their number corresponds to the time of first return of $x$ to $J^{(k)}$. See Figure~\ref{fig:Rokhlin}.
We call \textbf{towers of order $l$} the ones obtained when applying $l$ times the induction\footnote{Note that this notion depends on the induction we apply, Rauzy or Zorich one.}.

We denote by $A_{i,j}^{(l)}$ the number of subintervals stacked over the interval of continuity $I_j^{(l)}$ of $\mathfrak{R}^l(T)$ that were initially in the interval of continuity $I_i$ of $T$. It corresponds to the number of times that the orbit of $x\in I_j^{(l)}$ lies in the \emph{old} interval $I_i$ before first return to $I_j^{(l)}$. This number does not depend on $x\in I_j^{(l)}$ thanks to the choice of subintervals on which we induce.
We denote $A_T = \left( A_{i,j}^{(1)}\right)_{1\leqslant i,j \leqslant d}$ the matrix whose entries are these numbers. It is a transvection. This defines a cocycle $ T\in \text{IET}(I) \mapsto A_T \in \text{SL}_d(\RR)$.

Analogously we denote by $B$ the cocycle associated to Zorich induction.
We underline the fact that the time related to cocyle $B$ is not the same as the one related to cycle $A$, there is a speeding up factor. Moreover there is an exponential factor between the time of cocyle $B$ and the time of the iteration of $T$, see Lemma 4 in \cite{howdo}.

\subsection{Oseledets' theorem}\label{subsec_Oseledets}
We recall Oseledets' theorem on which rely Zorich's results that we will need later in this article.
\begin{thm}[Oseledets]
Let $X$ be a topological space, $\varphi : X\longrightarrow X$ be a measurable function, and $\mu$ be an ergodic $\varphi$-invariant probability mesure on $X$.
Let $A: \left\{ \begin{array}{ccc}
 X & \longrightarrow & \mathcal{M}_d (\RR) \\
 x & \mapsto & A_x
\end{array} \right\}$ be a cocycle.
Assume that the cocycle is log-integrable: $\int_X \log |||A_x||| \text{d}\mu < \infty$.

There exist $\theta_1>\theta_2>...>\theta_d$, called \emph{Lyapunov exponents}, such that for almost every $x\in X$, there exists a flag $ \RR^d = \mathcal{H}_k^{(x)} \supset \mathcal{H}_{k-1}^{(x)} \supset ... \supset \mathcal{H}_1^{(x)} $ such that:

$$
A_x \mathcal{H}_j^{(x)} = \mathcal{H}_j^{(\varphi(x))}
$$

and

$$
\forall v\in\mathcal{H}_j\setminus\mathcal{H}_{j-1}, \,   \frac{1}{n}\log || A_x^{(n)} v || \underset{n\rightarrow\infty}{\longrightarrow} \theta_j
$$
where $ A_x^{(n)} = A_{\varphi^{n-1}(x)}...A_{\varphi(x)}A_{x}$.
\end{thm}

In our context we apply Oseledets' theorem in the set 
$$
X=\underset{\sigma \in \left[\pi(T)\right]}{\bigsqcup} \Delta_{N-1}\times\{\sigma\}
$$  
of all IETs whose underlying permutation is in the Rauzy class $\left[\pi(T)\right]$ of
$$ \pi(T)~=~\left(
\begin{array}{cccccc}
I_1&I_2&...&I_{N-1}\\
I_{N-1}&...&I_2&I_1
\end{array}
\right).$$
The space $X$ is equipped with $\mu_Z$, which is an ergodic $\mathfrak{R}$-invariant probability measure \cite{Zorich_finite_measure}. The cocycle $B$ introduced in last section is invertible, $B$ and $B^{-1}$ are log-integrable. We consider here $B^{-1}$. This is used in \cite{howdo} to prove Proposition \ref{Prop_Zorich_lower_bound}.

In this article we call \textbf{Oseledets generic} both an IET for which the flag decomposition exists and the corresponding uplet of its lengths $(a_1,\dots, a_{N-1})$.


%
%
%
%
%
%

\subsection{Decomposition of the symbolic coding}

At each step of renormalized Rauzy induction, there is an associated substitution.
%
For example, the substitutions associated to the example of Figure~\ref{fig:Rokhlin} are:
$$
\sigma_1 : \left\{ \begin{array}{l}
a \\
b \\
c \\
d
\end{array} \mapsto \begin{array}{l}
ad \\
b \\
c \\
d
\end{array} \right.,
\qquad
\sigma_2 : \left\{ \begin{array}{l}
a \\
b \\
c \\
d
\end{array} \mapsto \begin{array}{l}
a \\
bd \\
c \\
d
\end{array} \right.,
\qquad
\sigma_3 = \sigma_4 : \left\{ \begin{array}{l}
a \\
b \\
c \\
d
\end{array} \mapsto \begin{array}{l}
a \\
b \\
c \\
cd
\end{array} \right.,
\qquad
\sigma_5 : \left\{ \begin{array}{l}
a \\
b \\
c \\
d
\end{array} \mapsto \begin{array}{l}
a \\
b \\
cd \\
d
\end{array} \right. .
$$
The composition of the substitutions gives the substitution associated to $\mathfrak{R}^5$:
$$
\sigma = \sigma_5 \circ \sigma_4 \circ \sigma_3 \circ \sigma_2 \circ \sigma_1 : \left\{ \begin{array}{l}
a \\
b \\
c \\
d
\end{array} \mapsto \begin{array}{l}
ad \\
bd \\
cccd \\
ccd
\end{array} \right. .
$$
Its associated matrix is
$$
 M_\sigma~=~\left( \begin{array}{cccc}
1 & 0 & 0 & 0 \\
0 & 1 & 0 & 0 \\
0 & 0 & 3 & 2 \\
1 & 1 & 1 & 1
\end{array} \right).
$$

See \cite{Delecroix} for more details about the substitution associated to the Rauzy induction. We mention also that Definition 8 of \cite{ArnouxIto} proposes a more general framework.
For more details about the matrix of a substitution we refer to \cite{QueffelecChap5}.

\begin{prop}\label{decomp_SD}
Let $T$ be a self-similar IET on an interval $I$. Let $q$ be the minimal integer such that $\mathcal{R}^q (T) = T$. 
Let $\sigma$ be the substitution associated to $\mathfrak{R}^q$ applied to $T$.
We define $K=\underset{a\in\mathcal{A}}{\text{max}} \, |\sigma (a)|$. 
Let $x\in I$. 
Let $\overline{x}=x_0x_1x_2...$ be the symbolic coding of $x$ under $T$. 
Let $n\in \mathbb{N}$. We can decompose the prefix of $\overline{x}$ of length $n$ as 
$$\overline{x}_n=s_0 \sigma(s_1) ... \sigma^{l-1}(s_{l-1}) \sigma^l (m_l) \sigma^{l-1}(p_{l-1}) ... \sigma(p_1) p_0 $$ 
where 
each word $s_i$ or $p_i$ is of length at most $K-1$ and the length of $m_l$ is between $1$ and $2K-2$. 
\end{prop}
\begin{rmk}
In this decomposition, $l$ is as big as possible, the words $s_i$ (resp. $p_i$) are suffixes (resp. prefixes) of words of $\sigma(\mathcal{A})$.
\end{rmk}

For a proof of this proposition, see Lemma 5.5 in \cite{QueffelecChap5}. 

\vspace*{0.5cm}

In terms of Rokhlin towers, this means that we consider towers of order $l$ (apply the substitution once corresponds to apply the Rauzy induction once, and consider one order further on the Rokhlin tower). If the point $x$ (whose symbolic coding begins with $\overline{x}_n$) belongs to the interval on which we induce, then all $s_i$ are empty words. The orbit of $x$ under $T$ can be read on the Rokhlin tower, starting at the bottom of one tower. If not, it starts in the middle of some tower, goes up until the top of this tower and then returns at the bottom of another (or maybe the same) tower (of order $l$). The length of the word $m_l$ corresponds to the number of towers (with multiplicity) totally wandered by the first $n$ steps of the orbit of $x$. If the orbit stops in the top of a tower, all $p_i$ are empty words. Otherwise, the orbit continues to wander a part of a tower of rank $l$, and this part is the concatenation of $|p_{l-1}|$ towers of rank $l-1$, $|p_{l-2}|$ towers of rank $l-2$, ..., $|p_1|$ towers of rank $1$ and finally $|p_0|$ towers of rank $0$ (i.e. non-towers).


\subsubsection*{Bounding the number $l$ of induction steps in terms of the size $n$ of the orbit}
We recall some classical facts.
\begin{lem} \label{lem_encadrement}
Let $\sigma$ be a primitive substitution on an alphabet $\mathcal{A}$ and $\lambda$ the Perron-Frobenius eigenvalue for the matrix $M_\sigma$ of the substitution. There exist constants $A, B>0$ such that for every letter $x$ in $\mathcal{A}$ and every $n\in\NN$, $A \lambda^n < |\sigma^n(x)| < B \lambda^n$.
\end{lem}
This lemma results from basic properties of primitive substitutions that are explained in \cite{QueffelecChap5}.

\begin{cor}\label{cor_nl}
Let $T$ be an IET on an interval $I$ and $x\in I$. 
Let $\overline{x}$ be the symbolic coding of $x$ under $T$. 
Assume that $\mathcal{R}^q(T)=T$. 
Denote by $\lambda$ the Perron-Frobenius eigenvalue of $M_\sigma$ where $\sigma$ is the substitution associated to $\mathfrak{R}^q$ applied to $T$.
There exists a constant $\kappa>0$ such that the decomposition 
$$\overline{x}_n=s_0 \sigma(s_1) ... \sigma^{l-1}(s_{l-1}) \sigma^l (m_l) \sigma^{l-1}(p_{l-1}) ... \sigma(p_1) p_0$$ verifies
$$ \forall l\in\NN, \, \left| l - \frac{\log(n)}{\log(\lambda)} \right| \leqslant \kappa.$$

\end{cor}

\begin{proof}
Let $A,B$ be the constants of Lemma~\ref{lem_encadrement}, and $K=\underset{a\in\mathcal{A}}{\max}|\sigma(a)|$. We set $C =\max \left( A, \frac{(2K-2)B\lambda}{\lambda - 1}\right)$ and $\kappa=\frac{\log C}{\log \lambda}$. Lemma~\ref{lem_encadrement} gives upper and lower bounds for  $n = \sum_{k=0}^{l-1} (|\sigma^k(s_k)| + |\sigma^k(p_k)|) + |\sigma^l(m_l)|$, which permits to show the corollary.
\end{proof}

\section{Estimation of the trajectory} \label{section_approx_traj} 
\subsection{Lengths of shifts}
The assumption (C1) implies that one out of two crossed sides is the side \textbf{$a_N$} (the biggest one).
Let us remark that every other polygon along the trajectory is a translation of the first one.  
Therefore we can approach two steps ($i$ and $i+1$) of the trajectory by a translation: the one that sends the polygon $P_i$ to the polygon $P_{i+2}$. 

\begin{figure}[h!]
     \centering
\begin{subfigure}[b]{0.35\textwidth}
\centering
\definecolor{qqqqff}{rgb}{0,0,1}
\definecolor{uququq}{rgb}{0.25,0.25,0.25}
\definecolor{zzttqq}{rgb}{0.6,0.2,0}
\definecolor{xdxdff}{rgb}{0.49,0.49,1}
\begin{tikzpicture}[line cap=round,line join=round,>=stealth,x=1.cm,y=1.cm]
\clip(-1,-1.8) rectangle (2.9,1.1);
\fill[color=zzttqq,fill=zzttqq,fill opacity=0.1] (1,0.1) -- (0.92,-0.39) -- (0.77,-0.64) -- (-0.76,-0.65) -- (-0.96,-0.28) -- cycle;
\fill[color=zzttqq,fill=zzttqq,fill opacity=0.1] (0.92,-0.39) -- (1,0.1) -- (1.15,0.35) -- (2.68,0.36) -- (2.88,-0.01) -- cycle;
\fill[color=zzttqq,fill=zzttqq,fill opacity=0.1] (2.88,-0.01) -- (2.8,-0.49) -- (2.65,-0.75) -- (1.12,-0.75) -- (0.92,-0.39) -- cycle;
\draw [color=zzttqq] (1,0.1) -- (0.92,-0.39) -- (0.77,-0.64) -- (-0.76,-0.65) -- (-0.96,-0.28) -- cycle;
\draw [color=zzttqq] (0.92,-0.39) -- (1,0.1) -- (1.15,0.35) -- (2.68,0.36) -- (2.88,-0.01) -- cycle;
\draw [color=zzttqq] (2.88,-0.01) -- (2.8,-0.49) -- (2.65,-0.75) -- (1.12,-0.75) -- (0.92,-0.39) -- cycle;
\draw (-0.51,-0.19)-- (0.97,-0.1);
\draw (0.97,-0.1)-- (1.46,-0.28);
\draw (1.46,-0.28)-- (1.06,-0.65);
\draw [color=qqqqff,->,>=stealth,thick] (1,0.1)-- (2.88,-0.01);
\begin{scriptsize}
\draw [color=qqqqff] (2.6,0.2) node {$f_1$};
\end{scriptsize}
\end{tikzpicture}
\caption{$f_1$}
\end{subfigure}
\hspace{-1cm}
\begin{subfigure}[b]{0.35\textwidth}
\centering
\definecolor{qqqqff}{rgb}{0,0,1}
\definecolor{zzttqq}{rgb}{0.6,0.2,0}
\definecolor{xdxdff}{rgb}{0.49,0.49,1}
\begin{tikzpicture}[line cap=round,line join=round,>=stealth,x=1cm,y=1cm]
\clip(-2.8,-1.6) rectangle (2.9,1.1);
\fill[color=zzttqq,fill=zzttqq,fill opacity=0.1] (1,0.1) -- (0.92,-0.39) -- (0.77,-0.64) -- (-0.76,-0.65) -- (-0.96,-0.28) -- cycle;
\fill[color=zzttqq,fill=zzttqq,fill opacity=0.1] (0.69,-1.13) -- (0.77,-0.64) -- (0.92,-0.39) -- (2.45,-0.38) -- (2.65,-0.75) -- cycle;
\fill[color=zzttqq,fill=zzttqq,fill opacity=0.1] (2.65,-0.75) -- (2.58,-1.23) -- (2.42,-1.49) -- (0.89,-1.49) -- (0.69,-1.13) -- cycle;
\draw [color=zzttqq] (1,0.1) -- (0.92,-0.39) -- (0.77,-0.64) -- (-0.76,-0.65) -- (-0.96,-0.28) -- cycle;
\draw [color=zzttqq] (0.69,-1.13) -- (0.77,-0.64) -- (0.92,-0.39) -- (2.45,-0.38) -- (2.65,-0.75) -- cycle;
\draw [color=zzttqq] (2.65,-0.75) -- (2.58,-1.23) -- (2.42,-1.49) -- (0.89,-1.49) -- (0.69,-1.13) -- cycle;
\draw (-0.5,-0.19)-- (0.89,-0.44);
\draw (0.89,-0.44)-- (1.32,-1);
\draw (1.32,-1)-- (1.19,-1.49);
\draw [color=qqqqff,->,>=stealth,thick] (1,0.1)-- (2.65,-0.75);
\begin{scriptsize}
\draw [color=qqqqff] (2.3,-0.7) node {$f_2$};
\end{scriptsize}
\end{tikzpicture}
\caption{$f_2$}
\end{subfigure}
\begin{subfigure}[b]{0.35\textwidth}
\centering
\definecolor{qqqqff}{rgb}{0,0,1}
\definecolor{zzttqq}{rgb}{0.6,0.2,0}
\definecolor{xdxdff}{rgb}{0.49,0.49,1}
\begin{tikzpicture}[line cap=round,line join=round,>=stealth,x=1cm,y=1cm]
\clip(-1,-1.8) rectangle (2.9,1.1);
\fill[color=zzttqq,fill=zzttqq,fill opacity=0.1] (1,0.1) -- (0.92,-0.39) -- (0.77,-0.64) -- (-0.76,-0.65) -- (-0.96,-0.28) -- cycle;
\fill[color=zzttqq,fill=zzttqq,fill opacity=0.1] (-0.99,-1.39) -- (-0.92,-0.9) -- (-0.76,-0.65) -- (0.77,-0.64) -- (0.97,-1.01) -- cycle;
\fill[color=zzttqq,fill=zzttqq,fill opacity=0.1] (0.97,-1.01) -- (0.89,-1.49) -- (0.74,-1.75) -- (-0.79,-1.76) -- (-0.99,-1.39) -- cycle;
\draw [color=zzttqq] (1,0.1) -- (0.92,-0.39) -- (0.77,-0.64) -- (-0.76,-0.65) -- (-0.96,-0.28) -- cycle;
\draw [color=zzttqq] (-0.99,-1.39) -- (-0.92,-0.9) -- (-0.76,-0.65) -- (0.77,-0.64) -- (0.97,-1.01) -- cycle;
\draw [color=zzttqq] (0.97,-1.01) -- (0.89,-1.49) -- (0.74,-1.75) -- (-0.79,-1.76) -- (-0.99,-1.39) -- cycle;
\draw (-0.17,-0.13)-- (0.41,-0.64);
\draw (0.41,-0.64)-- (-0.25,-1.25);
\draw (-0.25,-1.25)-- (0.8,-1.64);
\draw [color=qqqqff,->,>=stealth,thick] (1,0.1)-- (0.97,-1.01);
\begin{scriptsize}
\draw [color=qqqqff] (1.15,-0.9) node {$f_3$};
\end{scriptsize}
\end{tikzpicture}
\caption{$f_3$}
\end{subfigure}
\hspace{-1cm}
\begin{subfigure}[b]{0.35\textwidth}
\centering
\definecolor{qqqqff}{rgb}{0,0,1}
\definecolor{zzttqq}{rgb}{0.6,0.2,0}
\definecolor{xdxdff}{rgb}{0.49,0.49,1}
\begin{tikzpicture}[line cap=round,line join=round,>=stealth,x=1cm,y=1cm]
\clip(-2.8,-1.8) rectangle (2.8,1.46);
\fill[color=zzttqq,fill=zzttqq,fill opacity=0.1] (1,0.1) -- (0.92,-0.39) -- (0.77,-0.64) -- (-0.76,-0.65) -- (-0.96,-0.28) -- cycle;
\fill[color=zzttqq,fill=zzttqq,fill opacity=0.1] (-2.72,-1.03) -- (-2.64,-0.54) -- (-2.49,-0.29) -- (-0.96,-0.28) -- (-0.76,-0.65) -- cycle;
\fill[color=zzttqq,fill=zzttqq,fill opacity=0.1] (-0.76,-0.65) -- (-0.83,-1.14) -- (-0.99,-1.39) -- (-2.52,-1.4) -- (-2.72,-1.03) -- cycle;
\draw [color=zzttqq] (1,0.1) -- (0.92,-0.39) -- (0.77,-0.64) -- (-0.76,-0.65) -- (-0.96,-0.28) -- cycle;
\draw [color=zzttqq] (-2.72,-1.03) -- (-2.64,-0.54) -- (-2.49,-0.29) -- (-0.96,-0.28) -- (-0.76,-0.65) -- cycle;
\draw [color=zzttqq] (-0.76,-0.65) -- (-0.83,-1.14) -- (-0.99,-1.39) -- (-2.52,-1.4) -- (-2.72,-1.03) -- cycle;
\draw (0.56,0.01)-- (-0.85,-0.49);
\draw [color=qqqqff,->,>=stealth,thick] (1,0.1)-- (-0.76,-0.65);
\draw (-0.85,-0.49)-- (-1.16,-0.73);
\draw (-1.16,-0.73)-- (-0.79,-0.82);
\begin{scriptsize}
\draw [color=qqqqff] (-0.2,-0.55) node {$f_4$};
\end{scriptsize}
\end{tikzpicture}
\caption{$f_4$}
\end{subfigure}
\caption{Translations $f_i$ when crossing sides \textbf{$a_i$} followed by side \textbf{$a_5$}}\label{fig:f_i}
\end{figure}
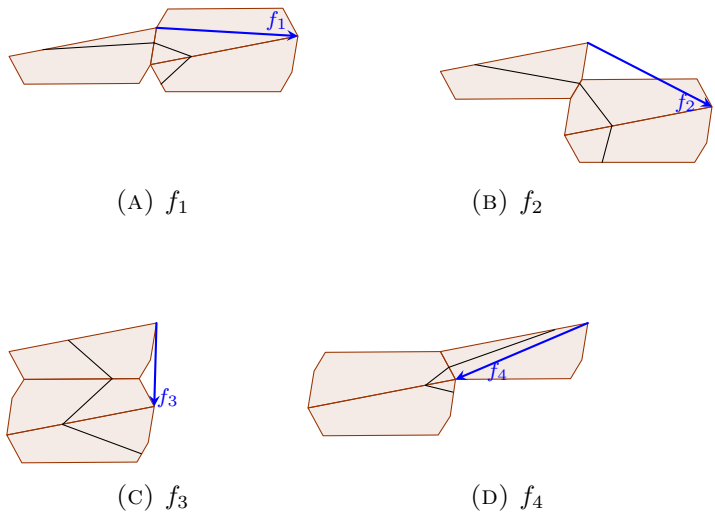

We identify the space of translation with $\CC$. 
Let us denote by $f_i$ 
the translation induced by crossing the side $\textbf{a}_\textbf{i}$ and then the side $\textbf{a}_\textbf{N}$ (see Figure~\ref{fig:f_i}).
We recall that $x_0$ denotes the position of the vertex $v_0$ of $P_0$ in $\mathcal{C}_0$. 
One has:
$$ \forall k\in\{1,..,N-1\}, \, f_k = r \left(e^{X_k} - e^{X_N} + e^{X_{k-1}} - e^{X_0}\right),$$ where
\label{Xi_more}
$$ X_k = -\frac{i}{r} \left( x_0 + \sum_{j=1}^{k} a_j \right) \quad \text{and} \quad r = \frac{1+a_{N}}{2\pi} \text{ is the radius of the circumcircle of } P.
$$


\subsection{Estimate the trajectory}
Let us define the piecewise constant complex function $f$ on the interval $[0,1)$ by:
$$f(x) =    f_i  \quad \text{if } \sum_{k=1}^{i-1} a_k \leqslant x < \sum_{k=1}^{i} a_k .
$$
After $2n$ steps in the plane, the trajectory is in the polygon $P_{2n}$, polygon which is equal to the first one ($P_0$) translated by
$$S_n f(x_0)= \sum_{k = 0}^{n-1} f(T^k (x_0)).$$
This sum is the \emph{Birkhoff sum} of $f$ with respect to $T$ and $x_0$.
We approach the trajectory with this value.
We highlight that we only need to know how many times the orbit $(T^k x)_{k\leqslant n}$ goes into each interval $I_i \subset (0,1)$ to compute Birkhoff sum.
The sequence $(S_n f(x_0))_{n\in\NN}$ of Birkhoff sums gives us (with uniformly bounded error) an estimate of the position $p_{2n}$ (of the $2n$-th side crossing of the trajectory) 
and hence also of $p_{2n+1}$.

Note that $f$ depends on $x_0$. However, a change of $x_0$ will only multiply $f$ by a complex number. Therefore all the translations $f_i$ depend on $x_0$ (i.e. on the position of $P_0$) only up to the same rotation.
Whether the trajectory has a limit or not hence does not depend on $x_0$. 
Nor the fact that the trajectory stays at bounded distance of a line. Only the asymptotic direction (if it exists) depends on $x_0$.



\subsection{Average shift}\label{subsec_average-shift}
We first show that the trajectory has almost always an asymptotic direction. 


\AsymptoticDir*

\begin{proof}
 
As soon as $T$ is uniquely ergodic (which is the generic case \cite{Masur}, \cite{Veech}), Birkhoff's theorem tells us that the main part of this sum will be $n$ times the mean value $m=\sum_{k=1}^{N-1} a_k f_k$ of $f$: 
$$ \forall x_0 \quad S_nf(x_0) \underset{n\rightarrow \infty}{=} n\sum_{k=1}^{N-1} a_k f_k + o(n).$$ 


Hence, in terms of the trajectory, $m$ corresponds to the asymptotic direction of the trajectory, as soon as it does not vanish. Let us show that this is almost always the case.

We have:
\begin{align*}
m &= \sum_{k=1}^{N-1} a_k f_k \\
 &= re^{\frac{i}{r}x_0}
 \left[ \sum_{k=1}^{N-2} (a_k+a_{k+1}) e^{\frac{i}{r}\sum_{j=1}^k a_j}
  - \left(\sum_{k=1}^{N-2} a_k \right) e^{\frac{i}{r}\sum_{j=1}^{N-1} a_j} 
  - \left(\sum_{k=2}^{N-1}a_k\right) \right].
\end{align*}

Let us denote $s=\frac{2\pi}{1+a_N}=\frac{1}{r}$ the inverse of the radius of the circumcircles $\mathcal{C}_i$. If $a_1,...,a_{N-1}$ are fixed, the map $s \in [0,\pi] \mapsto sm$ is analytic. As a non trivial sum of linearly independant functions ($s\mapsto e^{is\sum_{j=1}^{k}a_j}$ for $1\leqslant k\leqslant N$), it cannot be everywhere equal to zero.
Therefore the map has a finite number of zeros in the compact set $[0,\pi]$. We conclude that only a finite number of values of $a_N\in(1,+\infty)$ could make $m$ vanish.
\end{proof}

\subsection{Normalized deviations}
Our goal is to study the deviations from the mean displacement. We want to estimate their size and their direction, in particular we want to show that they are happening in the direction orthogonal to $m$.
To do so, we renormalize the situation so that the trajectory goes asymptotically in the horizontal direction at unit speed, namely we rotate and contract (or expand) the plane by the complex coefficient $\frac{1}{m}$. Now we have to show that the deviations have an unbounded imaginary part.

\begin{lem}
We set $G = \frac{1}{m}\left( \begin{array}{c}
     f_1\\
 \vdots \\ 
     f_{N-1}
\end{array}\right) - \left( \begin{array}{c}
     1\\
     \vdots\\
     1
\end{array}\right)  $ and  $H = \Im(G)$. Let $h: (0, 1) \longrightarrow \RR $ be the map that is constant and equal to $H_i$ on the subinterval $(\sum_{j=1}^{i} a_j, \sum_{j=1}^{i+1} a_j)$ for every $0 \leqslant i \leqslant N-1$.  
The following statements are equivalent:
\begin{enumerate}
\item
The trajectory in the plane admits unbounded deviations in the direction orthogonal to the asymptotic one.
\item
The trajectory does not stay at a bounded distance from a line.
\item
The sequence $(S_n(h))_{n\in\NN}$ is not bounded.
\end{enumerate}
Moreover, if $T$ is Osseledets-generic they are also equivalent to:
\begin{enumerate}[resume]
\item
The vector $H$ is not in the contracted plane of the flag decomposition of $\RR^{N-1}$ given by the Oseledets' theorem. 
\end{enumerate}
\end{lem}

\begin{proof}
The first two statements are clearly equivalent. 

Quantifying deviations consists in computing the second term of the asymptotic expansion of $S_n f(x_0)$, namely studying $S_n f(x_0)- nm$. We renormalize the situation by dividing by $m$: $\frac{1}{m}S_n f(x_0)- n $. This corresponds to contract the plane by a factor of $\frac{1}{m}$, the renormalized trajectory admits the real line as an asymptotic direction. Since we are interested in its deviations from this asymptotic line, we want to compute the imaginary part of the second term in the asymptotic expansion of $S_n f(x_0)$: $$ \Im \left(\frac{1}{m}S_n f(x_0)- n \right) = S_n h(x_0).$$ 
The function $h$ corresponds to the deviations from the renormalized mean displacement (equal to $(1,0)\in\RR^2$). 

If $H$ is in the contracted plane of the Oseledets' flag, then the sum converges (as a geometric sum), i.e. $S_n h(x_0)$ is bounded, and the trajectory stays at a bounded distance from the asymptotic line. On the contrary, if $H$ is not in the contracted plane of the Oseledets' flag, then the sum diverges (see Proposition \ref{Prop_Zorich_lower_bound}), i.e. $S_n h(x_0)$ is not bounded, and the trajectory does not stay at a bounded distance from the asymptotic line.
\end{proof} 

The following property of vector $G$ will be useful.
\begin{lem}
The vector $G$ is normal to the vector $V = \left( \begin{array}{c}
     a_1\\
 \vdots \\ 
     a_{N-1}
\end{array}\right)$. 
\end{lem}

\begin{proof}
This is a direct computation:
$$
 \prescript{t}{}{V} G = \frac{1}{m} \underset{=m}{\underbrace{\sum_{k=1}^{N-1} a_k f_k}} - \underset{=1}{\underbrace{\sum_{k=1}^{N-1} a_k}} = 0. 
 $$
\end{proof}

We show in the next section that in the generic case we can apply Zorich's result on deviations for IETs. This gives Theorem~\ref{th_generic} . We handle the non generic case of a self-similar IET in Section~\ref{section_selfsim}: one has $ B_T^{(nq)}= \left(B_T^{(q)}\right)^n$. Direct computations on $H$ and $M=B_T^{(q)}$ allow to prove Theorem~\ref{th_selfsim}.

\section{The generic case}
\label{section_G}


We assume here that $T$ is Oseledets-generic, which means (see Subsection \ref{subsec_Oseledets}) that there exists a flag depending on $T$ 
$$\mathcal{H}_1 \subset \mathcal{H}_2 \subset ... \subset \mathcal{H}_{\textbf{g}(T)} \subset \RR^{N-1} $$ 
with $2\textbf{g}(T)= \left\{ \begin{array}{ll} 
N-1 & \text{if } N \text{ is odd} \\ 
N-2 & \text{if } N \text{ is even}
\end{array} \right.$, 
such that
$$\forall v\in\mathcal{H}_k\setminus\mathcal{H}_{k-1}, \, \frac{1}{n}\log||B_T^{(n)}(v)|| \underset{n\rightarrow\infty}{\longrightarrow} \theta_k$$
where the numbers $\theta_1>\dots>\theta_{\textbf{g}(T)}$ do not depend on $T$. The $1+\frac{\theta_i}{\theta_1}$ are the Lyapunov exponents of the Teichmüller flow on the hypereliptic stratum ($\mathcal{H}(N-3)$ if $N$ is odd, (resp. $\mathcal{H}(\frac{N}{2}-2,\frac{N}{2}-2)$ if $N$ is even) of the moduli space of abelian differentials. 
%
%
We will apply Proposition 3 of \cite{howdo}:
\begin{prop}[Zorich]\label{Prop_Zorich_lower_bound}
For a generic IET on $(0,1)$, any point $x\in(0,1)$, and any function $\varphi \in Ann(\mathcal{H}_p) \setminus Ann(\mathcal{H}_{p+1})$, $1 \leqslant p \leqslant \emph{\textbf{g}}(T)-1$, the following holds: 
$$ \underset{n\rightarrow + \infty}{\limsup} \, \frac{\log|S_n(\varphi(x))|}{\log n} = \frac{\theta_{p+1}}{\theta_1}. $$
\end{prop}

We already know that the function $h$ is in $Ann(\mathcal{H}_1)$ (because it has $0$ mean). In this section we will show that it is not in $Ann(\mathcal{H}_3)$ for almost every choice of parameters (see Proposition~\ref{not_in_fixed_plane_N}). This will establish Theorems~\ref{th_dev} and~\ref{th_generic}. 

To do that, we give an expression of $H$. It involves a matrix $Q$ whose size increases with $N$: $Q\in\mathcal{M}_{\frac{(N-1)N}{2},N-1}(\RR)$. This is why we do the proof first for $N=5$ (in Sections \ref{subsec_G_5} and \ref{subsec_never_fixed_plane_5}) and then, with the very same proof but a bigger matrix, for $N>5$ (in Section \ref{subsec_G_never_fixed_plane_N}) .


\subsection{Expression of $H$, when $N=5$}\label{subsec_G_5}
Let us compute 
$$
H = Im\left(\frac{1}{m}
\left( \begin{array}{c}
     f_1\\
	f_2 \\
     f_3 \\
     f_{4}
\end{array}\right)
- \left( \begin{array}{c}
     1 \\
     1 \\  
     1 \\
     1 
\end{array} \right)\right)
=\frac{1}{|m|^2}Im\left(\overline{m}
\left( \begin{array}{c}
     f_1\\
	f_2 \\
     f_3 \\
     f_{4}
\end{array}\right)
\right)
$$
where $\overline{m}$ denotes the complex conjugate of $m$.
One can check that, with  $s=\frac{2\pi}{1+a_5}$,
$$
H = \frac{1}{|m|^2} Q \left( \begin{array}{c}
   \sin(a_1s) \\  
   \sin(a_2s) \\  
   \sin(a_3s) \\  
   \sin(a_4s) \\  
   \sin((a_1+a_2)s) \\  
   \sin((a_2+a_3)s)\\  
   \sin((a_3+a_4)s) \\  
   \sin((a_1+a_2+a_3)s) \\  
   \sin((a_2+a_3+a_4)s)\\ 
    \sin((a_1+a_2+a_3+a_4)s)     
    \end{array} \right)
    = H(a_1,a_2,a_3,a_4,s)
$$
where
$$
\prescript{t}{}{Q} = \left( \begin{array}{cccc}
   a_2+a_3+a_4 & -a_1+a_3+a_4 & -(a_1+a_2)  & -(a_1+a_2)\\
   a_2+a_3 & -a_1+a_3 & -(a_1+a_2) & 0 \\
    0 & a_3+a_4 & -a_2+a_4  & -(a_2+a_3)\\
    a_3+a_4 & a_3+a_4 & -(a_1+a_2) +a_4 & -(a_1+a_2+a_3)\\
    0 & a_4 & a_4 & -(a_2+a_3) \\
   a_3+a_4 & a_3+a_4 & -(a_1+a_2)& -(a_1+a_2) \\
    a_2+ a_3 & -a_1 & -a_1  & 0 \\
    0 & -(a_3+a_4) & a_2 & a_2\\
    -a_3 & -a_3 & a_1+a_2 & 0 \\
  -(a_2+a_3+a_4) & a_1-a_4 &  a_1-a_4 & a_1+a_2+a_3
      \end{array} \right).
$$

\subsection{Vector $H$ is almost never in a fixed plane, when $N=5$}\label{subsec_never_fixed_plane_5}

 

We show that $H$ does not stay in a fixed plane of $\RR^4$ when $s$ varies.
\begin{prop}\label{not_in_fixed_plane}
Let $\mathcal{P}$ be a plane in $\RR^4$. Let $a_1,a_2,a_3,a_4 \in \RR^+$ summing up to $1$.
Then the set $\{s\in(0,\pi) \ | \ H(a_1,a_2,a_3,a_4,s) \in \mathcal{P}\}$ is finite.
\end{prop}

We first show:
\begin{lem}
The matrix $Q$ has rank $3$.
\end{lem}

\begin{proof}
Since the vector $V = \ \prescript{t}{}{(\begin{array}{cccc}
     a_1 & a_2 & a_3 & a_4 
\end{array})}$ is normal to $G$, and $H = Im(G)$,
the relation
$\sum_{i=1}^{4}a_i R_i(Q) = 0$  holds, where $R_i(Q)$ denotes the $i$-th row of $Q$ (i.e. the $i$-th column of $\prescript{t}{}{Q}$). It means that $Q$ has rank at most $3$.

Moreover the minor of $\prescript{t}{}{Q}$ formed by the first three columns and the fourth, fifth and eighth rows of $\prescript{t}{}{Q}$ is:
$$\left| \begin{array}{ccc}
     a_3+a_4 & a_3+a_4 & -(a_1+a_2)+a_4 \\
     0 & a_4 & a_4 \\
     0 & -(a_3+a_4) & a_2 \\
\end{array}
\right| = (a_3+a_4)a_4(a_2+a_3+a_4) > 0.$$
The matrix $Q$ has hence rank at least $3$.
\end{proof}

\begin{proof}[Proof of Proposition~\ref{not_in_fixed_plane}]
The map 
$$\mathbf{H} : \begin{array}{ccc}
     [0,\pi] & \longrightarrow & \RR^4 \\
     s & \mapsto & H(a_1,a_2,a_3,a_4,s)
\end{array}$$ is analytic.
If $\mathbf{H}(s)$ were lying in $\mathcal{P}$ for all $s$ in a sub-interval $J$ of $[0,\pi]$, then we would have two linearly independent relations
$$
\left\{\begin{array}{c}
 \sum_{i} \mu_i \sin (\theta_i s) = 0  \\
\sum_{i} \lambda_i \sin (\theta_i s) = 0 
\end{array}
\right. 
$$
where the $\theta_i$'s are sums of $a_j$ with consecutive indices.
Since $G$ is normal to $v$, the first relation holds if $\mu = v$. But $Q$ has rank $3$ and $s\mapsto \sin ((a_1+a_2+a_2+a_3+a_4)s)=sin(s)$ is linearly independent of all other maps $s\mapsto \sin(\theta s)$ for $\theta \neq \pm 1$. Therefore the second relation cannot hold for all $s\in J$.

Moreover, the zeroes of the analytic map $s \in [0,\pi] \mapsto\sum_{i} \lambda \sin (\theta_i s)$ are discrete, and in finite number in the compact $[0,\pi]$. 
This proves Proposition~\ref{not_in_fixed_plane}

\end{proof}

\subsection{For $N>5$} \label{subsec_G_never_fixed_plane_N}

As in the case $N=5$, we consider $H$ the imaginary part of $G$:
$$H = Im\left(\frac{1}{m}
\left( \begin{array}{c}
     f_1\\
 \vdots \\ 
     f_{N-1}
\end{array}\right) 
- \left( \begin{array}{c}
     1 \\
     \vdots \\
     1 
\end{array} \right)\right) = \frac{1}{|m|^2}Im\left(\overline{m}
\left( \begin{array}{c}
     f_1\\
 \vdots \\
     f_{N-1}
\end{array}\right) 
\right).$$
All coordinates of $|m|^2 H$ are of the form $\sum \lambda_i \sin (\theta_i s)$ where $s=\frac{1}{r}$, $\theta_i$
are sums of $a_j$ with consecutive indices, and $\lambda_i$ are sums of $\pm a_j$.
We can write $|m|^2H$ as $|m|^2H = Q \Theta$ where the matrix $Q\in \mathcal{M}_{\frac{(N-1)N}{2},N-1}(\RR)$ has all of its coefficients in the set $\{0, \lambda_i\}$ and the vector $\Theta \in \mathcal{M}_{  \frac{(N-1)N}{2},1}$ has all of its coordinates of the form $\sin(\theta_i s)$. 
We order this coordinates by lexical order on (number of terms $a_j$ within the sum $\theta_i = \sum_j a_j$ ; number of the first index $j$), so that:
$$
\Theta
= \left( \begin{array}{c}
   \sin{(a_1 s)} \\  
   \vdots \\  
   \sin{(a_{N-1} s)} 
      \vspace{0.15cm} \\ 
   \sin{((a_1+a_2)s)} \\  
   \vdots \\  
   \sin{((a_{N-2}+a_{N-1})s)}
   \vspace{0.15cm}\\
   \sin{((a_1+a_2+a_3)s)} \\  
   \vdots \\ 
   \sin{((a_{N-3}+a_{N-2}+a_{N-1})s)}
   \vspace{0.15cm} \\
   \vdots \\
   \vdots 
   \vspace{0.15cm} \\
   \sin{((a_1+a_2+...+a_{N-2})s)} \\  
   \sin{((a_2+ a_3+ ...+a_{N-1})s)}
   \vspace{0.15cm}\\
    \sin{((a_1+a_2+...+a_{N-1})s)}
      \end{array} \right).
$$

One can check that lots of coordinates of $Q$ are zero.
In the following we only need to consider a submatrix of $Q$.

\begin{lem}\label{lem_rank_Q_N}
The matrix $Q\in\mathcal{M}_{\frac{(N-1)N}{2},N-1}(\RR)$ has rank at least $N-3$, at most $N-2$, and exactly $N-2$ when $N$ is odd.
\end{lem}

\begin{proof}
As in the case $N=5$, $G$ is normal to the vector $V=\ \prescript{t}{}{(\begin{array}{cccc}
     a_1 & \dots & a_{N-1}
\end{array})}$, which implies that $Q$ has rank at most $N-2$.

The submatrix of $\prescript{t}{}{Q}$ made of the first $N-1$ rows and $N-2$ columns is:
  $$
  \left(
     \raisebox{0.5\depth}{%
       \xymatrixcolsep{0.5ex}%
       \xymatrixrowsep{0.5ex}%
       \xymatrix{
      &  \sum_{k=2}^{N-1}a_k & -a_1+ \sum_{k=3}^{N-1}a_k & -(a_1+a_2)  \ar@{.}[rrrrrr] & & & & & & -(a_1+a_2)\\
   \ar@{-}@[gray][rrrrrrrr] \ar@{-}@[gray][dddddd] & & & & & & & & \ar@{-}@[gray][dddddd] &\\      
       &  a_2+a_3 & -a_1+a_3 & -(a_1+a_2) \ar@{.}[rrrrrrddd] & 0  \ar@{.}[rrrrr]  \ar@{.}[rrrrrdd] & & & & & 0\ar@{.}[dd] \\
        & 0 \ar@{.}[ddd] \ar@{.}[rrrrrddd]& a_3+a_4 & -a_2+a_4 \ar@{.}[rrrrdd] & & & & & & \\
         & & & a_4+a_5 \ar@{.}[rrrrdd] & & & & & & 0 \\
         & & & & & & & -a_{N-4}+a_{N-2} & & -a_{N-4}-a_{N-3}\\
        & 0\ar@{.}[rrrrr] & & & & & 0 & a_{N-2}+a_{N-1} & & -a_{N-3}+a_{N-1} \\
        \ar@{-}@[gray][rrrrrrrr] & & & & & & & & & \\
        & a_{N-2}+a_{N-1} \ar@{.}[rrrrrr] & & & & & & a_{N-2}+a_{N-1} & & -(\sum_{k=2}^{N-3}a_k)+a_{N-1}
       }%
     }
   \right)
  $$

Note that the framed submatrix formed by the rows $2,..,N-2$ and first $N-3$ columns of $\prescript{t}{}{Q}$ is an invertible diagonal matrix, so $Q$ has rank at least $N-3$.

We show now that if $N$ is odd, $Q$ has rank at least $N-2$.
We replace the $j$-th column $C_j$ of $\prescript{t}{}{Q}$ by $C_j + \sum_{k=1}^{j-1} (-1)^{j-k} C_k$. The rows between $2$ and $N-1$ and the first $N-2$ columns of this new matrix are equal to (changes are written in blue): 

  $$
  \left(
     \raisebox{0.5\depth}{%
       \xymatrixcolsep{0.5ex}%
       \xymatrixrowsep{0.5ex}%
       \xymatrix{
         a_2 + a_3 & -a_1-a_2 & \textcolor{blue}{0}  & 0  \ar@{.}[rrr]  \ar@{.}[rrrdd] & & & 0\ar@{.}[dd] \\
         0 \ar@{.}[ddd] \ar@{.}[rrrrddd]& a_3+a_4 & -a_2\textcolor{blue}{-a_3}  & \textcolor{blue}{0} \ar@{.}@[blue][rrrdd] &  & & \\
          & & a_4+a_5 \ar@{.}[rrrdd] & -a_3\textcolor{blue}{-a_4} \ar@{.}@[blue][rrd]& & & 0 \\
          &   & & & & -a_{N-4}\textcolor{blue}{-a_{N-3}} & \textcolor{blue}{0}\\
        0\ar@{.}[rrrr] & & & & 0 & a_{N-2}+a_{N-1} & -a_{N-3}\textcolor{blue}{-a_{N-2}} \\
        a_{N-2}+a_{N-1} & \textcolor{blue}{0} & a_{N-2}+a_{N-1}  \ar@{.}@[blue][rrr] & & & \textcolor{blue}{0} & -(\sum_{k=2}^{N-3}a_k)+a_{N-1}
       }%
     }
   \right)
  $$
  
We develop the determinant with respect to the last row of this submatrix and get:

\begin{align*}
& &(a_{N-2}+a_{N-1})[&(-a_1-a_2)(-a_2-a_3)...(-a_{N-3}-a_{N-2}) \\
&&&+(a_2+a_3)(a_3+a_4)(-a_3-a_4)...(-a_{N-3}-a_{N-2}) \\
&&&+(a_2+a_3)(a_3+a_4)(a_4+a_5)(a_5+a_6)(-a_5-a_6)...(-a_{N-3}-a_{N-2}) \\
&&&+(a_2+a_3)(a_3+a_4)...(a_7+a_8)(-a_7-a_8)...(-a_{N-3}-a_{N-2}) \\
&&&\vdots\\
&&&+(a_2+a_3)(a_3+a_4)...(a_{N-4}+a_{N-3})(-a_{N-3}-a_{N-2})] \\
&&+\left(-\sum_{k=1}^{N-3} a_k+a_{N-1}\right)&(a_2+a_3)(a_3+a_4)...(a_{N-3}+a_{N-2})(a_{N-2}+a_{N-1}) \\
&=& (a_{N-2}+a_{N-1})&\left[ (-1)^{N-3}\prod_{k=1}^{N-3} (a_k +a_{k+1}) + \sum_{\underset{l \text{ odd}}{l=3}}^{N-3} \left((a_l +a_{l+1}) (-1)^{N-3-l+1}\prod_{k=2}^{N-3} (a_k +a_{k+1})\right) \right] \\
&&+\left(-\sum_{k=1}^{N-3} a_k+a_{N-1}\right)&  \prod_{k=2}^{N-2} (a_k +a_{k+1})\\ 
&=& \prod_{k=2}^{N-2} (a_k+a_{k+1}) &\left[ \sum_{k=1}^{N-3}a_k - \sum_{k=1}^{N-3}a_k +a_{N-1} \right] \\
&=& a_{N-1}\prod_{k=2}^{N-2}(a_k+a_{k+1}) \neq 0.
\end{align*} 

We conclude that if $N$ is odd, then $Q$ has rank exactly $N-2$.

\end{proof}

As for the case $N=5$, the map
$$\mathbf{H} : \begin{array}{ccc}
     [0,\pi] & \longrightarrow & \CC^{N-1} \\
     s & \mapsto & H(a_1,...,a_{N-1},s) = Q \Theta 
\end{array}$$ is analytic.  With Lemma \ref{lem_rank_Q_N}, it allows us to prove the following proposition.

\begin{prop}\label{not_in_fixed_plane_N}
Let $\mathcal{P}$ be a subspace of dimension $N-3$ in $\CC^{N-1}$. Let $a =(a_1,...,a_{N-1}) \in (\RR_+^*)^{N-1}$ summing up to $1$ and Oseledets generic.
The set $\{s\in(0,\pi) \ | \ H(a,s) \in \mathcal{P}\}$ is finite.
Moreover, if $N$ is odd, it holds also for $\mathcal{P}$ a subspace of dimension $N-2$ in $\CC^{N-1}$.
\end{prop}

%
\subsection{Conclusion}
As a corollary we get Theorems~\ref{th_dev} and~\ref{th_generic}.
\ThDev*

\ThGeneric*
%
\ThDev*

\begin{proof}
Recall that the function $g$ is in $Ann(\mathcal{H}_1)$ (because it has $0$ mean).
Let $\mathcal{P} = Ann(\mathcal{H}_3)$ (or $Ann(\mathcal{H}_2)$ if $Q$ has rank $N-2$).
Proposition~\ref{not_in_fixed_plane_N} tells that for  almost every choice of parameters, $h$ is not in $\mathcal{P}$.
The proposition \ref{Prop_Zorich_lower_bound} applied to $h$ hence establishes Theorem~\ref{th_generic}:
\begin{align*}
|p_{2n} - 2nm| &= |S_n h(x_0)| + R_n  \quad \text{ where } R_n \text{ is uniformely bounded} \\
\underset{n\rightarrow + \infty}{\limsup} \frac{\log|p_{2n} - 2nm|}{\log(2n)} & = \underset{n\rightarrow + \infty}{\limsup} \frac{\log|S_n h(x_0)|}{\log(n)} \\
& = \frac{\theta}{\theta_2}
\end{align*}
where $\theta=\theta_2$ if $N$ is odd, and $\theta\in\{\theta_2,\theta_3\}$ if $N$ is even.

Moreover we have renormalized our problem to the case when the asymptotic direction is horizontal and we have studied the imaginary part of the deviations, this hence proves also Theorem~\ref{th_dev}.
\end{proof}

%

\section{A special case} \label{section_selfsim} 
The previous theorems tells nothing when our IET $T$ is not Oseledets generic. The goal of this section is to study a non generic  but important case for which we manage to do explicit computations, and prove Theorem \ref{th_selfsim}.

Let $M$ be the matrix associated to a loop in the Rauzy diagram. Then $M$ is positive. We assume that the loop is complex enough so that $M$ is primitive.
We assume that the second eigenvalue $\lambda_2$ of $M$ is simple, that $M$ has no other eigenvalue of the same modulus than $\lambda_2$, and $|\lambda_2|>1$. The article \cite{Hamenstadt} implies that this is the generic behavior. 
It follows that $\lambda\in\RR$.
We denote by $V$ the Perron Frobenius vector for $M$, normalized so that its coordinates sum up to 1 and set $(a_1,...,a_{N-1})=V$.
The IET with lengths $(a_1,...,a_{N-1})$ is self similar and uniquely ergodic. Section \ref{subsec_average-shift} applies: the mean $m$ of $f$ (corresponding to the asymptotic direction) does not vanish for almost every choice of $a_N>\sum_{k=1}^{N-1}a_k$.
We will now compute the order of magnitude of deviations.


We denote by $W_1$ the Perron Frobenius eigenvector for $\prescript{t}{}{M}$ (normalized so that its coordinates sum up to 1) and by $W_2$ its (normalized) eigenvector associated to the second eigenvalue $\lambda_2$. 
Note that $W_2$ is a real vector since $M$ and $\lambda$ are real.
We decompose $\RR^{N-1}=\text{Vect}(W_1)\oplus \text{Vect}(W_2)\oplus E_3$ where $E_3$ is the direct sum of characteristic spaces of $\prescript{t}{}{M}$ (associated to eigenvalues $\lambda_3, \dots, \lambda_d$) 
By assumption: 
$$\forall v\in E_3, \, \frac{|| \prescript{t}{}{M} v ||}{||v||} \leqslant |\lambda_3| < |\lambda_2| .
$$


We consider cases where $a_N>\sum_{k=1}^{N-1}a_k$ is such that $H \notin E_3$. According to Proposition \ref{not_in_fixed_plane_N}, this is all cases but finitely many ones.

\subsection{Upper bound}

If $w=w_0w_1...w_{n-1}$ is a finite word over the alphabet $\{1,...,N-1\}$, we set $e_w = \sum_{k=0}^{n-1} e_{w_k}$.

Let $x\in \mathbb{S}^1$ and let $\overline{x}_n = x_0x_1...x_{n-1}$ be the prefix of length $n$ of its symbolic coding $\overline{x}$ 
under $T$, decomposed as (see Proposition~\ref{decomp_SD})
$$\overline{x}_n=s_0 \sigma(s_1) ... \sigma^{l-1}(s_{l-1}) \sigma^l( m_l ) \sigma^{l-1} ( p_{l-1}) ... \sigma(p_1) p_0.$$
One has:
$$
e_{\overline{x}_n} = \sum_{k=0}^{l-1} M^k e_{s_k} + M^l e_{m_l} + \sum_{k=0}^{l-1} M^k e_{p_k}.
$$
%
%
%

The deviations from the mean displacement (normalized to be horizontal at unit-speed) are approximated by:
\begin{align*}
S_n h (x) &= \sum_{k=0}^{n-1} h(T^k(x))  \\
		  &= <H, e_{\overline{x}_n}> \\
       &=  \sum_{k=0}^{l-1} < \, \prescript{t}{}{M}^k H, e_{s_k} + e_{p_k} > +  < \, \prescript{t}{}{M}^l H, e_{m_l}>
\end{align*}
The vector $H = Im(G)$ is in the vector space $\text{Vect}(W_2)\oplus E_3$  because $G$ is normal to the Perron-Frobenius vector $V$ of $M$, and $H\notin E_3$.
. 
It can therefore be decomposed as $G_I = \alpha W_2 +U$ with $U \in E_3 $ 
and $\alpha\in\RR^*$. Then one has:

\begin{align*}
S_n h (x) &=  \alpha \left( \sum_{k=0}^{l-1} < \, \prescript{t}{}{M}^k W_2, e_{s_k} + e_{p_k} > +  < \, \prescript{t}{}{M}^l W_2, e_{m_l}> \right) \\
			& \quad + \sum_{k=0}^{l-1} < \, \prescript{t}{}{M}^k U, e_{s_k} + e_{p_k} > +  < \, \prescript{t}{}{M}^l U, e_{m_l}> .
\end{align*}
We write $P_1 = \alpha \left( \sum_{k=0}^{l-1} < \, \prescript{t}{}{M}^k W_2, e_{s_k} + e_{p_k} > +  < \, \prescript{t}{}{M}^l W_2, e_{m_l}> \right)$ and $ P_2 = \sum_{k=0}^{l-1} < \, \prescript{t}{}{M}^k U, e_{s_k} + e_{p_k} > +  < \, \prescript{t}{}{M}^l U, e_{m_l}> $.

The principal part of the deviations is given by:
$$
P_1 = \alpha \left( \sum_{k=0}^{l-1} \lambda_2^k  <W_2, e_{s_k}+e_{p_k}> + z \lambda_2^l  <W_2, e_{m_l}> \right)
$$
%
whose norm is dominated by: 
$$
|\alpha| |\lambda_2|^l\frac{|\lambda_2|}{|\lambda_2| -1} 2(K-1)||W_2||_\infty
$$
where $K$ is such that all $s_k$ and $p_k$ are at length at most $K-1$ and $m_l$ at most $2K-2$ (see Proposition~\ref{decomp_SD}).
The other part of the deviations, $P_2$, whose norm is dominated by $ |\lambda_3|^l\frac{|\lambda_3|}{|\lambda_3| -1} 2(K-1)||U||_\infty $, has a smaller order.
To summarize, we have, for $l$ big enough:
$$
|S_n g (x)|\leqslant D\lambda_2^l
$$ 
with
$D =  \alpha\frac{\lambda_2}{\lambda_2 -1} 2(K-1)||W_2||_\infty + 1.$

Finally, using Corollary~\ref{cor_nl}, one gets, for $n$ big enough:
$$|S_n g (x)| \leqslant D' n^{\frac{\log(\lambda_2)}{\log(\lambda_1)}} $$
with $D' = D \exp(\frac{\log(\lambda_2)\log(A)}{\log(\lambda_1)})$.
One can replace $D'$ by a larger constant $C_1$ so that the inequality 
$$|S_n g (x)|\leqslant C_1 n^{\frac{\log(\lambda_2)}{\log(\lambda_1)}} $$ 
is true for all $n\in\NN$.

This shows the first inequality of Theorem~\ref{th_selfsim}. We prove the second one in the following paragraph.

\subsection{Lower bound} 
We want to build an increasing sequence $(n_l)_{l\in\NN}$ of indices of "significant deviations", in other words such that
$$
(\star) \qquad \exists C_2, \forall l \geqslant l_0, |S_{n_l}h(x)|>C_2 n_l^{\frac{\log\lambda_2}{\log\lambda_1}}.
$$

Let $l$ be an integer. We decompose the symbolic coding of $x$ as 
$$\overline{x} = s_0\sigma(s_1)...\sigma^{l-1}(s_{l-1})\sigma^l(w)$$
 where the $s_i$ are finite words of length at most $K-1$ (see Proposition~\ref{decomp_SD}) and $w$ is an infinite word. This corresponds to decompose the orbit of $x$ with Rokhlin towers of order $l$. 

Let us denote by $\pi_j$ (for $1\leqslant j \leqslant 2)$ the projection onto Vect$(W_j)$ and $\pi_3$ the one onto $E_3$ with respect to the decomposition $\RR^{N-1} = \text{Vect}(W_1)\oplus\text{Vect}(W_2)\oplus E_3$.

Let us denote $m_l = |s_0\sigma(s_1)...\sigma^{l-1}(s_{l-1})|$, corresponding to the hight of the partial Rokhlin tower before the orbit of $x_0$ arrives at the bottom of a tower of order $l$.

If $|S_{m_l} h(x)| > |\lambda_2|^{l}$, then we define $n_l = m_l$ 
and we have:
$$
|S_{m_l} h (x)| >  |\lambda_2|^{l} \geqslant D |\lambda_2|^{\frac{\log(n_l)}{\log(\lambda_1)}} = D n_l^{\frac{\log(|\lambda_2|)}{\log(\lambda_1)}}
$$
with $D=exp(\frac{-C}{\log(\lambda_1)}) $ (where $C$ is the constant given by Corollary~\ref{cor_nl}).

Else ($|S_{m_l} h (x)| \leqslant |\lambda_2|^{l}$). At least one vector $e_k$ of the canonical basis of $\RR^{N-1}$ is such that $\pi_2(e_k)\neq 0$. 
Let $i$ be such that $
\pi_{2}(e_{w_i}) \neq 0 $ and for every index $j<i$, $\pi_{2}(e_{w_j}) = 0$. The matrix $M$ is primitive, so is the substitution $\sigma$. Hence there exists $L$ (which depends only on $M$, not on $w$) such that $i\leqslant L$. 
Let us denote $m'_l = |\sigma^l(w_1 ... w_{i-1})|$ and $m_l ''= |\sigma^l(w_i)|$.
We define $n_l$ as the length $  |s_0\sigma(s_1)...\sigma^{l-1}(s_{l-1})\sigma^l(w_1...w_i)| = m_l + m'_l + m_l ''$.
Then we have:
\begin{align*}
|S_{n_l} h(x) | &= |S_{m_l} h(x) + S_{m'_l}h (T^{m_l}(x)) + S_{m_l''} h (T^{m_l+m'_l} (x))| \\
&\geqslant |S_{m_l''} h (T^{m_l+m'_l} (x) )| - |S_{m_l} h(x)| - | S_{m'_l}h (T^{m_l}(x)) |   
\end{align*}
By hypothese, one has
$$|S_{m_l} h(x)| \leqslant |\lambda_2|^{l}  . $$
By definition of $i$, if we set $c=\underset{1\leqslant k \leqslant N-1}{\max} | \pi_{3}(e_k) | $, one has
 $$| S_{m'_l}h (T^{m_l}(x)) | \leqslant c |\lambda_{3}|^i \leqslant c |\lambda_{3}|^L  .$$ 
Finally
\begin{align*}
\left|S_{m_l''} h \left(T^{m_l+m'_l} (x) \right)\right| &= \left|< \, \prescript{t}{}{M}^l H, e_{w_i} > \right| \\
			&\geqslant \left| \alpha \right| \left|\lambda_2^l < W_2, e_{w_i} >\right| - \left|< \, \prescript{t}{}{M}^l U, e_{w_i} > \right| \\
			&\geqslant \left| \alpha \right| |\lambda_2|^l - c|\lambda_{3}|^l ||U||.
\end{align*}
Putting these inequalities together 
we conclude that, for $l$ big enough:
$$
|S_{n_l} h(x) | \geqslant \frac{ \left| \alpha \right|}{2} |\lambda_2|^l \geqslant \frac{ \left| \alpha \right|}{2}D |\lambda_2|^{\frac{\log(n_l)}{\log(\lambda_1)}} = \frac{ \left| \alpha \right|}{2}D n_l^{\frac{\log(|\lambda_2|)}{\log(\lambda_1)}}.
$$
%


In the first case $|S_{m_l} h (x)|>|\lambda|^l$  so $m_l$ goes to infinity as $l$ increases and 
in the second case $n_l \geqslant m_l'' \geqslant A\lambda_1^l$ (Lemma~\ref{lem_encadrement}). Therefore the sequence $(n_l)_{l\in\NN}$ goes to infinity.
If we set $$C_2 = \frac{ \left| \alpha \right|}{2}D= \frac{ \left| \alpha \right|}{4} exp(\frac{-C}{\log(\lambda_1)}),$$ the inequality $(\star)$ holds for $l\in\NN$ big enough, which gives the second inequality of Theorem~\ref{th_selfsim}.


\section{Concluding remarks}
In this article we assumed that the $N$-gons are inscribed in a half-circle (that means a proportion of $\frac{N}{2^N}$ of the cyclic $N$-gons) and we studied trajectories passing near the circumcenter of the polygons. In this section, we discuss briefly the obstacles to remove such assumptions.
We recall the key elements to prove Theorems \ref{th_dev}, \ref{th_generic} and \ref{th_selfsim}.
\begin{enumerate}
\item
We can study the trajectory via an IETF $\Phi$ because the $N$-gon is cycic.
\item\label{item_square}
The square $\Phi^2$ of this map stabilizes a subinterval $J$, and $\Phi^2_{|J}$ is an IET with $N-1$ intervals of continuity. We do not need assumption to ensure that $Phi^2$ does not flip any interval, but we do need the assumption on the parameter $\tau$ of the trajectory to ensure that  $\Phi^2_{|J}$ has only $N-1$ intervals of continuity.
\item
We can define a function $f$ of displacement, that is analytic.
\item
We can show that for any fixed codimension 3 (or 2 if $N$ is odd) subspace $E$ of $\RR^{N-1}$, the vector $H$, whose coordinates are the values of the function $h = Im(\frac{1}{m}f - 1)$, is not in $E$. This relies strongly on the fact that we have an explicit expression of $H$.
\end{enumerate}

We could overcome point \ref{item_square} by studying the Rauzy-Nogueira diagram of the underlying permutation of $\Phi$. We would get many open sets of parameters, each one associated to a new map, that is a renormalization of $\Phi$.
So far, we can again define a function $f$ of displacement and even ensure that it is analytic.
But if we want to conclude, we need to have an explicit expression of the function $f$ for each open set of parameters. These open sets depend on paths in the Rauzy-Nogueira diagram in question.
The trouble is that this diagram, for $N=5$, has more than 1 millon of vertices and more than 3 millions of edges! Hence establishing the apropriate open sets seems out of reach without another idea, either about how we cut the set of parameters, or about how we show that $H$ is not in $E$.


\printbibliography 




%
%

\vspace{0.5cm}

Magali Jay, \textsc{Aix-Marseille Université, I2M, Marseille, France}

E-mail address: \url{magali.jay@univ-amu.fr}

\end{document}